\documentclass[psamsfonts]{amsart}
\usepackage{amssymb}

\input{diagrams}
\diagramstyle[scriptlabels]

\newcommand\query[1]{}

\def\ns{\noalign{\bigskip}}

\def\overrightarrow#1{\vbox{\ialign{##\crcr
      \rightarrowfill\crcr\noalign{\kern-\p@\nointerlineskip}
      $\hfil\displaystyle{#1}\hfil$\crcr}}}
\def\overleftarrow#1{\vbox{\ialign{##\crcr
      \leftarrowfill\crcr\noalign{\kern-\p@\nointerlineskip}
      $\hfil\displaystyle{#1}\hfil$\crcr}}}
\def\lrarrow{\relbar\joinrel\longrightarrow}
\def\llrarrow{\relbar\joinrel\relbar\joinrel\longrightarrow}

\def\dh{\begin{array}{l}\kern -.10em \hbox{$\scriptstyle{\cap\;}$}\\  [-5pt]  
      \hbox{$\Biggl\downarrow$}\end{array}}

\def\longhookrightarrow{\mathop{\raisebox{.6ex}{$\scriptscriptstyle\subset$}\kern -.6em        
      \relbar\joinrel\relbar\joinrel\relbar\joinrel\relbar\joinrel\joinrel\longrightarrow}}
\def\lhookrightarrow{\mathop{\raisebox{.6ex}{$\scriptscriptstyle\subset$}\kern -.6em        
      \relbar\joinrel\relbar\joinrel\joinrel\longrightarrow}}
\def\lhookleftarrow{\mathop{\longleftarrow\joinrel\relbar\joinrel\relbar\joinrel      \kern -.6em\raisebox{.6ex}{$\scriptstyle\supset$}}}

\def\hookupleftarrow{\begin{array}{l}\Biggl\uparrow\\ [-8pt]
  \kern -.10em\scriptstyle{\cup}\end{array}}
\def\hookuprightarrow{\begin{array}{l}\kern -.20em\Biggl\uparrow\\ [-8pt]  
   \scriptstyle{\cup}\end{array}}

\newarrow{Equal}{}{=}{}{=}{}

\def\lto{{\longrightarrow}}
\def\hto{{\hookrightarrow}}
\def\xto{\xrightarrow}
\def\lllrarrow{\relbar\joinrel\relbar\joinrel\relbar\joinrel%
\relbar\joinrel\longrightarrow}

\newcommand{\fa}{{\mathfrak a}}

\newcommand{\fm}{{\mathfrak m}}

\newcommand{\cala}{{\mathcal A}}
\newcommand{\calb}{{\mathcal B}}
\newcommand{\calc}{{\mathcal C}}

\newcommand{\cale}{{\mathcal E}}
\newcommand{\calf}{{\mathcal F}}
\newcommand{\calg}{{\mathcal G}}
\newcommand{\calh}{{\mathcal H}}
\newcommand{\cali}{{\mathcal I}}
\newcommand{\calj}{{\mathcal J}}
\newcommand{\calk}{{\mathcal K}}

\newcommand{\calm}{{\mathcal M}}
\newcommand{\caln}{{\mathcal N}}
\newcommand{\calo}{{\mathcal O}}
\newcommand{\calp}{{\mathcal P}}
\newcommand{\calq}{{\mathcal Q}}
\newcommand{\calr}{{\mathcal R}}
\newcommand{\cals}{{\mathcal S}}

\newcommand{\calv}{{\mathcal V}}
\newcommand{\calw}{{\mathcal W}}

\newcommand{\calos}{{\calo_S}}
\newcommand{\calot}{{\calo_T}}

\newcommand{\calow}{{\calo_W}}
\newcommand{\calox}{{\calo_X}}
\newcommand{\caloy}{{\calo_Y}}
\newcommand{\caloz}{{\calo_Z}}

\def\cExt{{\mathcal E}xt}
\def\cTor{{\mathcal T}or}
\def\cHom{{\mathcal H}om}

\newcommand{\bbbl}{{\mathbb L}}
\newcommand{\bbbq}{{\mathbb Q}}
\newcommand{\bbbc}{{\mathbb C}}

\newcommand{\bbbs}{{\mathbb S}}
\newcommand{\bbbn}{{\mathbb N}}
\newcommand{\bbbp}{{\mathbb P}}

\newcommand{\bbbz}{{\mathbb Z}}

\newcommand{\bB}{\ensuremath{\mathbf B}}

\newcommand{\bE}{\ensuremath{\mathbf E}}

\newcommand{\vp}{\varphi}

\DeclareMathOperator{\At}{At}

\DeclareMathOperator{\ch}{ch}

\DeclareMathOperator{\coker}{coker}

\DeclareMathOperator{\Con}{Con}

\DeclareMathOperator{\Der}{Der}

\DeclareMathOperator{\dotimes}{\underline{\otimes}}

\DeclareMathOperator{\End}{ End}
\DeclareMathOperator{\Ex}{ Ex}

\DeclareMathOperator{\Ext}{Ext}

\DeclareMathOperator{\HH}{H}
\DeclareMathOperator{\Hilb}{Hilb}
\DeclareMathOperator{\Hom}{Hom}
\DeclareMathOperator{\calhom}{{\mathcal Hom}}
\DeclareMathOperator{\id}{id}
\DeclareMathOperator{\Imm}{Im}

\DeclareMathOperator{\Ker}{Ker}

\DeclareMathOperator{\ob}{ob}
\DeclareMathOperator{\Ob}{Ob}

\DeclareMathOperator*{\Projlim}{\underleftarrow{\,\vphantom{(}\lim}}
\DeclareMathOperator{\Quot}{Quot}

\DeclareMathOperator{\RHom}{RHom}

\DeclareMathOperator{\Sym}{{\bbbs}}
\DeclareMathOperator{\Tor}{Tor}

\DeclareMathOperator{\Tr}{Tr}

\newcommand{\sbullet}{{\scriptstyle\bullet}}

\DeclareMathOperator{\Ans}{{\mathbf{\An}}_\Sigma}
\DeclareMathOperator{\An}{{\mathbf{An}}}
\DeclareMathOperator{\Ab}{{\mathbf{Ab}}}
\DeclareMathOperator{\Coh}{\mathbf{Coh}}

\DeclareMathOperator{\Sets}{\ensuremath{ \mathbf{Sets}}}

\newcommand{\catc}{\ensuremath{ \mathbf C}}
\newcommand{\cate}{\ensuremath{ \mathbf E}}

\def\Arts{\mathop{{\bf Art}_\Sigma}\nolimits}

\def\Formals{\mathop{{\bf An\hat{}}_\Sigma}\nolimits}

\def\Cech{\v{C}ech}

\theoremstyle{definition}
\newtheorem{defn}{Definition}[section]
\newtheorem{exam}[defn]{Example}
\newtheorem{rem}[defn]{Remark}
\newtheorem{rems}[defn]{Remarks}
\newtheorem{sit}[defn]{}

\theoremstyle{plain}
\newtheorem{prop}[defn]{Proposition}
\newtheorem{theorem}[defn]{Theorem}
\newtheorem{lem}[defn]{Lemma}
\newtheorem{cor}[defn]{Corollary}

\theoremstyle{remark}
\newtheorem{prob}[defn]{Problem}

\begin{document}

\title[Semiregularity Map and Deformations]{A Semiregularity Map for
Modules and Applications to Deformations}

\author{Ragnar-Olaf Buchweitz}
\address{Dept.\ of Math., University of
Tor\-onto, Tor\-onto, Ont.\ M5S 3G3, Canada}
\email{ragnar@math.utoronto.ca}

\author{Hubert Flenner}
\address{Fakult\"at f\"ur Mathematik der Ruhr-Universit\"at,
Universit\"atsstr.\ 150, Geb.\ NA 2/72, 44780 Bochum, Germany}
\email{Hubert.Flenner@ruhr-uni-bochum.de}

\subjclass{14C25, 32G05, 13D03}
\thanks{R.O.B.~was partly supported by NSERC grant
3-642-114-80.  
\endgraf 
Both authors were partly supported by a grant from the Volkswagen
Foundation under the Research in Pairs program at
Math.~Forschungsinstitut Oberwolfach.}  

\date{\today}

\begin{abstract} 
We construct a general semiregularity map for algebraic cycles as
asked for by S.\ Bloch \cite{Blo} in 1972. The existence of such a
semiregularity map has well known consequences for the structure of
the Hilbert scheme and for the variational Hodge conjecture. Aside
from generalizing and extending considerably previously known results
in this direction, we give new applications to deformations of modules
that encompass, for example, results of Artamkin \cite{Art} and Mukai
\cite{Muk}.

The formation of the semiregularity map here involves powers of the
cotangent complex, Atiyah classes, and trace maps, and is defined not
only for subspaces of manifolds but for perfect complexes on arbitrary
complex spaces. It generalizes in particular Illusie's \cite{Ill}
treatment of the Chern character to the analytic context and
specializes to Bloch's earlier description of the semiregularity map
for locally complete intersections as well as to the infinitesimal
Abel-Jacobi map for submanifolds.
\end{abstract}

\maketitle

\tableofcontents

\section{Introduction} 

Let $Z$ be a closed complex subspace of a compact complex manifold $X$
and let $[Z]$ denote the corresponding point in the Douady space $H_X$
of $X$, the complex analytic analogue of the Hilbert scheme. It is a
classical fact that the tangent space of $H_X$ at $[Z]$ is naturally
isomorphic to $H^0(Z,\caln_{Z/X})$, where
$\caln_{Z/X}=\cHom_X(\calj,\caloz)$ denotes the normal sheaf of $Z$ in
$X$ with $\calj\subseteq \calox$ the ideal sheaf of $Z$. Moreover, if
$Z$ is locally a complete intersection in $X$ then the vanishing of
$H^1(Z,\caln_{Z/X})$ implies that $[Z]$ is a smooth point of $H_X$. It
was, however, already observed by Severi \cite{Sev} that the converse
is not true in general. He introduced the notion of a semiregular
curve on a surface to mean that the restriction map
$H^0(X,\omega_X)\to H^0(Z,\omega_X|Z)$ is surjective or, dually, that
the {\em semiregularity map\/} $H^1(Z,\caln_{Z/X})\to H^2(X,\calox)$
is injective, and showed that the point of the Hilbert scheme
corresponding to a semiregular curve is always smooth. This result was
extended to divisors in arbitrary projective complex manifolds by
Kodaira and Spencer \cite{KSp}. In 1972, S.~Bloch \cite{Blo} was able
to define more generally for every locally complete intersection $Z$
of codimension $q$ in $X$ a semiregularity map
$\tau:H^1(Z,\caln_{Z/X})\to H^{q+1}(X,\Omega^{q-1}_X)$ to show again
that the injectivity of $\tau$ implies that $[Z]$ is a smooth point of
$H_X$. His semiregularity map admits a simple description using Serre
duality. However, in case of an arbitrary subspace $Z$ of $X$ the
obstructions for extending embedded deformations lie in the tangent
cohomology group $T^2_{Z/X}(\caloz)$, and it is no longer possible to
apply duality. Thus the problem arises to define such a semiregularity
map by other means.

In our approach we will more generally assign first a semiregularity
map $\sigma:\Ext^2_X(\calf,\calf)\to \prod_{q\ge
0}H^{q+2}(X,\Omega_X^q)$ to every coherent $\calox$--module $\calf$ on
a complex manifold. Indeed this map will be defined for any coherent
$\calox$--module $\calf$ of locally finite projective dimension, or
even for perfect complexes of modules, on arbitrary complex spaces,
and it occurs as the component $\sigma = \sigma^{(2)}$ of a family of
maps
$$
\sigma^{(r)}:\Ext^r_X(\calf,\calf)\to \prod_{q\ge
0}H^{q+r}(X,\Lambda^q\bbbl_X)\,,\quad {r\ge 0}\,,
$$
where $\bbbl_X$ denotes the cotangent complex of $X$. We will
outline in brief our construction when $X$ is smooth, which special
case is also subject of our survey \cite{BFl2}.

To begin with, we assign to $\calf$ its Atiyah class, as originally
defined in \cite{At} for locally free sheaves. Following \cite{Ill}, a
possible way of construction for any coherent $\calox$--module is to
use the extension on $X\times X$ that defines the module of analytic
differential forms,
\begin{equation}
\label{eqn:1}
0\to \calj/\calj^2\cong\Omega^1_X\to\calo_{X\times X}/\calj^2\to
\calox\to 0\,,
\end{equation}
where $\calj\subseteq \calo_{X\times X}$ is the ideal of the diagonal.
With $p_i:X\times X\to X$ the $i^{th}$ projection for $i=1,2$, we
tensor (\ref{eqn:1}) with $p_1^*(\calf)$ and consider the resulting
sequence
$$
0\to\calf \otimes_\calox\Omega^1_X  \to
p_1^*(\calf)
\otimes_{\calo_{X\times X}} \calo_{X\times X}/\calj^2
\to \calf\to 0
$$
as an extension of $\calox$--modules via $p_{2*}$ so that it defines
an element
$$
\At(\calf)\in\Ext^1_X(\calf, \calf\otimes_\calox \Omega^1_X)\,,
$$
 the {\em Atiyah class} of $\calf$. Taking
powers gives elements $\At^q(\calf)\in \Ext^q_X(\calf,
\calf\otimes\Omega^q_X)$. Now the $q^{th}$ component of the
semiregularity map
$\sigma$ is the composition of the two maps
\begin{diagram}
\sigma_q:
\Ext^2_X(\calf,\calf)&
\rTo^{ *\cdot  (-\At(\calf))^q/q!}&
\Ext_X^{q+2}(\calf,\calf\otimes\Omega_X^q)&
\rTo^{\Tr}&
H^{q+2}(X,\Omega_X^q)\,,
\end{diagram}
where $\Tr$ is the trace map as defined in \cite{Ill,OTT}. Finally, to
get a semiregularity map for subspaces $Z\subseteq X$, we observe that
there is a natural homomorphism from $T^k_{Z/X}(\caloz)$ into
$\Ext^k_{X}(\caloz,\caloz)$ for each $k\ge 0$. Composing this for
$k=2$ with the map $\sigma$ above in case $\calf=\caloz$ gives the
desired semiregularity map
$$
\tau = (\tau_q)_{q\ge 0}: T^2_{Z/X}(\caloz)\to \prod_{q\ge 0}
H^{q+2}(X,\Omega_X^q)
$$
for subspaces. We will verify in Section 8 that for a locally complete
intersection $Z$ in $X$ of codimension $q$ the component $\tau_{q-1}$
of our semiregularity map coincides with the classical one defined by
Bloch.

To understand some of the geometrical implications of such
semiregularity maps, let us restrict further to the case of coherent
modules on a compact algebraic manifold, the case of subspaces being
similar. With respect to the Hodge decomposition,
$H^k(X,\bbbc)=\bigoplus_{p+q=k}H^q(X,\Omega^p_X)$, the Chern character
of a coherent sheaf is obtained from its Atiyah class by the formula
$$
\ch(\calf)=\Tr \exp(-\At(\calf))=\sum_{k\ge 0} 
\frac{(-1)^{k}}{k!}\Tr (\At^k(\calf))\, ,
$$
see \cite{At} for the case of vector bundles and \cite{Ill, OTT} for
the general case.  

Assume now given an infinitesimal deformation of $X$ represented by a
class $\xi\in H^1(X,\Theta_X)$, where $\Theta_X$ denotes the tangent
bundle. By Bloch's interpretation of Griffiths' transversality
theorem, for fixed $k\ge 0$ the unique horizontal extension of the
cohomology class $\ch_{k+1}(\calf)$ relative to the Gau\ss-Manin
connection stays of Hodge type $(k+1, k+1)$ if and only if the
contraction of this class by $\xi$ vanishes, $\langle \xi ,
\ch_{k+1}(\calf)\rangle = 0$ in $H^{k+2}(X, \Omega^k_{X})$. On the
other hand, let us consider the deformations of $\calf$ itself instead
of just extending its Chern character horizontally. The deformations
of $\calf$ are controlled by the differential graded Lie algebra
underlying $\Ext_{X}^\sbullet(\calf,\calf)$ so that the space of
infinitesimal deformations is given by $\Ext_{X}^1(\calf,\calf)$ and
the obstructions to extend deformations live in
$\Ext_{X}^2(\calf,\calf)$. Contracting against the negative of the
Atiyah class serves as an obstruction map
$$
\ob :=\langle *, -\At(\calf)\rangle: H^1(X,\Theta_{X})\lto
\Ext_{X}^2(\calf,\calf)
$$
so that $\calf$ admits a deformation into the direction of $\xi$ if
and only if $\ob(\xi)=0$; see \cite{Ill} for the algebraic and
\ref{3.4} for the analytic case. A key observation is now that the
maps just described fit into a commutative diagram
\begin{diagram}[w=7mm,h=8mm,midshaft]
H^{1}(X,\Theta_X)
&& \rTo^{\langle *, -\At(\calf)\rangle} && \Ext_{X}^2(\calf,\calf)\\
&\rdTo<{\langle *, \ch_{k+1}(\calf)\rangle}&&
\ldTo>{\sigma_{k}}\\
&& H^{k+2}(X,\Omega^k_X)\,.
\end{diagram}
As a consequence of (a suitable generalization) of this fact and and
generalizing the arguments of \cite{Blo} we obtain for instance in
Section 5 that {\em the variational Hodge conjecture holds for cycles
that are representable as $(k+1)^{st}$ component of the Chern
character of a $k$-semiregular sheaf} $\calf$, where we mean by {\em
$k$-semiregular\/} that the component $\sigma_{k}$ of the
semiregularity map for $\calf$ is injective. An analogous result holds
for subspaces $Z\subseteq X$ that are ($k$-)semiregular in the
corresponding sense.

Other important applications are to deformations of modules. In
analogy with the aforementioned results of Bloch we will show that
{\em the basis of the semiuniversal deformation of $\calf$ is smooth
if the semiregularity map $\sigma$ is injective}. We will deduce this
result more generally for arbitrary singular complex spaces $X$ and
also for relative situations. We derive as well analogous smoothness
criteria for the Douady space and the $\Quot$-space, and give
applications to deformations of holomorphic mappings.

Ideally, the semiregularity map, say, for a module $\calf$ should
correspond to a morphism between two deformation theories so that it
maps the obstruction space $\Ext^2_{X}(\calf, \calf)$ into the
obstruction space of some other deformation theory. It seems quite
clear that this second deformation theory should be given in terms of
the intermediate Jacobians, or, more naturally, by Deligne cohomology.
The map from the deformations of $\calf$ into the intermediate
Jacobian, say $J^k(X)$, should be a generalized Abel-Jacobi map that
associates to a deformation of $\calf$ over a germ $(S,0)$ the map of
$S$ into the intermediate Jacobian given by integration over a
topological cycle whose boundary is the difference of the $k^{th}$
Chern characters of the special fibre and the fibre over $s$. As the
intermediate Jacobian is smooth this would provide a satisfactory
explanation of the fact that the injectivity of $\sigma$ implies the
smoothness of the versal deformation of $\calf$, and it would show
that all obstructions of $\calf$ vanish under $\sigma$ and not merely
the curvilinear ones as we show here; for the special case of
submanifolds instead of modules see \cite{Cle, Ran4}. Such an
interpretation indeed applies for the lowest component of the
semiregularity map: the work of Artamkin \cite{Art} and Mukai
\cite{Muk} interprets $\sigma_0 : \Ext^2_{X}(\calf,\calf) \to
H^{2}(X,\calox)$ as the map between obstruction spaces for the
deformations of $\calf$ versus those of its determinant line bundle.
As a further clue that such interpretation might be true in general we
verify in Section 9 that for a submanifold $Z$ of $X$ the differential
of the Abel-Jacobi map admits the same homological description as the
semiregularity map. \medskip

A few remarks about the contents of the various sections: In Section 2
we review the technique of Forster-Knorr systems, originally used in
\cite{FKn} and further exploited in \cite{Pal,Fle1}, to construct a
cotangent complex for complex spaces. Most of that material is a
largely generalized version of parts of \cite{Fle1}. As this source is
not easily accessible we use the opportunity to give an exposition of
that technique of simplicial spaces of Stein compact sets in the
generality we need. Following \cite{Pal} we will review in brief the
notion of resolvents and give the relevant descriptions of tangent
(co)homology as used later on.

In Section 3 we construct the Atiyah class of a coherent sheaf $\calf$
as a class in $\Ext^1_{X}(\calf, \calf\dotimes\bbbl_X)$, thereby
generalizing the construction by Illusie \cite{Ill} to the analytic
case. First we do this for modules on simplicial spaces of Stein
compact sets and then use the results of the previous section to
descend to actual complex spaces. Following the classical approach of
Atiyah \cite{At}, see also \cite{ALJ}, we will construct these classes
using connections, in our case on modules over the resolvent of a
complex space, thus verifying the basic funtorial properties by
explicit computation.

Section 4 contains the construction of the semiregularity map for
modules as well as for subspaces. We give an interpretation of the
obstruction map for modules or subspaces in terms of Atiyah classes
and derive the aforementioned commutative diagram. This is the basic
tool in Section 5, where we prove the variational Hodge conjecture for
the special case described above.

In Section 6 we prove some general criteria for the smoothness of the
basis of a semiuniversal deformation. We give a new and transparent
proof of the $T^1$-lifting criterion of Ran and Kawamata
\cite{Kaw1,Kaw2,Ran2,Ran3} and show how to deduce their results by a
simple argument from the well known Jacobian criterion for smoothness.
Our method also yields new generalizations to the relative case. These
results, together with the existence of the semiregularity map, have
many applications in deformation theory, some of which we treat in
Section 7. In the first part we deduce applications to deformations of
modules as mentioned above. In the second part we turn to the Douady
space and give various criteria for its smoothness. The more general
case of the Quot-scheme is treated in the third part, and in the last
part we apply our constructions to deformations of mappings and define
a semiregularity map under very general circumstances. For the special
case of stable curves, results in this direction were independently
obtained by K.~Behrend and B.~Fantechi.

In the final Section 8 we compare our semiregularity map with the one
constructed by S.~Bloch. This requires an explicit description of the
trace map via a Cousin type resolution. Moreover, we show how the
infinitesimal Abel-Jacobi map fits into this framework.

In an Appendix we collect some results on integral dependence and
infinitesimal deformations of complex spaces that are needed in
Section 6. Especially the (elementary) characterizations of the
subspaces of $T^1_X$ given by the curvilinear extensions, respectively
by $\Ext_{X}^1(\Omega_X^1,\calox)$, seem to be new.

\subsection*{General notation}
We explain some notation used throughout this paper.

Categories are written in boldface and categories like $\Sets$
should need no further explanation. Whenever we talk about
isomorphism classes of objects from a category $\catc$, it will be
assumed that those classes form a set. Such a category $\catc$ is
sometimes called {\em essentially small\/}.

A germ of a (formal) complex space is denoted $(S,0)$ or often simply
$S$. As a rule, every germ has $0$ as its basepoint, and the same
symbol represents the (reduced) point. For a (formal) complex space
$X$ the sheaf of holomorphic functions is as usual denoted $\calox$,
whereas for a germ $S=(S,0)$ the symbol $\calos$ indicates the local
ring $\calo_{S,0}$ and $\fm_S$ its maximal ideal.

If $X$ is a complex space then $\Coh(X)$ will be its category of
coherent modules. Similarly, if $S=(S,0)$ is the germ of a (formal)
complex space, $\Coh(S)$, $\Coh_{art}(S)$ will denote the categories
of finite, respectively finite artinian $\calo_{S,0}$--modules. A
closed embedding $S\hookrightarrow S'$ of complex spaces is an {\em
extension\/} of $S$ by $\calm\in\Coh(S)$ if the ideal
$\cali:=\Ker(\calo_{S'}\to\calo_S)$ defining $S$ in $S'$ is of square
zero and isomorphic to $\calm$ as $\calo_S$--module under a fixed
isomorphism. In particular, $S[\calm]$ indicates the {\em trivial
extension\/} whose structure sheaf $\calo_{S[\calm]}$ is the direct
sum $\calo_S\oplus \calm\varepsilon $ endowed with the multiplication
$(s+m\varepsilon)(s'+m'\varepsilon)=ss'+(sm'+ms')\varepsilon$ so that
$\varepsilon^2=0$.

To reduce complexity of display, we use unadorned tensor products,
such as $\calm\otimes \caln$, whenever the sheaf or ring over which
the tensor product is taken should be clear from the context. We also
use $\dotimes$ instead of $\otimes^{\bbbl}$ to denote derived tensor
products.

\section{Homological algebra on simplicial schemes of Stein
compact sets}

Let $X$ be a complex space and $\calf$ a coherent $\calox$--module. In
the introduction we reviewed for $X$ smooth the construction of the
{\em Atiyah class\/}, $\At(\calf)\in\Ext^1_X(\calf,
\calf\otimes_\calox \Omega^1_X)$, of $\calf$ that is given by the
extension
$$
0\to\calf \otimes_\calox\Omega^1_X  \to
p_1^*(\calf)
\otimes_{\calo_{X\times X}} \calo_{X\times X}/\calj^2
\to \calf\to 0\,.
$$
In the singular case, however, we need to produce more generally the
Atiyah class in
$$
\Ext^1_X(\calf, \calf\dotimes_\calox\bbbl_X)\,,
$$
where $\bbbl_X$ is the cotangent complex of $X$ and $\dotimes$ 
denotes the derived tensor product.

In the algebraic setting this was done by Illusie \cite{Ill} using
simplicial methods. Those do not immediately generalize to the
analytic situation due to the lack of global resolutions. Instead we
use the technique of Forster-Knorr systems \cite{FKn} on simplicial
schemes of Stein compact sets and construct resolvents on complex
spaces as in \cite{Fle1, Pal}. Alternatively one might be tempted to
use the method of twisted cochains as developed in \cite{OTT}, but
this theory is so far only available for manifolds, in particular the
theory of cotangent complexes on singular spaces has not yet been
established in that framework.

We first state and (indicate how to) prove the results on the
homological algebra of Forster-Knorr systems that we will use. Key
references are \cite{Fle1}, \cite{Pal}; see also \cite{BMi}. All our
complex spaces are assumed to be paracompact.

\begin{sit}\label{2.1}
A subset $K$ of a complex space $X$ is called {\em Stein compact}\/ if
it is compact, semianalytic and admits arbitrary small open
neighbourhoods that are Stein. We equip $K$ with the structure sheaf
$\calo_K:=\calox|_{K}$ so that the ring $\Gamma(K,\calo_K)$ consists
of all $K$-germs of functions that are analytic in an open
neighbourhood of $K$ in $X$. In the extreme case that $K=\{x\}$
consists just of a point in $X$ one retrieves the local analytic
algebra $\calox_{,x}$. Every coherent $\calo_K$--module for a Stein
compact set $K$ satisfies Cartan's Theorems A and B, and by a
fundamental result of Frisch \cite{Fri} the ring $\Gamma(K,\calo_K)$
is noetherian. These facts imply that a Stein compact set behaves like
a noetherian affine scheme: for example, the category $\Coh\calo_K$ of
coherent $\calo_K$--modules is equivalent to the category of finite
$\Gamma(K,\calo_K)$--modules, in particular it contains enough
projectives. Again as for affine schemes, these projectives are
usually not projective objects in the category of all
$\calo_K$--modules.

Note that the {\em dimension\/} of a Stein compact set is given by
$$
\dim K := \sup\{\dim \calo_{K,x}\mid x\in K\} = \inf\{\dim U\mid
K\subseteq U \subseteq X, U \text{\ open}\}\,.
$$

If $K\subseteq X$, $L\subseteq Y$ are Stein compact sets with $X$, $Y$
complex spaces over some complex space $Z$, then the product $K\times
L\subseteq X\times Y$ as well as the fibre product $K\times_Z
L\subseteq X\times_Z Y$ are again Stein compact sets.

A covering of a complex space $X$ by Stein compact sets $\{X_i\}_{i
\in I}$ is called {\em locally finite\/} if every point in $X$ admits
an open neighbourhood that meets only finitely many $X_i$. Any two
such coverings clearly admit common refinements and, as $X$ is
paracompact, each covering by open sets can be refined to a locally
finite covering by Stein compact sets.
\end{sit}

\begin{sit}
\label{2.2}
Let $X$ be a complex space over $Y$ and $\{X_i\}_{i \in I}$ a locally
finite covering of $X$ by Stein compact sets. The {\em nerf\/} of the
covering is the simplicial set
$$
A = \{\alpha \subseteq I \mid
X_\alpha:=\textstyle
\bigcap\limits_{i \in\alpha} X_i \ne \emptyset\}\ ,
$$
whose simplices are thus finite subsets of $I$, and
$$
X_* = (X_\alpha)_{\alpha \in A}\,, \text{with the inclusions}\quad
X_\beta
\hto X_\alpha\quad\text{for $\alpha\subseteq\beta$,}
$$
forms the corresponding {\em simplicial scheme of Stein compact
sets\/} over $A$. Denote $|\alpha|=k$ the dimension of the
simplex $\alpha=\{i_0,...,i_k\}$, not to be confused with $\dim
X_{\alpha}$, the complex dimension of the Stein compact $X_{\alpha}$.

Aside from simplicial schemes of Stein compact sets arising from
locally finite coverings we will also need the following type. Let $A$
be again a simplicial set of finite subsets of some index set $I$ and
assume given for any vertex $i\in I$ a Stein compact set $L_i\subseteq
\bbbc^{n_i}\times Y$ for some $n_i$. Set
$$
\textstyle W_\alpha:=\prod^Y_{i\in\alpha}L_i\,,\quad \alpha\in A\,,
$$
where $\prod^Y$ denotes the fibre product over $Y$, and let the
natural projections $p_{\alpha\beta}: W_\beta\to W_\alpha$ for
$\alpha\subseteq\beta$ serve as transition maps. The natural maps from
$W_\alpha$ onto $Y$ are smooth and compatible with the transition
maps. The resulting simplicial scheme $W_*\to Y$ over $Y$ is called
{\em smooth\/}.

Returning to the complex space $X$ with its given covering choose
closed $Y$-embeddings $X_i\hto L_i\subseteq \bbbc^{n_i}\times Y$ with
$L_i$ a Stein compact subset. These data yield diagonal embeddings
$X_\alpha \hto W_\alpha,\alpha\in A,$ that define a morphism of
simplicial schemes $X_*\hto W_*$ over $Y$. We will refer to it as a
(simplicial) {\em smoothing} of the given map $X\to Y$.
\end{sit}

\begin{sit}
Now we consider (negatively graded, simplicial) DG algebras $\calr_*
=\bigoplus_{i\le 0} \calr_*^i$ over $W_*$. The differential of such an
algebra is a derivation of degree $+1$ and will be denoted
$\partial_{\calr_*}$ or simply $\partial$. The simplicial structure
consists of a system of compatible maps
$p_{\alpha\beta}^{-1}(\calr_\alpha) \to \calr_\beta$ for
$\alpha\subseteq\beta$ that are morphisms of DG algebras over
$W_{\beta}$. If $f$ is a local homogeneous section of $\calr_*$, then
$|f|$ denotes its degree. All our DG algebras will be (graded) {\em
commutative\/} so that the product on $\calr_*$ satisfies the sign
rule
$$
fg=(-1)^{|f||g|}gf\,.
$$

Graded modules over $\calr_*$ are defined in the obvious way: An {\em
$\calr_*$--module} $\calm_*$ consists of a family of graded right
$\calr_\alpha$--modules
$\calm_\alpha=\bigoplus_{j\in\bbbz}\calm^j_{\alpha}$ together with a
compatible collection of transition maps for $\alpha\subseteq\beta$
that are degree preserving homomorphisms of right
$\calr_\beta$--modules,
$$
p_{\alpha\beta}^*(\calm_\alpha):= p_{\alpha\beta}^{-1}(\calm_\alpha)
\otimes_{p_{\alpha\beta}^{-1}(\calr_\alpha)} \calr_\beta \to
\calm_\beta\,.
$$
As we only consider commutative DG algebras, the sign rule allows to
switch the module structure from right to left, in particular one can
form the tensor product $\calm_*\otimes_{\calr_*}\caln_*$ of
$\calr_*$--modules, defining it simplex by simplex through
$(\calm_*\otimes_{\calr_*}\caln_*)_\alpha =
\calm_\alpha\otimes_{\calr_\alpha}\caln_\alpha$.

If $\calm_*$ is equipped with a differential
$\partial=\partial_{\calm_*}$ such that $(\calm_*,\partial)$ becomes a
DG module over $\calr_*$ then we call $\calm_*$ a DG $\calr_*$--module
in brief. Such a module is said to have {\em coherent cohomology} if
$\calh^i(\calm_\alpha)$, the cohomology with respect to the
differential on $\calm_\alpha$, is a coherent
$\calo_{W_\alpha}$--module for each simplex $\alpha$ and each integer
$i$; it is said to be {\em (locally) bounded above\/} if for each
simplex $\calh^i(\calm_\alpha)=0$ for $i\gg 0$, and it is said to {\em
vanish (locally) above\/} if already $\calm^{i}_{\alpha}=0$ for $i\gg
0$.
\end{sit}

\begin{sit}
Next we comment upon and fix some of the usual notations and
conventions from homological algebra that extend in a straightforward
manner to (DG) $\calr_*$--modules. We include these details in order
to be able later to calibrate Atiyah classes against Chern classes:
Over time, the conventions and signs in constructing mapping cones,
distinguished triangles and their associated extension classes have
changed, say from \cite{Ha} over \cite{Ver} to \cite{Alg}; see also
the comments in \cite[Exp.\ XVII]{SGA4} and \cite{Del3}. Whereas
classically, e.g.\ in \cite[]{At} or \cite[V.5.4.1, 5.3.3]{Ill}, the
first Chern class is the trace of the opposite of the Atiyah class, it
appears in \cite{ALJ} as the trace of the Atiyah class itself.

If $i$ is an integer, $\calm_*[i]$ is the shifted module with
$\calm_*[i]^n = \calm_*^{i+n}$ and the {\em same right
$\calr_*$--module structure\/}, whence the {\em left} structure
becomes $f(m[i]) = (-1)^{|i||f|}(fm)[i]$. In case $\calm_*$ is a DG
module, $\calm_*[i]$ becomes a DG module with respect to the
differential $\partial_{\calm_*[i]}=(-1)^i\partial_{\calm_*}$. Writing
the shift functor on the left, say as $T^i\calm_*= \calm_*[i]$, and
considering $T$ as an operator of degree $-1$, the conventions just
introduced obey the usual sign rule.

A { \em morphism of $\calr_*$--modules} $\calm_*\to\caln_*$ of degree
$i$ is a collection of homomorphisms
$f_{\alpha}:\calm_\alpha\to\caln_\alpha$ of right
$\calr_{\alpha}$--modules that satisfy
$f_{\alpha}(\calm_{\alpha}^j)\subseteq \calm^{i+j}_{\alpha}$ for all
$j\in\bbbz$ and are compatible with the transition maps. By
convention, if no degree is specified, a homomorphism of
$\calr_{*}$--modules is assumed to be of degree 0. We set
$$
\Hom_{\calr_*}(\calm_*,\caln_*):=\bigoplus_{i\in\bbbz}
\Hom^i_{\calr_*}(\calm_*,\caln_*)\,,
$$
where $\Hom^i_{\calr_*}(\calm_*,\caln_*)$ denotes the morphisms of
degree $i$. If $\calm_*$, $\caln_*$ are DG modules with differential
$\partial$ then $\Hom_{\calr_*}(\calm_*,\caln_*)$ is a complex of
vector spaces with differential $h\mapsto [\partial,h]:= \partial h
-(-1)^{|h|}h\partial$. Note that for any integer $j$ the maps $h[j]$
and $h$ are identical once the grading is ignored, and so one usually
suppresses the shift in the notation of morphisms.
\end{sit}

\begin{sit}
The morphisms of degree $i$ are thus the (degree preserving)
homomorphisms $\calm_*\to\caln_*[i]$ of $\calr_*$--modules, the cycles
of degree $i$ in the Hom--complex are the homomorphisms
$\calm_*\to\caln_*[i]$ of DG $\calr_*$--modules, whereas the
boundaries are the homomorphisms that are homotopic to zero. The {\em
homotopy category\/} $K(\calr_*)$ of DG $\calr_*$--modules has as its
objects the DG $\calr_*$--modules but its morphisms are the homotopy
classes, $\Hom_{K(\calr_*)}(\calm_*,\caln_*)=
\HH^0(\Hom_{\calr_*}(\calm_*,\caln_*))$. The {\em derived category\/}
$D(\calr_*)$ is obtained from $K(\calr_*)$ as usual by inverting all
quasiisomorphisms. Adornements such as $(\ )^{a,b}_c$ or $(\ )^{b}_c$
determine full subcategories of DG modules, where $a$ bounds the
underlying graded objects on each simplex, $b$ bounds the cohomology,
and $c$ refers to special structure of the cohomology modules. For
example, $D^{-}_{coh}(\calr_*)$ denotes the derived category of DG
$\calr_*$--modules that are bounded above with coherent cohomology. In
case $\calr_*$ is a structure sheaf, such as $\calo_{X_*}$, we write
$\Hom_{X_*}$ and $D({X_*})$ to reduce clutter.

The morphisms in $D(\calr_*)$ are denoted
$$
\Ext^k_{\calr_*}(\calm_*,\caln_*):=
\Hom_{D(\calr_*)}(\calm_*,\caln_*[k])\,,\quad k\in\bbbz\,.
$$
Similarly, if $X$ is a complex space and $\calm$, $\caln$ are
complexes of $\calox$--modules then $\Ext^k_{X}(\calm,\caln)$ will
represent the set of morphisms of degree $k$ in the derived category
$D(X)$. Note that if $\caln$ is bounded below and so admits an
injective resolution, say, $\cali$ then these Ext-groups are given as
usual by the cohomology of the complex $\Hom_X(\calm,\cali)$. This
definition allows in particular to define
$H^k(X,\caln):=\Ext^k_{X}(\calox,\caln)$ for any complex $\caln$ of
$\calox$--modules. Below we will show how to compute these groups
using projective ``resolutions" on $X_*$.
\end{sit}

\begin{sit}
\label{triangulated structure}
For a homomorphism $f:\caln'_{*}\to\caln_{*}$ of DG
$\calr_*$--modules, its {\em mapping cone\/}, $\Con_*(f):=\caln_{*}
\oplus\caln'_{*}[1]$, is formed simplex by simplex according to the
conventions of \cite[X.36ff]{Alg} so that for local sections $n$ of
$\caln_{\alpha}$ and $Tn'=n'[1]$ of $\caln'_{\alpha}[1]$ the
differential in $\Con_*(f)$ maps $(n,Tn')$ to
$(\partial_{\caln_{*}}(n) -f(n'), -T\partial_{\caln'_{*}}(n'))$. The
mapping cone is again a DG $\calr_*$--module, and the triangulated
structure of either $K(\calr_*)$ or $D(\calr_*)$ is now defined in
terms of the distinguished triangles arising from mapping cones.
\end{sit}

\begin{sit}
\label{total complex}
Generalizing the construction of mapping cones, if
$\calm^\sbullet_*\equiv( \cdots\to
\calm^{(i)}_*\xto{\delta^i}\calm^{(i+1)}_*\to \cdots)$ is a complex of
DG $\calr_*$--modules, then the associated total complex
${\bar\calm}_*=\prod_i \calm^{(i)}_*[i]$ is a DG $\calr_*$--module as
well. The following simple fact will be used throughout: If the
complex $\calm^\sbullet_*$ is acyclic and locally vanishes above, then
the associated DG module ${\bar\calm}_*$ is again acyclic.
\end{sit}

\begin{sit}
\label{left adjoint}
As just done, we will often argue {\em simplex by simplex\/} and such
arguments are alleviated by the following: For each simplex $\alpha$,
the restriction functor $\calm_*\mapsto \calm_\alpha$ is exact and
admits a left adjoint $p^*_\alpha$ defined through
$$
p^*_\alpha(\calm_\alpha)_\beta := \left\{
\begin{array}{ll}
p^*_{\alpha\beta} (\calm_\alpha) &\quad \mbox{for $\alpha \subseteq
\beta$}\\
\ns 0 &\quad \mbox{otherwise,}\end{array}\right.
$$
whence one has
\begin{equation}
\label{2.2.1}
\Hom_{\calr_*}\left(p^*_\alpha(\calm_\alpha), \caln_*\right) \cong
\Hom_{\calr_\alpha} (\calm_\alpha, \caln_\alpha)
\end{equation}
for each $\calr_\alpha$--module $\calm_\alpha$ and each
$\calr_*$--module $\caln_*$, see \cite[\S 2]{Fle1}, \cite[4.1]{BMi},
or, even more generally, \cite[VI.5.3]{Ill}. The functor $p^*_\alpha$
is compatible with differentials, it is as well a left adjoint to the
restriction functor on the category of DG $\calr_*$--modules.
\end{sit}

\begin{sit}
\label{skeleton filtration}
Inductive arguments based on the dimension of simplices often use the
{\em skeleton filtration} on an $\calr_*$--module $\calm_*$ given by
$$
\calm_{\le k}= \Imm\bigg( \bigoplus_{|\alpha|\le
k}p^*_\alpha(\calm_\alpha)\to \calm_*\bigg)\subseteq \calm_*\,.
$$
For a DG module this filtration is by DG submodules.
\end{sit}

\begin{sit}
\label{degree filtration}
Analogous facts and notions allow to argue {\em degree by degree\/} on
$W_\alpha$ for each individual simplex $\alpha$. Consider, in general,
(DG) $\calr$--modules $\calm =\bigoplus_{i\in \bbbz} \calm^i$ over a
(negatively graded DG) $\calow$--algebra $\calr =\bigoplus_{i\le 0}
\calr^i$ with $(W,\calow)$ some ringed space. For each fixed degree
$i\in\bbbz$, the restriction functor $\calm\mapsto \calm^i$ to
$\calow$--modules is exact and admits as left adjoint the functor
$\calv\mapsto \calv\otimes_{\calow}\calr[-i]$. The corresponding
canonical $\calr$--hom\-om\-orph\-ism $\bigoplus_{i\in
\bbbz}\calm^i\otimes_{\calow}\calr[-i]\to \calm$ is an epimor\-phism
of (DG) $\calr$--modules. The analogue to the skeleton filtration is
the {\em degree filtration\/}
$$
\calm^{\ge k}= \Imm\bigg( \bigoplus_{i\ge
k}\calm^i\otimes_{\calow}\calr[-i]\to\calm\bigg) \subseteq\calm\,,
$$
whence $\calm^{\ge k}$ is just the (DG) submodule of $\calm$ generated
by all homogeneous components in degrees at least $k$.
\end{sit}

Now we turn to the existence and construction of resolutions for DG
$\calr_*$--modules.

\begin{sit}
We will say in brief that an $\calr_*$--module
$\calm_*=\bigoplus_{i\in\bbbz}\calm_*^i$ has {\em coherent homogeneous
components} if each $\calm_\alpha^i$ is coherent as
$\calo_{W_\alpha}$--module. {\em By convention, our DG algebras will
always have coherent homogeneous components.} The category of all
$\calr_*$--modules with coherent homogeneous components is then
abelian and will, slightly abusively, be denoted $\Coh\calr_*$. The
restriction to a simplex $\alpha$ and its left adjoint $p^*_\alpha$
preserve coherence of homogeneous components and the $\calr_*$--module
$\bigoplus_{\alpha\in A}p^*_\alpha(\calm_\alpha)$ is in $\Coh\calr_*$
along with $\calm_*$.
\end{sit}

\begin{sit}
    \label{2.3}
A module $\calp_*\in\Coh \calr_*$ will be called {\em projective} if
it is a projective object in that category, equivalently, if the
functor $\Hom_{\calr_*}(\calp_*,-)$ is exact on $\Coh \calr_*.$ As
remarked above, such a projective module is in general not a
projective object in the category of all $\calr_*$--modules. Note also
that, in general, $\calr_*$ is not a projective module over itself. In
contrast, $\calr_\alpha$, being a module with coherent homogeneous
components over the Stein compact set $W_{\alpha}$, is projective as a
module over itself, and if $\calp_\alpha$ is projective in
$\Coh\calr_\alpha$ then $p^*_\alpha(\calp_\alpha)$ is a projective
$\calr_*$--module since a left adjoint of an exact functor preserves
projectives. A module $\calf_*\in\Coh \calr_*$ is called {\em graded
free} if it is a direct sum of modules $p^*_\alpha(\calp_\alpha)$ with
graded free $\calr_\alpha$--modules $\calp_\alpha$. Thus graded free
modules are always projective. Each category $\Coh\calr_\alpha$ has
enough projectives, for example, the graded free modules. As the
canonical $\calr_*$--hom\-om\-orph\-ism $\bigoplus_{\alpha\in
A}p^*_\alpha(\calm_\alpha)\to \calm_*$ is an epimorphism, there are
enough projectives in $\Coh\calr_*$ as well.

Using the simplicial structure, and generalizing slightly \cite[2.3,
(3)]{Fle1}, one obtains the following precise description of
projective $\calr_*$--modules.
\end{sit}

\begin{lem}
\label{2.4}
{\em 1.} An $\calr_*$--module $\calp_*$ is projective if and only if
$$
\textstyle \calp_* \cong\bigoplus_{\alpha \in A} p^*_\alpha
(\calq_\alpha)
$$
for some projective $\calr_\alpha$--modules $\calq_\alpha$ with
coherent homogeneous components.

{\em 2.} If $\calq_\alpha$ is a projective $\calr_\alpha$--module with
$\calq^i_\alpha=0$ for $i\gg 0$, then $\calq_\alpha \cong \bigoplus_{j
\in \bbbz} \calq^{(j)}_\alpha$, where $\calq^{(j)}_\alpha$ is a
projective $\calr_\alpha$--module generated in degree $j$.
\end{lem}

\begin{proof}
By induction on the dimension of simplices, it clearly suffices to
show the following. If $A'\subseteq A$ consists of all simplices
$\alpha'$ such that $\calp_\alpha =0$ for each strict subset $\alpha
\subset \alpha'$, then $\bigoplus_{\alpha' \in A'} p^*_{\alpha'}
(\calp_{\alpha'})$ is a direct summand of $\calp_*$. As $\calp_*$ is
projective the canonical epimorphism $\vp:\bigoplus_{\alpha\in
A}p^*_\alpha(\calp_\alpha)\to \calp_*$ admits a section $\psi$ . The
definition of $A'$ yields for each $\alpha'\in A'$ that
$\vp_{\alpha'}$, and thus also $\psi_{\alpha'}$, is the identity on
$(\bigoplus_{\alpha\in A}p^*_\alpha(\calp_\alpha))_{\alpha'}=
\calp_{\alpha'}$. Composing the projection $\bigoplus_{\alpha\in
A}p^*_\alpha(\calp_\alpha)\to p^*_{\alpha'}(\calp_{\alpha'})$ with
$\psi$ then retracts the natural homomorphism
$p^*_{\alpha'}(\calp_{\alpha'})\to\calp_*$ and the claim (1) follows.

Repeating the argument with respect to degrees instead of simplices
yields the second assertion.
\end{proof}

If $\calp_*$ is a projective $\calr_*$--module which also carries a DG
structure then we will say in brief that $\calp_*$ is a projective DG
$\calr_*$--module, although, of course, $\calp_*$ is not necessarily a
projective object in the category of all DG modules with coherent
homogeneous components.

For a projective module $\calp_*$ the skeleton filtration
\ref{skeleton filtration} is by direct summands: In the notation of
the preceding lemma, $\calp_{\le k}\cong \bigoplus_{|\alpha|\le k}
p^*_\alpha(\calq_\alpha)$, whence each such submodule as well as each
subquotient $\calp_{\le k}/\calp_{\le k-1}\cong \bigoplus_{|\alpha|=k}
p^*_\alpha(\calq_\alpha)$ is again projective. Analogously, if
$\calp_\alpha$ vanishes above, each submodule or subquotient with
respect to the degree filtration is again projective. If $\calp_*$,
resp. $\calp_\alpha$, is a projective DG module, then $\calp_{\le k}$,
resp. $\calp^{\ge k}_\alpha$, is not necessarily a direct summand as
DG module, but one has the following simple observations that reflect
typical reduction steps when dealing with projective DG
$\calp_*$--modules that locally vanish above.

\begin{lem}
    \label{2.5}
If $\calm_*$ and $\calp_*$ are DG $\calr_*$--modules with $\calp_*$
projective, then the natural map
$$
H^p(\Hom_{\calr_*}(\calp_*,\calm_*))\lto \Projlim_k
H^p(\Hom_{\calr_*}(\calp_{\le k},\calm_*))
$$
is bijective for each integer $p$.

If $\calm_\alpha$ and $\calp_\alpha$ are DG $\calr_\alpha$--modules
with $\calp_\alpha$ projective and $\calp^i_\alpha=0$ for $i\gg 0$,
then the natural map
$$
H^p(\Hom_{\calr_\alpha}(\calp_\alpha,\calm_\alpha))\lto \Projlim_k
H^p(\Hom_{\calr_\alpha}(\calp^{\ge k}_{\alpha},\calm_\alpha))
$$
is bijective for each integer $p$.

If $\calp_\alpha$ is generated in a single degree $k$, then there is a
natural number $s$ such that the complex
$\Hom_{\calr_\alpha}(\calp_\alpha,\calm_\alpha)$ becomes a direct
summand of $\Gamma(W_\alpha,\calm_\alpha[k])^{s}$.
\end{lem}

\begin{proof}
As said above, $0\to \calp_{\le k}\to \calp_*\to \calp_*/\calp_{\le
k}\to 0$ is an exact sequence of projective $\calr_*$--modules and so
splits, whence the induced map $\Hom_{\calr_*}(\calp_*,\calm_*)\to
\Hom_{\calr_*}(\calp_{\le k},\calm_*)$ is surjective. The first claim
follows immediately and the second one is obtained replacing skeleton
by degree filtration.

For the final assertion, note that the hypothesis yields a natural
epimorphism of DG $\calr_\alpha$--modules
$\calp^k_\alpha\otimes_{\calow_\alpha}\calr_\alpha[-k]\to
\calp_\alpha$. As $\calp_\alpha$ is projective it has in particular
coherent homogeneous components and so $\calp^k_\alpha$ can be
generated by a finite number of sections over $W_\alpha$, which in
turn provides an epimorphism $\calr_\alpha[-k]^s\to\calp_\alpha$.
Employing now fully that $\calp_\alpha$ is projective, this
epimorphism splits as a homomorphism of $\calr_\alpha$--modules. As
$\calp_\alpha$ is generated in a single degree, the differential
necessarily vanishes on the generators and the splitting is already a
morphism of DG modules. Finally use that
$\Hom_{\calr_\alpha}(\calr_\alpha[-k],\calm_\alpha) =
\Gamma(W_\alpha,\calm_\alpha[k])$.
\end{proof}

An $\calr_*$--module $\calm_*$ will be called $W_*$-{\em acyclic} if
$H^i(W_\alpha,\calm_\alpha)=0$ for all $\alpha\in A$. For instance, by
Theorem B every module $\calm_*\in\Coh \calr_*$ is $W_*$-acyclic. The
point of the definition is of course that a quasiisomorphism
$\calm_*\to \caln_*$ of $W_*$-acyclic DG modules induces
quasiisomorphisms $\Gamma(W_\alpha,\calm_\alpha)\to
\Gamma(W_\alpha,\caln_\alpha)$ over each simplex.

\begin{lem}
    \label{2.6}
For each DG $\calr_*$--module $\calm_*$ there is a $W_*$-acyclic
resolution, thus a quasiisomorphism $ \calm_*\hto\widetilde\calw_* $
from $\calm_*$ into a $W_*$-acyclic module $\widetilde\calw_*$. The
construction is functorial in $\calm_*$.
\end{lem}

\begin{proof} 
Consider on each simplex $\alpha\in A$ the canonical flabby
resolution, say, $\calm_\alpha\to\calw_\alpha^\sbullet$. As this
resolution is functorial, each component of the complex
$\calw_\alpha^\sbullet$ is naturally again a DG $\calr_\alpha$--module
and the resulting system $\calm_*\to\calw_*^\sbullet$ is a resolution
by DG $\calr_*$--modules. Now cut this resolution on the simplex
$\alpha$ at the place $d(\alpha):=\sum_{\beta\subseteq\alpha}\dim
W_\beta$ to obtain $\bar\calw_\alpha^i:=\calw_\alpha^i$ for $i<
d(\alpha)$, $\bar\calw_\alpha^i:=0$ for $i> d(\alpha)$ and
$\bar\calw_\alpha^i:=\Ker (\calw_\alpha^i\to \calw_\alpha^{i+1})$ for
$i=d(\alpha)$.

The choice of the cut-off point guarantees first, as $d(\alpha)\ge
\dim W_\alpha$, that $\bar\calw^{d(\alpha)}_\alpha$ is $W_*$-acyclic
along with the other components. Secondly, as $d(\alpha)\le d(\beta)$
for $\alpha\subseteq \beta$, the truncated complexes still form a
simplicial system. Finally, $\bar\calw^\sbullet_*$ is now a complex
that locally vanishes above and thus the canonical injection
$\calm_*\hto\widetilde\calw_*$ into the DG module obtained as the
total complex associated to $\bar\calw^\sbullet_*$ provides the
desired functorial quasiisomorphism into a $W_*$-acyclic
$\calr_*$--module.
\end{proof}

\begin{cor}
    \label{acyclic suffices}
The inclusion of the full subcategory of $D(X_{*})$ consisting of all
$W_{*}$-acyclic DG $\calr_{*}$--modules into $D(X_{*})$ is an
equivalence of triangulated categories. An inverse is given by
associating to a DG $\calr_*$--module its $W_*$-acyclic resolution.

To derive a functor on all DG $\calr_*$--modules it suffices thus to
derive it on the $W_*$-acyclic ones. \qed
\end{cor}

\begin{prop}
\label{2.7}
Let $\calm_*\to\caln_*$ be a quasiisomorphism of DG
$\calr_*$--modules and $\calp_*$ a projective DG
$\calr_*$--module that locally vanishes above.
\begin{enumerate}
\item
\label{2.7.1}
The map $\calp_*\otimes_{\calr_*}\calm_*
\to \calp_*\otimes_{\calr_*}\caln_*$ is a quasiisomorphism.

\item
\label{2.7.2}
If $\calm_*,\,\caln_*$ are $W_*$-acyclic, then
$\Hom_{\calr_*}(\calp_*,\calm_*)\lto
\Hom_{\calr_*}(\calp_*,\caln_*)$
is a quasiisomorphism. In particular, any morphism
$\calp_*\to \caln_*$ of DG $\calr_*$--modules lifts through
the given quasiisomorphism to a morphism of DG
$\calr_*$--modules $\calp_*\to \calm_*$.
\end{enumerate}
\end{prop}

\begin{proof}
Assertion (\ref{2.7.1}) can be verified locally and there
\cite[X.66,\S 4 no.3]{Alg} applies.

To prove (\ref{2.7.2}), we may assume by the first part of lemma
\ref{2.5} that $\calp_*=\calp_{\le k}$. In this case the spectral
sequence associated to the skeleton filtration on $\calp_*$ converges
and so it is sufficient to show the claim in case that $\calp_*\cong
p^*(\calp_\alpha)$ with a projective DG $\calr_\alpha$--module
$\calp_\alpha$. In view of \ref{left adjoint}(\ref{2.2.1}) this
requires to show that the corresponding map
$$
\Hom_{\calr_\alpha}(\calp_\alpha,\calm_\alpha)\to
\Hom_{\calr_\alpha}(\calp_\alpha,\caln_\alpha).\leqno (*)
$$
is a quasiisomorphism. Because of the second part of lemma \ref{2.5}
we may reduce to the situation where $\calp_\alpha$ is generated in
finitely many degrees, in which case the spectral sequence associated
to the degree filtration on $\calp_\alpha$ converges. It remains to
deal with the case that $\calp_\alpha$ is generated in a single degree
$k$ and then the final part of lemma \ref{2.5} exhibits the map $(*)$
as a direct summand of $\Gamma(W_\alpha,\calm_\alpha[k]^s)\to
\Gamma(W_\alpha,\caln_\alpha[k]^s)$ for some $s$. As $\calm_*,
\caln_*$ are $W_*$-acyclic, this last map is a quasiisomorphism and
the claim follows.
\end{proof}

\begin{cor}\label{2.8}
Any quasiisomorphism $\calp_*\to\calq_*$ between projective DG
$\calr_*$--modules that locally vanish above is a homotopy 
equivalence. 
In
particular, for every  DG $\calr_*$--module $\calm_*$  and any
quasiisomorphism $\calp_*\to\calq_*$ the induced maps
\begin{align*}
\Hom_{\calr_*}(\calq_*,\calm_*)&\lto
\Hom_{\calr_*}(\calp_*,\calm_*)\qquad \text{and}\\
\calp_* \otimes_{\calr_*}\calm_* &\lto\calq_*
\otimes_{\calr_*}\calm_*
\end{align*}
are quasiisomorphisms.
\end{cor}

\begin{proof}
It follows from \ref{acyclic suffices} and \ref{2.7}(\ref{2.7.2}) that
any quasiisomorphism between projective DG $\calr_*$--modules that
locally vanish above is represented by an actual morphism
$f:\calp_*\to\calq_*$ of DG $\calr_*$--modules.  Such a morphism is a
homotopy equivalence if the identity on the mapping cone
$\calc_*=\Con_*(f)$ is homotopic to 0.  With $\calp_*,\calq_*$ also
$\calc_*$ is projective, thus $W_*$-acyclic, and locally vanishes
above.  So the preceding result \ref{2.7}(\ref{2.7.2}) applies to
yield that $\Hom_{\calr_*}(\calc_*,\calc_*)\to \Hom_{\calr_*}(\calc_*,
0)=0$ is a quasiisomorphism.  The rest follows immediately.
\end{proof}

A quasiisomorphism $h:\calp_*\to \calm_*$ of DG $\calr_*$--modules
will be called a {\em projective approximation} of $\calm_*$ if
$\calp_*$ is a projective DG $\calr_*$--module that locally vanishes
above. The preceding corollary shows that projective approximations
of a DG $\calr_*$--module are unique up to homotopy equivalence.
Their existence is settled next.

\begin{prop}
    \label{2.9}
A DG $\calr_*$--module $\calm_*$ admits a projective approximation
if and only if it is bounded above with coherent cohomology.
\end{prop}

\begin{proof}
If $h:\calp_*\to \calm_*$ is a projective approximation, then
$\calm_*$ is necessarily bounded above with coherent cohomology as
this holds for $\calp_{*}$ and $h$ is a quasiisomorphism.
Conversely,
it suffices by \ref{2.6} to show the existence of a
projective approximation when $\calm_*$ is furthermore
$W_*$-acyclic.  But if $\calm_*$ is $W_*$-acyclic with
coherent, thus $W_*$-acyclic, cohomology, then its submodule
of boundaries, $\calb_* = \partial(\calm_*)\subseteq
\calm_*$ is also $W_*$-acyclic: the exact sequences
\begin{align*}
0\to \calk_*^i\to &\calm_*^i \to \calb_*^i\to 0\\
0\to \calb_*^{i-1}\to &\calk_*^i \to \calh^i(\calm_*)\to 0
\end{align*}
yield for each $j\ge 1$ and each simplex $\alpha\in A$
isomorphisms
$$
H^j(W_\alpha, \calb^i_\alpha)\xto{\cong\ } H^{j+1}(W_\alpha,
\calk^i_\alpha)\xleftarrow{\cong} H^{j+1}(W_\alpha,
\calb^{i-1}_\alpha),
$$
whence $H^j(W_\alpha, \calb^i_\alpha)\cong H^{j+\dim
W_\alpha}(W_\alpha,\calb^{i-\dim W_\alpha}_\alpha) = 0$ for
each $i$.

In view of this observation, the classical construction of a
projective resolution applies: as $\calh_{*}(\calm_*)$ is coherent
there exists a surjection from a graded free $\calr^0_*$--module
$\calq_*$ with coherent homogeneous components onto $\calh(\calm_*)$
such that $\calq^{i}_{\alpha}$ vanishes in degrees $i>d$ if
$\calh^{i}(\calm_{\beta})$ vanishes for those degrees for all
simplices $\beta\subseteq\alpha$. Due to the acyclicity of the
boundaries $\calb_{*}\subseteq\calm_{*}$, this homomorphism lifts from
the cohomology to the cycles $\calk_*\subseteq\calm_*$, and the
resulting morphism of DG $\calr_*$--modules
$\calp_*^{(0)}:=\calq_{*}\otimes_{\calr_{*}^0}\calr_{*} \to\calm_*$ is
surjective in cohomology. Now $\calp_*^{(0)}$ is projective, thus any
submodule is $W_*$-acyclic, and the usual process, see e.g.\
\cite[I.4.6]{Ha}, produces a complex $\cdots\to
\calp_*^{(i)}\to\cdots\to \calp_*^{(0)}\to \calm_*\to 0$ with
projective DG $\calr_*$--modules $\calp_*^{(i)}$ that becomes acyclic
when $\calh$ is applied. Moreover, on each simplex $\alpha$ the
projective modules $\calp_{\alpha}^{(i)}$ can be choosen to vanish
uniformly above so that the DG $\calr_*$--module $\calp_*$ associated
to the total complex of the $\calp_*^{(i)}$ is again a projective
$\calr_*$--module. The induced morphism $h:\calp_*\to \calm_*$ then
resolves, cf.\ \ref{total complex}, and constitutes a desired
projective approximation.
\end{proof}

\begin{exam}
    \label{exam2.4}
As we mentioned earlier in \ref{2.3}, and as can be seen easily from
\ref{2.4}, $\calr_*$ is in general not a projective module over
itself. Instead, $\calr_*$ admits the following explicit, and natural,
projective approximation. First consider the co-\v{C}ech complex
$$
\cdots\to \bigoplus_{|\alpha|=p} p_\alpha^*(\calr_\alpha) e(\alpha)
\xto{\delta}
\cdots\xto{\delta}
\bigoplus_{|\alpha|=1} p_\alpha^*(\calr_\alpha) e(\alpha)
\xto{\delta}
\bigoplus_{|\alpha|=0} p_\alpha^*(\calr_\alpha) e(\alpha)
\to 0\,,
$$
where the direct sums are indexed by ordered simplices, $\alpha =
(\alpha_0,\ldots,\hat \alpha_j,\ldots ,\alpha_p)$, and the
differential $\delta$ is dual to the differential of the \v{C}ech
complex so that $$\delta(e(\alpha_0,\ldots,\alpha_p))=\sum (-1)^j
e(\alpha_0,\ldots,\hat \alpha_j,\ldots ,\alpha_p)\,.
$$ 
This differential is a homomorphism of DG modules and each term of the
complex is a projective $\calr_*$--module by \ref{2.4}. The
augmentation map $\bigoplus_{|\alpha|=0} p_\alpha^*(\calr_\alpha)
e(\alpha) \xto{\epsilon} \calr_{*}$, sum of the natural maps
$p_\alpha^*(\calr_\alpha)\to \calr_{*}$, realizes the co-\v{C}ech
complex as a projective resolution of $\calr_*$ as an
$\calr_*$--module: the restriction of the augmented complex to any
simplex is evidently contractible.

The total complex associated to this co-\v{C}ech complex gives thus a
projective approximation $\calp_*$ of $\calr_*$, see \ref{total
complex} or the proof of \ref{2.9}. Note that $\calp_*=\bigoplus
p_\alpha^*(\calr_\alpha)[|\alpha|]$, where the sum is taken over all
ordered simplices $\alpha$.
\end{exam}

\begin{rems}
    \label{modified projective approximations}
1. It follows that a projective approximation $h:\calp_{*}\to
\calm_{*}$ can be realized as a pair of morphisms of DG
$\calr_{*}$--modules $\calp_{*}\to {\widetilde\calm_{*}}\gets
\calm_{*}$, where the morphism $\calm_{*}\to
{\widetilde\calm_{*}}$ is a $W_{*}$-acyclic resolution and
$\calp_{*}\to {\widetilde\calm_{*}}$ is a projective resolution as
constructed in the proof.

2. If so desired, one may clearly modify the above construction
to obtain in the end a projective approximation $\calf_*\to \calm_*$
with $\calf_*$ a {\em graded free\/} DG $\calr_*$--module.

3. With the notation of the proposition, if $A'\subseteq A$ is a
set of simplices such that $\calh(\calm_\alpha) = 0$ for all simplices
$\alpha$ that do not contain a simplex from $A'$, then the above
construction produces a projective approximation with $\calp_\alpha =
0$ for the same simplices.

Analogously, if $\calh^{i}(\calm_{*})=0$ for $i>i_{0}$, then the
construction provides a projective approximation $h:\calp_{*}\to
\calm_{*}$ with $\calp_{*}^{i}=0$ for $i > i_{0}$.
\end{rems}

\begin{sit}\label{Ext-functors}
Now we can describe the Ext-functors of DG $\calr_*$--modules when the
contravariant argument admits a projective approximation. Recall that
by our conventions for $\calr_*$--modules $\calm_*$, $\caln_*$ the
$\Ext^k$--functors are given by the set of morphisms $\calm_*\to
\caln_*$ of degree $k$ in the derived category so that
$$
\Ext^k_{\calr_*}(\calm_*,\caln_*)=
\Hom_{D(\calr_*)}(\calm_*,\caln_*[k])\,.
$$
The following proposition summarizes \ref{2.9},
\ref{2.7}(\ref{2.7.2}), and \ref{acyclic suffices}. The reader should
compare the result with \cite[VI.10.2.4]{Ill} in the algebraic case.
\end{sit}

\begin{prop}\label{2.22.a}
If $\calm_*$ is a DG $\calr_*$--modules that is bounded above
with coherent cohomology,  and $\caln_*$ is any DG $\calr_*$--module
then
$$
\Ext^k_{\calr_*}(\calm_*,\caln_*)\cong
H^k(\Hom _{\calr_*}(\calp_*,\tilde\caln_*))\quad\text{for}\quad 
k\in\bbbz\,,
$$
where $\calp_*\to \calm_*$ is a projective
approximation and $\caln_*\to\tilde \caln_*$ is a quasiisomorphism
into a $W_*$-acyclic DG $\calr_*$--module.
\qed
\end{prop}

\begin{sit}
Projective approximations allow as well to derive tensor products of
DG $\calr_*$--modules. Calling a DG $\calr_*$--module $\calm_{*}$ {\em
flat\/} if $\calm_{*}\otimes_{\calr_*}(\,\ )$ preserves
quasiisomorphisms of DG $\calr_*$--modules, a projective module that
locally vanishes above is flat by \ref{2.7}(\ref{2.7.1}). Indeed, for
flatness it suffices already that the module locally vanishes above
and that its restriction to each simplex is projective. For example,
$\calr_{*}$ itself is always flat. If $\calp_*\to\calm_*$ is a
projective approximation of a DG $\calr_*$--module $\calm_*$,
necessarily bounded above with coherent cohomology, and if $\caln_*$
is any DG $\calr_*$--module, then $\calp_*\otimes_{\calr_*}\caln_*$
represents $\calm_*\dotimes_{\calr_*}\caln_*$, the {\em derived tensor
product \/} of $\calm_*$ with $\caln_*$ over $\calr_{*}$.

By \ref{2.8}, the derived tensor product is well defined up to
homotopy equivalence, and by \ref{2.7}(\ref{2.7.1}) and \ref{2.8} a
pair of quasiisomorphisms $\calm_*\to\calm'_*,\,\caln_*\to\caln'_*$ of
DG $\calr_*$--modules induces a quasiisomorphism
$\calm_*\dotimes_{\calr_*}\caln_*\xto{\quad}
\calm'_*\dotimes_{\calr_*}\caln'_*$. If $\caln_*$ admits a projective
approximation as well, say $\calq_*\to \caln_*$, then
$\calm_*\otimes_{\calr_*}\calq_*$ represents
$\calm_*\dotimes_{\calr_*}\caln_*$ too. If $\caln_{*}$ is flat and
locally vanishes above then, again by \cite[X.Prop.4]{Alg}, the given
quasiisomorphism $\calp_*\to\calm_*$ induces a quasiisomorphism
$\calm_*\dotimes_{\calr_*}\caln_*\xto{\quad}
\calm_*\otimes_{\calr_*}\caln_*$ and so, for example,
$\calm_*\dotimes_{\calr_*}\calr_*\cong \calm_*$.
\end{sit}

A useful consequence of the preceding considerations is the
following result that is again well known in the algebraic
case, see \cite[VI.10.3.15]{Ill}.

\begin{cor}\label{2.11}
Let $\calr_*\to\cals_*$ be a quasiisomorphism of DG
algebras over $\calo_{W_*}$.
\begin{enumerate}
\item
\label{2.11.2}
If $\calm_*$, $\caln_*$  are  DG
$\calr_*$--modules that are bounded above with coherent
cohomology, then there are natural isomorphisms
$$
\Ext^k_{\calr_*}(\calm_*,\caln_*)\xto{\ \cong\ }
\Ext^k_{\cals_*}(\calm_*\dotimes_{\calr_*}\cals_*,
\caln_*\dotimes_{\calr_*}\cals_*)
$$
for each integer $k$.
\item
\label{2.11.3}
If $\calm_*$ is in $D^-_{coh}(\calr_*)$ and $\caln_*$
any DG $\cals_*$--module, then for each integer $k$ there is a 
natural 
isomorphism
$$
\Ext^k_{\calr_*}(\calm_*,\caln_*)\xto{\ \cong\ }
\Ext^k_{\cals_*}(\calm_*\dotimes_{\calr_*}\cals_*,
\caln_*)\,.
$$
\item
\label{2.11.4}
Restriction of scalars from $\cals_*$ to $\calr_*$ and $(\
\,)\dotimes_{\calr_*}\cals_*$ form a pair of inverse exact
equivalences between $D^-_{coh}(\cals_*)$ and $D^-_{coh}(\calr_*)$.
\end{enumerate}
\end{cor}

\begin{proof}
Let $\calp_*\to \calm_*$ and $\calq_*\to \caln_*$ be projective
approximations as DG $\calr_*$--modules.
By \ref{2.7}(\ref{2.7.1}),  $\calq_*\to \calq_* 
\otimes_{\calr_*}\cals_*$
is a quasiisomorphism and so
$$
\Ext^k_{\calr_*}(\calp_*,\calq_*)\stackrel{\cong}{\lto}
\Ext^k_{\calr_*}(\calp_*,
\calq_*\otimes_{\calr_*}\cals_*)\,.
$$
As $\Hom_{\calr_*}(\calp_*,\calq_*\otimes_{\calr_*}\cals_*) \cong
\Hom_{\cals_*}(\calp_*\otimes_{\calr_*}\cals_*,
\calq_*\otimes_{\calr_*}\cals_*)$ and
$\calp_*\otimes_{\calr_*}\cals_*$ is a projective DG $\cals_*$--module
that locally vanishes above, assertion (\ref{2.11.2}) follows from
\ref{2.22.a}. 

To obtain (\ref{2.11.3}), let $\calp_*\to \calm_*$ be as before and
$\caln_{*}\to
\widetilde\caln_{*}$ a $W_{*}$--acyclic resolution that is a morphism
of DG $\cals_*$--modules. One has then
\begin{alignat*}{2}
\Ext^k_{\calr_*}(\calm_*,\caln_*) &\cong
H^{k}(\Hom_{\calr_*}(\calp_*,\widetilde\caln_{*}))&\qquad & \text{by
\ref{2.22.a},}\\
&\cong H^{k}(\Hom_{\cals_*}(\calp_*\otimes_{\calr_*}\cals_*,
\widetilde\caln_{*}))&\qquad & \text{by adjunction,}\\
&\cong \Ext^k_{\cals_*}(\calm_*\dotimes_{\calr_*}\cals_*,
\caln_*)&\qquad &\text{by \ref{2.22.a} again.}
\end{alignat*}
For (\ref{2.11.4}), note that $(\ \,)\dotimes_{\calr_*}\cals_*$ on
$D^-_{coh}(\calr_*)$ is fully faithful by (\ref{2.11.2}). This functor
takes its values in $D^-_{coh}(\cals_*)$, and to establish it is an
equivalence with inverse as indicated, it suffices to remark that for
each $\caln_*$ in $D^-_{coh}(\cals_*)$ the natural morphism
$\caln_*\dotimes_{\calr_*}\cals_*\to \caln_*$ obtained from
(\ref{2.11.3}) is an isomorphism in $D^-_{coh}(\cals_*)$. Indeed, let
$\caln_{*}\to \widetilde\caln_{*}$ be again a $W_{*}$--acyclic
resolution that is a morphism of DG $\cals_*$--modules and choose a
morphism of DG $\calr_{*}$--modules $\calq_{*}\to \widetilde\caln_{*}$
that constitutes a projective approximation. As $\widetilde\caln_{*}$
is already a $\cals_*$--module, this quasiisomorphism factors as
$\calq_{*}\to
\calq_{*}\otimes_{\calr_*}\cals_*\to\widetilde\caln_{*}$, the morphism
$\calq_{*}\to \calq_{*}\otimes_{\calr_*}\cals_*$ is a quasiisomorphism
by \ref{2.7}(\ref{2.7.1}), and thus so is
$\calq_{*}\otimes_{\calr_*}\cals_*\to\widetilde\caln_{*}$. Now the
pair of quasiisomorphisms
$\calq_{*}\otimes_{\calr_*}\cals_*\to\widetilde\caln_{*}\gets
\caln_{*}$ represents the morphism
$\caln_*\dotimes_{\calr_*}\cals_*\to \caln_*$, whence the latter is an
isomorphism in $D^-_{coh}(\cals_*)$.
\end{proof}

\begin{sit}
\label{2.10}
Considering $\calo_{X_*}$ as a DG algebra concentrated in degree
$0$, (projective) DG $\calo_{X_*}$--modules are just complexes of
(projective) $\calo_{X_*}$--modules. In this situation, the 
restriction
functor to a simplex $\alpha$ is easily seen to admit a right adjoint,
given by
$$
p_{\alpha*}(\calm_\alpha)_\beta := \left\{
\begin{array}{ll}
p_{\beta\alpha*} (\calm_\alpha) &\quad \mbox{for $\beta \subseteq
\alpha$}\\
\ns
0 &\quad \mbox{otherwise,}
\end{array}\right.
$$
cf.\ \cite[\S 2]{Fle1}, or, again more generally,
\cite[VI.5.3]{Ill}. As it is right adjoint to an exact functor,
$p_{\alpha*}$ transforms an injective $\calo_{X_\alpha}$--module
$\cali_\alpha$ into the injective $\calo_{X_*}$--module
$p_{\alpha*}(\cali_\alpha)$. The canonical map
$\calm_*\hto \prod_{\alpha\in A}p_{\alpha*}(\calm_\alpha)$ is a
monomorphism for each $\calo_{X_*}$--module $\calm_*$, and embedding 
in
turn each $\calm_\alpha$ into an injective $\calo_{X_\alpha}$--module
$\cali_\alpha$ yields by composition a monomorphism $\calm_*\hto
\prod_{\alpha\in A}p_{\alpha*}(\cali_\alpha)$ into an injective
$\calo_{X_*}$--module. Thus the category of $\calo_{X_*}$--modules has
enough injectives. In particular, for complexes $\calm_*, \caln_*$ of
$\calo_{X_*}$--modules one may calculate $\Ext^k_{X_*}(\calm_*,
\caln_*)$ in the ``classical'' way as
$\HH^k(\Hom_{X_*}(\calm_*, \cali_*))$ if $\caln_*$ admits an
injective resolution $\caln_*\to \cali_*$. By \ref{2.22.a}, if
$\calm_*$ is a complex of $\calo_{X_*}$--modules that is bounded
above with coherent cohomology, then these groups can be calculated
using a projective approximation of $\calm_*$.
\end{sit}

\begin{sit}
\label{Cech functor}
Restricting a given $\calox$--module $\calm$ to the Stein
compact sets of the given covering defines the
$\calo_{X_*}$--module $\calm_*=j^*\calm$ with
$\calm_\alpha:=\calm|X_\alpha$. This functor is exact
and so induces directly a functor
$j^*:D(X)\to D(X_*)$ between the derived categories.

To describe a right adjoint, denote $j_\alpha:X_\alpha\hto X$ the
inclusion and associate to a module $\calm_*$ on $X_*$ the {\em
\v{C}ech complex\/} $C^\sbullet(\calm_*)$ with terms
$$
C^p (\calm_*):= \prod_{|\alpha|=p} {j_\alpha}_*(\calm_\alpha)\,,
$$
where the product is over all ordered simplices, and differential
defined in the usual way by means of the transition maps for $\calm_*$
and the given ordering on the simplices. The functor $j_*(\calm_*):=
\calh^0(C^\sbullet(\calm_*))$ is a right adjoint to $j^*$ on the
category of $\calo_{X_{*}}$--modules, and the canonical homomorphism 
of
$\calox$--modules $\calm\to j_*j^*(\calm)$ is an isomorphism.

As the given covering is locally finite, the complex
$C^\sbullet(\calm_*)$ is locally bounded; the localization at a point
$x\in X$ vanishes in degrees greater than $\max\{|\alpha|, x\in
X_\alpha\}$. As the covering is by closed sets, the functors
${j_\alpha}_*$, and then also $C^p (\ \,)$, are exact. These two facts
together imply that the total complex (associated to)
$C^\sbullet(\calm^{\sbullet}_*)$ is acyclic whenever
$\calm^{\sbullet}_*$ is an acyclic complex of $\calo_{X_*}$--modules.
Accordingly, $C^\sbullet$ can be viewed as a functor from $D(X_*)$ to
$D(X)$. Note that the terms of the complex
$C^\sbullet(\calm^{\sbullet}_*)$ are flat $\calox$--modules whenever
$\calm^\sbullet_\alpha$ is flat over $\calo_{X_\alpha}$ for each
simplex $\alpha$.

We now show that $C^\sbullet$ represents $Rj_*$, the right derived
functor of $j_*$: for an $\calo_{X_*}$--module of the form
$p_{\alpha*}(\calm_\alpha)$, the complex
$C^\sbullet(p_{\alpha*}(\calm_\alpha))$ is nothing but the usual
(sheafified) \v{C}ech complex of $\calm_\alpha$ on $X_\alpha$ with
respect to the trivial covering $\{X_\alpha\cap
X_i=X_\alpha\}_{i\in\alpha}$, then extended by zero to the rest of
$X$. Clearly $C^\sbullet(p_{\alpha*}(\calm_\alpha))$ resolves
$j_*(p_{\alpha*}(\calm_\alpha))\cong j_{\alpha*}(\calm_\alpha)$, the
extension of $\calm_\alpha$ by zero. By \ref{2.10}, each
$\calo_{X_*}$--module $\calm_*$ admits a resolution by
$\calo_{X_*}$--modules of the form $\prod_{\alpha\in
A}p_{\alpha*}(\caln_\alpha)$, whence the natural morphism of functors
$Rj_*\to C^\sbullet$, induced by the universal property of the derived
functor, is indeed an isomorphism. The relationship between $D(X)$ and
$D(X_*)$ can now be summarized as follows.
\end{sit}

\begin{prop}
\label{2.16}
The functor $j^*:D(X)\to D(X_*)$ embeds $D(X)$ as a full and exact
subcategory into $D(X_*)$ and $C^\sbullet\cong Rj_*$ is an exact right
adjoint. In particular, for
$\calm, \caln \in D(X)$ and
$\caln'_*\in D(X_*)$ there are functorial isomorphisms
\begin{align*}
    \calm&\cong C^\sbullet(j^{*}\calm)\,,\\
\Ext^k_{X}(\calm,\caln)&\cong
\Ext^k_{X_*}(j^*(\calm),j^*(\caln))\,,\quad\text{and}\\
\Ext^k_{X_*}(j^*(\calm),\caln'_*)&\cong
\Ext^k_{X}(\calm, C^\sbullet(\caln'_*))\,.
\end{align*}
The adjoint pair $j^*, Rj_*$ satisfies the {\em projection formula:\/}
the canonical morphism
$$
\calm\dotimes_{\calo_{X}}Rj_*(\caln_*)\to
Rj_*(j^*\calm\dotimes_{\calo_{X_*}}\caln_*)
$$
is an isomorphism in $D(X)$ for $\calm \in D(X)$ and
$\caln_*\in D^-_{coh}(X_*)$.
\end{prop}

\begin{proof}
As mentioned above, the natural morphism of functors $\id\to j_*j^*$ 
is 
an
isomorphism on the level of modules. It induces
an isomorphism $\id \to (Rj_*)j^*\cong C^\sbullet j^*$ of
the corresponding derived functors on $D(X)$, whence the functor 
$j^*$ 
is
still fully faithful on $D(X)$.

To prove the projection formula, observe first that the natural map
$$
j^*(\calm\otimes_{\calox}\calm')\to
j^*(\calm)\otimes_{\calo_{X_*}}j^*(\calm') 
$$ 
is an isomorphism as this holds on each simplex. As $j^*$ is exact,
this isomorphism passes to the derived tensor products as soon as
those exist. Adjunction then yields a morphism
$\calm\dotimes_{\calox}\calm'\to
Rj_*(j^*(\calm)\dotimes_{\calo_{X_*}}j^*(\calm'))$. Now set $\calm' =
Rj_*\caln_*$ and compose the corresponding morphism with the one
induced by the adjunction map $j^*Rj_*(\caln_*)\to \caln_*$ to obtain
the morphism in the projection formula. Given the existence of this
natural morphism, to establish it as an isomorphism for $\caln_*\in
D^-_{coh}(X_*)$, we first replace $\caln_*$ by a projective
approximation $\calp_*$, then use \ref{2.4} to reduce to the case
$\caln_* = p^*_\alpha(\calo_{X_\alpha})$ for some simplex $\alpha$.
Now $Rj_*(p^*_\alpha(\calo_{X_\alpha}))\cong
C^\sbullet(\calox|X_\alpha)$ is a finite complex of flat
$\calox$--modules that resolves $\calox|X_\alpha$ and
$Rj_*(j^*\calm\dotimes_{\calo_{X_*}}p^*_\alpha(\calo_{X_\alpha}))\cong
C^\sbullet(\calm|{X_\alpha})$ resolves $\calm|{X_\alpha}$. The desired
isomorphism in $D(X)$ follows thus from the obvious one for
$\calox$--modules, $\calm|{X_\alpha}\cong
\calm\otimes_{\calox}(\calo_{X}|{X_\alpha})$.
\end{proof}

\begin{exam}
If $\caln$ is any complex on $X$ then, by our conventions, its
cohomology groups $H^k(X,\caln)$ are the groups
$\Ext^k_X(\calox,\caln)$ and so can be computed by the complex
$\Hom_{X_*}(\calp_*, \tilde\caln_*)$, where $\calp_*\to\calo_{X_*}$ is
a projective resolution of $\calo_{X_*}$ as a module over itself and
$\caln_*\to\tilde\caln_*$ is a quasiisomorphism into a $W_*$-acyclic
complex. Using the projective approximation $\calp_*$ of $\calo_{X_*}$
constructed in \ref{exam2.4} it follows that $\Hom_{X_*}(\calp_*,
\tilde\caln_*)$ is the usual \Cech-complex $\Gamma(X,
C^\sbullet(\tilde\caln_*))$.
\end{exam}

\begin{rem}
    \label{2.28}
   
1.\ Note that the preceding constructions work more generally on the
category of all sheaves of abelian groups $\Ab_{X}$ on $X$ and
$\Ab_{X_{*}}$ on $X_{*}$. In particular, later on we will need the
\v{C}ech functor in this context, where it is still exact and defines
with the same arguments as before the derived functor $Rj_{*}$ from
$D(\Ab_{X_{*}})$ into $D(\Ab_{X})$.

2.\ As $j^{*}$ is fully faithful, the essential image of $D(X)$ under 
this functor is easily seen to consist of all complexes $\caln_{*}$ in
$D(X_{*})$ for which the transition maps 
$p_{\alpha\beta}^*(\caln_\alpha)\to 
\caln_{\beta}$ are quasiisomorphisms for all simplices $\alpha, 
\beta$ with
$\alpha\subseteq\beta$. In this case the natural map
$j^{*}C^\sbullet (\caln_{*})\to \caln_{*}$ is a quasiisomorphism, and
\ref{2.16} implies that for each complex $\calm\in D(X)$, there are
isomorphisms
$$
\Ext^i_{X_*}(\caln_{*}, j^{*}\calm)\cong
\Ext^i_{X}(C^\sbullet(\caln_{*}), \calm)\, , \quad i\in\bbbz\,.
$$
\end{rem}

We now turn to the construction of algebra resolutions and resolvents.

\begin{sit}
    A morphism $\calr_*\to\widetilde\calr_*$ of DG algebras
    over ${W_*}$ is called {\em free\/} if there is a graded
    free DG $\calr_*$--module $\calf_*$ such that
    $\widetilde\calr_*\cong \Sym_{\calr_*}(\calf_*)$ as
    $\calr_*$-algebras, where the symmetric algebra functor
    is to be understood in the graded context.  Note that
    the composition of free morphisms of DG algebras is
    again a free morphism.

    If $\calr_*\to\cals_*$ is any morphism of DG algebras over 
${W_*}$, 
a
    factorization $\calr_*\to\widetilde\calr_*\to \cals_*$ with
    $\calr_*\to\widetilde\calr_*$ free and $\widetilde\calr_*\to 
\cals_*$ a
    surjective quasiisomorphism of DG algebras will be called a DG
    {\em algebra resolution\/} of $\cals_*$ over $\calr_*$.

    The result \ref{2.11} together with the next one are the crucial
    ingredients in the construction of cotangent complexes following 
Quillen's
    original approach \cite{Qui1}, \cite{Qui2}. In his framework of 
closed
    model categories, the free morphisms are the cofibrations and the
    surjective quasiisomorphisms are the acyclic fibrations.
\end{sit}

\begin{prop}
\label{models}
Every morphism $\calr_*\to\cals_*$ of DG algebras admits a DG algebra
resolution.
\end{prop}

\begin{proof}
According to our general assumption on DG algebras, the
$\calr_*$--module $\cals_*$ is in $\Coh\calr_*$, and so \ref{2.9} and
\ref{modified projective approximations}(a),(b) guarantee a projective
approximation of $\cals_*$ in form of a morphism of DG
$\calr_*$--modules from a graded free DG $\calr_*$--module, say
$\calf^{(0)}_*\to \cals_*$.  The induced morphism $\calr^{(0)}_* :=
\Sym_{\calr_*}(\calf^{(0)}_*)\to \cals_*$ of DG algebras induces a
surjection in cohomology and the structure map
$\calr_*\to\calr^{(0)}_*$ is free.

Now assume constructed for some integer $k\ge 0$ a morphism of DG
$\calr_*$--algebras $\calr^{(k)}_*\to\cals_*$ with $\calr^{(k)}_*$
free over $\calr_*$ that is surjective in cohomology and such that
$\calh^i(\calr^{(k)}_*)\to \calh^i(\cals_*)$ is an isomorphism for
$i>-k$. As the kernel of the surjection $\calh^{-k}(\calr^{(k)}_*)\to
\calh^{-k}(\cals_*)$ is coherent, one may choose a graded free
coherent $\calo_{W_*}$--module $\calf^{-k-1}_*$ that is concentrated
in degree $-k-1$ and a morphism $\vp: \calf^{-k-1}_*\to \calr^{(k)}_*$
of graded $\calo_{W_*}$--modules of degree 1 that maps
$\calf^{-k-1}_*$ into the cycles of $\calr^{(k)}_*$ and such that the
sequence of $\calo_{W_*}$--modules
$$
\calf^{-k-1}_* \lto
\calh^{-k}(\calr^{(k)}_*)\lto \calh^{-k}(\cals_*)\lto 0
$$
is exact. Now set
$$
\calr^{(k+1)}_*\cong \Sym_{\calr^{(k)}_*}(\calf^{-k-1}_*
\otimes_{\calo_{W_*}}\calr^{(k)}_*) \cong \Sym_{\calr_*}
(\calf^{-k-1}_*\otimes_{\calo_{W_*}}\calr_{*})
\otimes_{\calr_*}\calr^{(k)}_*
$$
and use $\vp$ to extend the given differential. This DG algebra is
free over $\calr^{(k)}_*$ and the structure map is an isomorphism in
degrees greater than $-(k+1)$. The composed map, say, $\gamma:
\calf^{-k-1}_*\to \calr^{(k)}_* \to \cals_{*}$ maps $ \calf^{-k-1}_*$
into the boundaries, hence we can find a lifting $\tilde\gamma:
\calf^{-k-1}_*\to \cals_{*}^{-k-1}$ so that $\gamma=\partial
\tilde\gamma$. There is a unique homorphism of
$\calr^{(k)}_*$-algebras $\calr^{(k+1)}_*\to\cals_{*}$ that restricts
to $\tilde \gamma$ on $\calf^{-k-1}_*$. By construction it is also a
morphism of DG algebras, and it is surjective in cohomology with
$\calh^i(\calr^{(k+1)}_*)\to \calh^i(\cals_*)$ an isomorphism for each
$i> -(k+1)$. Finally set $\widetilde\calr_*=
{\underrightarrow{\lim\vphantom{p}}}_k\calr^{(k)}_*$.
\end{proof}

To remember the graded free DG $\calr_*$--modules
$\calf^{-k}_*\otimes_{\calo_{W_*}}\calr_{*}$ that were successively
adjoined in the construction of the algebra resolution we also write $
\widetilde\calr_*=\calr_*[\calf^{-k}_* \otimes_{\calo_{W_*}}\calr_{*};
k \ge 0]\,, $ and the sequence of free DG algebra morphisms
$$
\calr_*\hto\cdots\hto \widetilde\calr^{(k)}_* :=
\calr_*[\calf^{-i}_*\otimes_{\calo_{W_*}}\calr_{*} ; k \ge i
\ge 0]\hto \cdots\hto\widetilde\calr_*
$$
is sometimes called the associated {\em Postnikov tower\/}.

\begin{sit}
    \label{2.17}
Following \cite{Pal,Fle1} we will introduce the notion of a resolvent.
Given a morphism of complex spaces $X\to Y$, a {\em resolvent} for $X$
over $Y$ consists in a triple $(X_{*}, W_*,\calr_*)$ satisfying the
following conditions.
\begin{enumerate}
\item 
$X_*$ is the simplicial space associated to some locally finite
covering $(X_i)_{i\in I}$ of $X$ by Stein compact subsets, see
\ref{2.2};

\item 
$X_*\hookrightarrow W_*$ is a smoothing of $X\to Y$ as in \ref{2.2};

\item 
$\calr_*$ is a free DG algebra resolution of $\calo_{X_*}$ over
$\calo_{W_*}$.
\end{enumerate}
\end{sit}
Given $X_{*}$ and $W_{*}$, we will sometimes also refer to $\calr_{*}$
as a resolvent of $\calox_{*}$. The preceding results have the
following application.

\begin{cor}
\label{resolvent}
A resolvent $(X_{*}, W_*,\calr_*)$ of $X$ over $Y$ exists and one may
assume that $\calo_{W_*}\to \calr^0_*$ is an isomorphism, thus that
$\calr_* = \calo_{W_*}[\calf^{-k}_*; k \ge 1]$ with suitable graded
free $\calo_{W_*}$--modules $\calf^{-k}_*$ concentrated in degree
$-k$.

The induced functor $(\ \,)\dotimes_{\calr_*}\calo_{X_*}:
D^-_{coh}(\calr_*)\to D^-_{coh}(X_*)$ is an exact equivalence of
categories. \qed
\end{cor}

\begin{rem}\label{2.33}
Up to this point, all results are valid with insignificant 
modifications
for a finitely presented morphism $X\to Y$ of arbitrary locally 
noetherian
schemes if one replaces ``Stein compact set" by ``affine scheme".
If one further replaces ``coherent'' by ``quasi-coherent'', the 
preceding
results hold even for any morphism of schemes.
\end{rem}

We finish this section with the relevant results on cotangent
complexes, and here we need characteristic zero, as otherwise DG 
algebra
resolutions are not sufficient, \cite{Qui2}, \cite{Qui3}.

\begin{sit}
    \label{2.34}
Let $\calr_*$ be a DG algebra over $W_*$ as before. By definition,
$\calo_{W_\alpha\times_YW_\alpha}\cong
\calo_{W_\alpha}\tilde\otimes_{\caloy}\calo_{W_\alpha}$, where
$\tilde\otimes$ denotes the {\em analytic tensor product\/}.
Abusively, we set
$$
\calr_*\otimes_{\caloy}\calr_* :=
\calr_*\otimes_{\calo_{W_*}}(\calo_{W_*}
\tilde\otimes_{\caloy}\calo_{W_*})
\otimes_{\calo_{W_*}}\calr_*\,,
$$
and note that $\calr_*\otimes_{\caloy}\calr_*$ is naturally
a DG algebra over the smooth simplical scheme $W_*\times_Y
W_* = \{W_\alpha\times_Y W_\alpha\}_{\alpha\in A}$ of Stein
compact sets, see \ref{2.2}.

Let $\mu: \calr_*\otimes_{\caloy}\calr_*\to \calr_*$ denote
the multiplication map and set $\cali_{*}:= \ker\mu \subseteq
\calr_*\otimes_{\caloy}\calr_*$.
The DG $\calr_*$--module $\Omega^1_{\calr_*/Y}
=\cali_{*}/\cali_{*}^2$ is the {\em module of (analytic) differential
1-forms\/} of $\calr_*$ over $Y$. As $\calr_*$ has coherent 
homogeneous
components by hypothesis, $\Omega^1_{\calr_*/Y}$ is a DG
$\calr_*$--module in $\Coh(\calr_*)$. The {\em universal
derivation\/}  $d:\calr_* \to \Omega^1_{\calr_*/Y}$ maps
 a local section $f$ of $\calr_{\alpha}$ to the class
$$
df = 1\otimes f-f\otimes 1\in \cali_{\alpha}\bmod \cali_{\alpha}^2\,.
$$
It is a map of degree zero that has the desired universal property
with respect to homogeneous $Y$-derivations into graded 
$\calr_*$--modules
in $\Coh \calr_*$. More precisely, with
$\Der_{Y}^{i}(\calr_{*},\caln_{*})$ the group of $Y$-derivations of 
degree $i$ into the $\calr_{*}$--module $\caln_{*}$, one has a 
natural inclusion
\begin{diagram}
    \Hom^{i}_{\calr_{*}}(\Omega^1_{\calr_*/Y}, \caln_{*})
    &\rInto^{\quad(\ \,)\circ d\quad} 
&\Der_{Y}^{i}(\calr_{*},\caln_{*})
\end{diagram}
that becomes an isomorphism for $\caln_{*}\in \Coh \calr_*$. As the
classes $df$ locally generate $\Omega^1_{\calr_*/Y}$, its
differential, inherited from $\calr_*\otimes_{\caloy}\calr_*$, is
uniquely determined through
$$
\partial(df) = d(\partial f)\,.
$$

The module of (analytic) differential forms of degree $k \ge 0$ is the
DG $\calr_*$--module $\Omega^k_{\calr_*/Y} :=\Lambda^k_{\calr_*}
\Omega^1_{\calr_*/Y}$, the {\em alternating\/} or {\em (graded)
exterior power\/} of $\Omega^1_{\calr_*/Y}$, see, for example,
\cite[5.4.3]{Lod}. These DG modules are again in $\Coh \calr_*$. For
later use we recall that (graded) symmetric and exterior powers are
related through $\Sym^k_{\calr_*} (\Omega^1_{\calr_*/Y}[1]) \cong
\Omega^k_{\calr_*/Y}[k]$.
\end{sit}

\begin{sit}
\label{2.35}
The differential module $\Omega^1_{\calr_*/Y}$ of a free algebra over
$\calo_{W_*}$ is a graded free $\calr_*$--module and so are the
exterior powers $\Omega^k_{\calr_*/Y}$. The natural maps 
$$
\Omega^k_{\calr_*/Y}\lto
\Omega^k_{\calr_*/Y}\otimes_{\calr_*}\calo_{X_*} \quad \text{for
$k\ge 0$}
$$ 
are thus quasiisomorphisms by \ref{2.6}. The \v{C}ech-complex
$C^\sbullet(\Omega^1_{\calr_*/Y} \otimes_{\calr_*}\calo_{X_*})$ on $X$
is a {\em cotangent complex\/} of $X$ over $Y$ and is denoted
$\bbbl_{X/Y}$. By construction, it is a complex in $D^-_{coh}(X)$
whose terms are flat $\calox$--modules. The isomorphism class of
$\bbbl_{X/Y}$ in $D(X)$ is well defined, see \cite{Fle1}, in the sense
that it does not depend on the choice of the resolvent
$(X_{*},W_{*},\calr_{*})$. We set furthermore
$$
\Lambda^k\bbbl_{X/Y}:=C^\sbullet
(\Omega^k_{\calr_*/Y}\otimes_{\calr_*}\calo_{X_*})\quad \text{for
$k\ge 0$}\,.
$$
These complexes are also in $D^-_{coh}(X)$ and their isomorphism
classes are again well defined.  Indeed, they represent the {\em
derived exterior powers\/} of $\bbbl_{X/Y}$ in the sense of
\cite{DPu,Qui2} or \cite[I.4.2.2]{Ill}.

The functoriality of the formation of the cotangent complex and its
powers has the following consequence: the natural
morphisms of complexes of $\calo_{X_{\alpha}}$--modules
\begin{equation}
    \label{localization of L}
\left.\Lambda^k\bbbl_{X/Y}\right|_{X_{\alpha}}\to
\Omega^k_{\calr_{\alpha}/Y}\otimes_{\calr_{\alpha}}\calo_{X_{\alpha}}
\end{equation}
are quasiisomorphisms for each simplex $\alpha$ and each $k\ge 0$.
\end{sit}

Recall that the tangent cohomology functors of $X$ over $Y$ are
defined by
$$
T^i_{X/Y}(\caln):=\Ext^i_X(\bbbl_{X/Y},\caln)\,,\quad i\in \bbbz\,,
$$
for any $\calox$--module, or, more generally, for any
complex $\caln$ in $D(X)$. The
following more explicit description in terms of resolvents will be
frequently used in this paper.

\begin{prop}
    \label{2.36}
If $(X_*, W_*,\calr_*)$ is a resolvent of $X$ over $Y$ then for every
complex $\caln\in D(X)$ and each integer $i$ there are canonical
isomorphisms
$$
T^i_{X/Y}(\caln) \cong
\Ext^i_{\calr_*}(\Omega^1_{\calr_*/Y},\caln_*) \cong H^i(
\Hom_{\calr_*}(\Omega^1_{\calr_*/Y},\tilde\caln_*))\,,
$$
where $\caln_{*}\to\tilde\caln_{*}$ is a quasiisomorphism of
$\caln_{*}:=j^*(\caln)$ into a $W_{*}$-acyclic complex of
$\calo_{X_{*}}$--modules. Moreover, if $\caln$ is a complex of 
coherent
$\calox$--modules then these groups are as well isomorphic to
$H^i(\Der_{Y}(\calr_*,\caln_*))$.
\end{prop}

\begin{proof}
In view of \ref{2.28}(2), the quasiisomorphism (\ref{localization of 
L}) for $k=1$ shows that
$$
\Ext^i_{\calo_{X_*}}(
\Omega^1_{\calr_*/Y}\otimes_{\calr_*}\calo_{X_*},\caln_*)
\cong
\Ext^i_X(\bbbl_{X/Y},\caln)\,,
$$
and the term on the right is isomorphic to
$T^i_{X/Y}(\caln)$. According to \ref{2.11} (3), the
term on the left is isomorphic to
$\Ext^i_{\calr_*}(\Omega^1_{\calr_*/Y},\caln_*)$, and so the first
isomorphism follows. The second one follows from the fact
that $\Omega^1_{\calr_*/Y}$ is a projective $\calr_*$--module. The
final assertion is a consequence of the usual universal property of
the module of analytic differentials, see \ref{2.34}.
\end{proof}

\section{The Atiyah class}

We first define and investigate Atiyah classes and their powers on
$D^{-}_{coh}(\calr_{*})$ for any DG algebra $\calr_{*}$ over $W_{*}$
in terms of connections, then descend to $D^{-}_{coh}(X)$ by means of
the \Cech\ functor. We keep the notations of the preceding section,
and unadorned tensor products will be over $\calr_{*}$. Moreover
$\partial$ will indiscriminately denote the differentials of the
respective DG modules.

\subsection*{Atiyah classes via connections}
Recall that a {\em connection\/} on an $\calr_*$--module $\calm_*$ is
a map of degree 0,
$$
\nabla:\calm_* \lto \calm_* \otimes \Omega^1_{\calr_*/Y}\,,
$$
that satisfies the usual product rule, $\nabla(mf)= \nabla(m)f+
m\otimes df$, for local sections $m$ in $\calm_\alpha$ and $f$ in
$\calr_\alpha$.

We first collect some basic and simple facts about connections.

\begin{lem}\label{2.12}
Every projective $\calr_*$--module that locally vanishes above admits
a connection.
\end{lem}

\begin{proof}
As the direct sum of a family of connections is again a connection, we
may restrict by \ref{2.3} to the case that $\calp_*\cong
p^*_{\alpha}(\calp_\alpha)$, where $\calp_\alpha$ is a projective
$\calr_\alpha$--module generated in a single degree, say $k$. If
$\calp_\alpha$ is graded free, then $\calp_\alpha \cong
V[k]\otimes_\bbbc\calr_\alpha$ with $V$ a finite dimensional vector
space over $\bbbc$, and the collection of maps
$$
1\otimes d: V[k]\otimes_{\bbbc} \calr_\beta \to V[k]
\otimes_{\bbbc} \Omega^1_{\calr_\beta/Y}\,,\quad
\alpha\subseteq\beta\,,
$$
defines a connection on $V[k]\otimes_{\bbbc}
p^*_{\alpha}(\calr_\alpha)$. In the general case $\calp_\alpha$ embeds
into a free module such that $\calf_\alpha \cong
V[k]\otimes_\bbbc\calr_\alpha\cong \calp_\alpha\oplus\calq_\alpha$. If
$\nabla:p^*_{\alpha}(\calf_\alpha)\to
p^*_{\alpha}(\calf_\alpha)\otimes \Omega^1_{\calr_*/Y}$ is a
connection on $p^*_{\alpha}(\calf_\alpha)$, then the composition
$$
p^*_{\alpha}(\calp_\alpha)\xto{incl.}
p^*_{\alpha}(\calf_\alpha)\xto{\nabla}
p^*_{\alpha}(\calf_\alpha)\otimes \Omega^1_{\calr_*/Y}
\xto{proj.}
p^*_{\alpha}(\calp_\alpha)\otimes \Omega^1_{\calr_*/Y}
$$
is easily seen to be a connection on $p^*_{\alpha}(\calp_\alpha)$.
\end{proof}

\begin{prop} \label{2.13}
For any connection $\nabla: \calm_* \to \calm_*\otimes
\Omega^1_{\calr_*/Y}$ on a DG $\calr_*$--module $\calm_*$, the map
$[\partial,\nabla]$ of degree $1$ is a homomorphism of DG
$\calr_*$--modules, so that
$$
[\partial,\nabla] = \partial\nabla - \nabla\partial \in
\Hom^1_{\calr_*} (\calm_*, \calm_* \otimes \Omega^1_{\calr_*/Y})
$$
is a cycle. Its cohomology class in $\HH^1(\Hom_{\calr_*} (\calm_*,
\calm_* \otimes \Omega^1_{\calr_*/Y}))$ is independent of the choice
of connection.
\end{prop}

\begin{proof}
That $[\partial,\nabla]$ is a homomorphism
of right ${\calr_*}$--modules is easily verified by explicit
calculation. Moreover,
$$
\bigl[\partial,[\partial,\nabla]\bigr] = \partial[\partial,\nabla] +
[\partial,\nabla]\partial
= \partial(\partial\nabla - \nabla\partial) +
(\partial\nabla - \nabla\partial)\partial = \partial^2 \nabla -
\nabla \partial^2 = 0\,,
$$
whence $[\partial,\nabla]$ is a homomorphism of DG modules, thus a
cycle.

If $\nabla_1, \nabla_2: \calm_* \to
\calm_*\otimes\Omega^1_{\calr_*/Y}$ are connections, then
$\nabla_1 - \nabla_2$ is ${\calr_*}$--linear and so
$[\partial,\nabla_1] = [ \partial,\nabla_2] + [\partial,\nabla_1 -
\nabla_2]$, which means that the cycles $[\partial,\nabla_1],
[\partial,\nabla_2]$ are cohomologous.
\end{proof}

This construction of a well defined cycle from a connection can be
iterated: for each $k\ge 0$, the cycle $[\partial,\nabla]$ defines a
morphism  of DG $\calr_*$--modules of degree $k$,
$$
[\partial,\nabla]^k:\calm_* \lto
\calm_*\otimes(\Omega^1_{\calr_*/Y})^{\otimes k}
\xto{1\otimes \wedge^k}
\calm_*\otimes\Omega^k_{\calr_*/Y}\,,
$$
where $\wedge^k: (\Omega^1_{\calr_*/Y})^{\otimes k}\to
\Omega^k_{\calr_*/Y}$ is the natural projection.  The class of
$[\partial,\nabla]^k$ in $\HH^{k}(\Hom_{\calr_*} (\calm_*, \calm_*
\otimes \Omega^k_{\calr_*/Y}))$ is again independent of the chosen
connection as follows from the case $f = \id_{\calm_*}$ in the next
result that deals more generally with the functoriality of these
iterated classes.

\begin{lem}\label{2.14}
Let $f:\calm_*\to\calm'_*$ be a morphism of DG $\calr_*$--modules. If
$$
\nabla: \calm_* \to \calm_*\otimes\Omega^1_{\calr_*/Y}
\quad\text{and}\quad
\nabla': \calm'_* \to \calm'_*\otimes\Omega^1_{\calr_*/Y}
$$
are connections, then $(f\otimes 1)\circ[\partial,\nabla]^{k}$ and
$[\partial,\nabla']^{k}\circ f$ represent the same class in
$\HH^{k}(\Hom_{\calr_*} (\calm_*, \calm'_* \otimes
\Omega^k_{\calr_*/Y}))$.
\end{lem}

\begin{proof}
If $k=0$, the classes in question are equal to $f$ itself. If $k>0$,
the map
$$
g:= (f\otimes 1)\circ\nabla\circ [\partial,\nabla]^{k-1} -
\nabla'\circ [\partial,\nabla']^{k-1}\circ f:
\calm_*\to\calm'_*\otimes\Omega^k_{\calr_*/Y}
$$
of degree $k-1$ is $\calr_*$-linear. As $[\partial,g]$ is equal to
$(f\otimes 1)\circ[\partial,\nabla]^{k}-
[\partial,\nabla']^{k} \circ f$ the claim follows.
\end{proof}

Applying the preceding result when $f$ is a homotopy equivalence
between projective approximations of the same DG $\calr_*$--module
shows that the following definition is independent of the choice of
projective approximation or connection on it.

\begin{defn}\label{2.15}
Let $\calm_*$ be a DG $\calr_*$--module that is bounded above with
coherent cohomology.  Let $\calp_*\to\calm_*$ be a projective
approximation as in \ref{2.9} and let $\nabla$ be a connection on
$\calp_{*}$ that exists by \ref{2.12}. The {\em Atiyah class\/} of
$\calm_{*}$ with respect to $\calr_{*}/Y$ is the image of
$[\partial,\nabla]$ under the isomorphism
\begin{align*}
    \HH^{1}(\Hom_{\calr_*} (\calp_*, \calp_* \otimes
    \Omega^1_{\calr_*/Y}))
    &\cong
    \Ext^1_{\calr_*}(\calm_*,\calm_*\dotimes\Omega^1_{\calr_*/Y})\\
    [\partial,\nabla]\quad&\mapsto \quad\At(\calm_*)\,,
\end{align*}
and the class of $[\partial,\nabla]^{k}$ is mapped to the {\em $k$-th
power of the Atiyah class \/} of $\calm_*$,
$$
[\partial,\nabla]^{k}\mapsto\At^{k}(\calm_*)\in
\Ext^k_{\calr_*}(\calm_*,\calm_*\dotimes\Omega^k_{\calr_*/Y})\,.
$$
\end{defn}

\begin{sit}
Now let $\vp:\calr_*\to\cals_*$ be a morphism of DG algebras over
${W_*}$ and $\calm_{*}$ an $\calr_*$--module with connection $\nabla$.
The composition
$$
\calm_{*}\xto{\nabla\otimes 1}
\calm_{*}\otimes\Omega^1_{\calr_*/Y}\otimes \cals_*\xto{1\otimes d\vp}
\calm_{*}\otimes\Omega^1_{\cals_*/Y}\xto{\cong} (\calm_{*}\otimes
\cals_*)\otimes_{\cals_*}\Omega^1_{\cals_*/Y}
$$
extends by the product rule to a connection $\nabla_{\cals_*}$ on the
$\cals_*$--module $\calm_{*}\otimes \cals_*$.  If $\calm_{*}$ is a DG
$\calr_*$--module, one verifies easily that
$$
[\partial, \nabla_{\cals_*}]^{k} = (1\otimes
\wedge^{k}d\vp)\circ([\partial,\nabla]^{k}\otimes \cals_*)\,,
$$
where $ \wedge^{k}d\vp :
(\Lambda^{k}_{\calr_*}\Omega^1_{\calr_*/Y})\otimes \cals_* \cong
\Lambda^{k}_{\cals_*}(\Omega^1_{\calr_*/Y}\otimes \cals_*)\to
\Lambda^{k}_{\cals_*}\Omega^1_{\cals_*/Y}$ is the morphism of DG
$\cals_*$--modules induced by $d\vp$.  These considerations imply the
following result.
\end{sit}

\begin{prop}\label{prop:base change Atiyah}
Let $\vp:\calr_*\to\cals_*$ be a morphism of DG algebras over ${W_*}$
and $\calm_{*}\in D^{-}_{coh}(\calr_*)$.  Under the natural maps
induced by $d\vp$,
$$
(\wedge^{k}d\vp)_{*}:
\Ext^{k}_{\calr_*}(\calm_{*},\calm_{*}\dotimes\Omega^k_{\calr_*/Y})
\xto{\quad} \Ext^{k}_{\cals_*}(\calm_{*}\dotimes \cals_*,
\calm_{*}\dotimes\Omega^k_{\cals_*/Y})\,,
$$
the powers of the Atiyah class of $\calm_{*}$ are mapped to those of
$\calm_{*}\dotimes \cals_*$.
\qed
\end{prop}

If $\vp$ is a quasiisomorphism, then $(\wedge^{k}d\vp)_{*}$ is an
isomorphism for each $k$ and in this sense the exact equivalence of
triangulated categories $(\quad)\dotimes \cals_*:
D^{-}_{coh}(\calr_*)\to D^{-}_{coh}(\cals_*)$ from
\ref{2.11}.\ref{2.11.4} commutes with powers of Atiyah classes.

\begin{sit}
    \label{2.19}
The definition of a connection does not involve the grading of the
underlying $\calr_{*}$--module.  In particular, if $\nabla$ is a
connection on the DG $\calr_*$--module $\calm_{*}$, it is as well one
on the shifted module $\calm_{*}[i]$ for each integer $i$.  As
$\partial_{\calm_{*}[i]}= (-1)^{i}\partial_{\calm_{*}}$, the canonical
identification
$$
\Ext^k_{\calr_*}(\calm_*,\calm_*\dotimes\Omega^k_{\calr_*/Y})
\xto[\cong]{\displaystyle \ [i]\ }
\Ext^k_{\calr_*}(\calm_*[i],\calm_*[i]\dotimes\Omega^k_{\calr_*/Y})
$$
for $\calm_*$ in $D^{-}_{coh}(\calr_*)$ maps $\At^{k}(\calm_*)$ to
$(-1)^{ki}\At^{k}(\calm_*[i])$, thus, in short,
$$
\At^{k}(\calm_*[i]) = (-1)^{ki}\At^{k}(\calm_*)[i]\,.
$$
\end{sit}

\begin{sit}
The sign in this last equality disappears if one changes slightly the
point of view: recall that for any DG $\calr_*$--modules $\calm_{*},
\caln_{*}$ and any integer $k$ one has a natural isomorphism of DG
$\calr_*$--modules that moves the shift functor $T^{k}$ to the
second factor,
$$
(\calm_{*}\otimes \caln_{*})[k]\xto{\cong}
\calm_{*}\otimes (\caln_{*}[k])
\,,\quad
T^{k}(m\otimes n) \mapsto (-1)^{k|m|}m\otimes T^{k}n\,,
$$
where $m$ is a local homogeneous section of $\calm_{*}$ and $n$ is a
local section of $\caln_{*}$.
    
Composing $\At^{k}(\calm_{*})$ with this isomorphism yields then a
morphism
$$
\calm_{*} \xto{\At^{k}(\calm_{*})}
(\calm_*\dotimes\Omega^k_{\calr_*/Y})[k]\xto
{\cong}\calm_*\dotimes(\Omega^k_{\calr_*/Y}[k]) \xto{\cong}
\calm_*\dotimes \bbbs^{k}(\Omega^1_{\calr_*/Y}[1])
$$
that we denote, abusively, again by $\At^{k}(\calm_{*})$. In this 
form,
\ref{2.14} and \ref{2.19} translate into the following
\end{sit}

\begin{prop}\label{thm:Atiyah}
The powers of the Atiyah class define morphisms of exact
functors
$$
\At^{k}(\ \,): \id \xto{\quad} (\ \,)\dotimes
\bbbs^{k}(\Omega^{1}_{\calr_{*}/Y}[1])\,,\quad k \ge 0\,,
$$
on $D_{coh}^{-}(\calr_{*})$ that commute with the shift functor.\qed
\end{prop}

\subsection*{Atiyah classes of coherent $\calox$--modules}
Now we descend to $X$. 
Let $\calf\in D^-_{coh}(X)$ be a complex and $(X_{*}, W_*,\calr_*)$
a resolvent of the morphism of complex spaces $X\to Y$.  Associating
to these data the complex of $\calo_{X_*}$--modules $\calf_*$ with
$\calf_\alpha=\calf|X_\alpha$, the powers of the Atiyah classes of
$\calf$ are defined to be the images of $\At^k(\calf_*)$ under the
isomorphism
\begin{align*}
\Ext^k_{X_*}\left(\calf_*, \calf_*\dotimes\Omega^k_{\calr_*/Y}\dotimes
\calo_{X_{*}}\right)
&\cong
\Ext^k_X\left(\calf,\calf\dotimes\Lambda^k\bbbl_{X/Y}\right)\,,\\
\At^k(\calf_*)&\mapsto \At^k(\calf)\,.
\end{align*}

\begin{theorem} \label{2.18}
The Atiyah classes $\At^k(\calf) \in \Ext^k_X(\calf,
\calf\otimes\Lambda^k\bbbl_{X/Y})$ are well defined for each $\calf\in
D^-_{coh}(X)$.
\end{theorem}

\begin{proof}
To compare the Atiyah classes formed with DG algebra resolutions
$\calr_*$ and $\calr'_*$ of $\calo_{X_*}$ over $W_*$, note first that 
by
\ref{models} there is a free
$\calr_*\otimes_{\calo_{W_*}}\calr'_*$ algebra that constitutes a
resolvent of $\calo_{X_*}$. Thus we may suppose that
$\calr'_*$ is a free $\calr_*$-algebra. 
In this case \ref{prop:base change Atiyah} gives the independence 
from the choice of DG algebra resolution.

The construction is as well independent of the embedding $X_*\subseteq
W_*$. With similar arguments as above it suffices to
compare two embeddings $X_*\subseteq W_*$ and $X_*\subseteq W'_*$
that are related by a smooth map $p:W_*'\to W_*$, meaning that
$\calo_{W_\alpha,p(x)}\subseteq \calo_{W'_\alpha,x}$ is smooth for
every $x\in W_*'$.  If $\calr_*$ is a free $\calo_{W_*}$-algebra
forming a resolvent then we can take a free $p^*\calr_*$ algebra, say,
$\calr'_*$ as a resolvent on $ W_*$.  Now a projective 
$\calr_*$-resolution
$\calp_*$ of $\calf_*$ gives a projective $\calr_*'$-resolution
$\calp_*'=\calp_*\otimes_{\calr_*}\calr'_*$, and \ref{prop:base change
Atiyah} applies again.

Finally, the independence from the choice of locally finite coverings
by Stein compact sets is easily seen considering refinements; we leave
the simple details to the reader.
\end{proof}

We now translate the earlier results on the naturality of Atiyah 
classes.  
The following is an almost immediate consequence of \ref{2.14}.

\begin{prop}\label{3.14}
For every morphism of complexes $\alpha: \calf\to\calg$ of
degree $0$ in $D^-_{coh}(X)$ the diagram
\begin{diagram}[h=8mm]
\calf & \rTo^{\alpha} & \calg\\
\dTo^{\At^k(\calf)} &&\dTo_{\At^k(\calg)}\\
\calf\dotimes \Lambda^k\bbbl_{X/Y} & \rTo^{\alpha\otimes \id} &
\calg\dotimes \Lambda^k\bbbl_{X/Y}
\end{diagram}
commutes.
\end{prop}

\begin{proof}
Consider as in \ref{2.17} a resolvent of $X$ over $Y$ and choose
projective approximations $\calp_*\to \calf_*$ and $\calq_*\to
\calg_*$.  There is a morphism
$\tilde\alpha_*:\calp_*\to\calq_*$ lifting the given morphism
$\alpha$ and the assertion follows now easily from \ref{2.14}.
\end{proof}

For $\calf\in D^-_{coh}(X)$ as before,
$\Ext^\sbullet_X(\calf,\calf):=\bigoplus_{i}\Ext^i_X(\calf, \calf)$
carries a natural algebra structure given by Yoneda product, and
$\Ext^\sbullet_X(\calf, \calf\dotimes\Lambda^k \bbbl_{X/Y})$ is a
bimodule over $\Ext^\sbullet_X(\calf, \calf)$.

\begin{prop}
    \label{2.20}
The power $\At^k(\calf)$ of the Atiyah class of $\calf$ 
is a (graded) central element of degree $k$ in the bimodule 
$\Ext^\sbullet_X(\calf, \calf\dotimes\Lambda^k\bbbl_{X/Y})$,
which means that
$$
\xi\cdot \At^k(\calf)=(-1)^{ik}\At^k(\calf)\cdot\xi
$$
for every element $\xi\in \Ext^i_X(\calf, \calf)$.
\end{prop}

\begin{proof}
By the preceding result and by \ref{2.19} the diagram
\begin{diagram}[s=8mm]
\calf & \rTo^\xi & \calf[i] \\
\dTo^{\At^k(\calf)} && \dTo_{(-1)^{ik}\At^k(\calf)[i]}\\
\calf\dotimes\Lambda^k \bbbl_{X/Y} & \rTo^{\xi\otimes\id} &
\calf[i]\dotimes \Lambda^k\bbbl_{X/Y}
\end{diagram}
commutes.
\end{proof}

Using the isomorphisms $\bbbs^{k}(\bbbl_{X/Y}[1])\xto{\cong}(\Lambda^k
\bbbl_{X/Y})[k]$, this compatibility with morphisms can be summarized
in analogy to \ref{thm:Atiyah} as follows.

\begin{cor}
The powers of the Atiyah class define morphisms of exact
functors
$$
\At^{k}(\ \,): \id \xto{\quad} (\ \,)\dotimes
\bbbs^{k}(\bbbl_{X/Y}[1])\,,\quad k \ge 0\,,
$$
on $D_{coh}^{-}(X)$ that commute with the shift functor.\qed
\end{cor}

The Atiyah classes are as well compatible with mappings in the
following sense.

\begin{prop}\label{2.21}
Let
\begin{diagram}[h=7mm]
X' & \rTo^f & X\\
\dTo && \dTo\\
Y' & \rTo & Y
\end{diagram}
be a diagram of complex spaces and $\calf\in D^-_{coh}(X)$. 
Under the natural map
$$
\Ext^k_X(\calf, \calf\dotimes \Lambda^k\bbbl_{X/Y})\lto
\Ext^k_{X'}(Lf^*(\calf), Lf^*(\calf)\dotimes
\Lambda^k\bbbl_{X'/Y'})
$$
the Atiyah class $\At^k(\calf)$ is mapped onto the Atiyah class
$\At^k(Lf^*(\calf))$.
\end{prop}

Before giving the proof we have to choose compatible resolvents
for $X$ and $X'$.  This is done as follows. Consider locally
finite coverings $\{X_j\}_{j\in J}$,
$\{X_i'\}_{i\in I}$ of $X$ resp.\ $X'$ by Stein compact sets
such that there is a map $\sigma:I\to J$ with
$f(X_i')\subseteq X_{\sigma(i)}$ for all $i\in I$.  We may
assume that $I=J$ and $\sigma=\id_I$ since otherwise we can
replace $I$ and $J$ by the disjoint union $K=I\cup J$ and
consider the coverings $\{X_k\}_{k\in K}$
resp.\ $\{X_k'\}_{k\in K}$, where we set $X_i:=X_{\sigma(i)}$
for $i\in I$ and $X_j':=\emptyset$ for $j\in J$.

Let $X_*=(X_\alpha)_{\alpha\in A}$ and
$X'_*=(X'_\alpha)_{\alpha\in A}$ be the associated simplicial
spaces. There is a natural functor
$f^*:\Coh(X_*)\to \Coh (X_*')$ which associates to
$\calm_*\in \Coh(X_*)$ the module with
$f^*(\calm_*)_\alpha:=(f|X_\alpha)^* (\calm_\alpha)$. It is easy to
see that $f^*$ transforms projective modules into
projective modules (for instance, using \ref{2.4} it is
sufficient to check this for modules of type
$p^*(\calp_\alpha)$).

Choose embeddings $X_i\hto L_i$, $X_i'\hto V_i$ into Stein
compact sets in $\bbbc^{n_i}\times Y$ resp.\ $\bbbc^{n_i}\times
Y'$ and take the diagonal embedding $X_i'\subseteq
L_i':=L_i\times_Y V_i$. From the data
$$
X_i\subseteq L_i\quad \mbox{and}\quad  X_i'\subseteq L_i'
$$
we construct smoothings $X_*\subseteq W_*$, $X_*'\subseteq
W_*'$ of $X\to Y$ and $X'\to Y'$, respectively, as explained
in \ref{2.2}. The projections induce a system of compatible maps
$\tilde f:W_\alpha'\to W_\alpha$ restricting to $f$ on
$X_\alpha'$. As above there is a natural functor $\tilde
f^*:\Coh W_*\to\Coh W_*'$ transforming projective modules into
projective modules.  Therefore, if  $\calr_*\to \calo_{X_*}$ is
a resolvent then
$\tilde f^*(\calr_*)$ is a free DG algebra over $W_*'$. The
projection $\tilde f^*(\calr_*)\to \calo_{X_*'}$ can be
factored through a quasiisomorphism $\calr_*'\to\calo_{X_*'}$
such that $\calr_*'$ is a free DG $\tilde f^*(\calr_*)$-algebra,
and we take this as a resolvent for $X'$ over $Y'$.

We remark that one has a natural functor
$\tilde f^{-1} :\mathbf {Ab}(W_*)\to \mathbf{Ab}(W_*')$ on
the category of simplicial systems of abelian groups on $W_*$
with $\tilde f^{-1}(\cala)_\alpha:=\tilde
f^{-1}(\cala_\alpha)$.

After these preparation we can easily deduce \ref{2.21}.

\begin{proof}[Proof of \ref{2.21}]
Let $\calf_*\to\tilde\calf_*$ be a quasiisomorphism into a
$W_*$-acyclic module, $\calp_*\to\tilde\calf_*$ a projective
approximation, and let $\nabla:\calp_*\to \calp_* \otimes
\Omega^1_{\calr_*/Y}$ be a connection.
Using the product rule, the composed map
\begin{diagram}
\tilde f^{-1}\calp_*
&\rTo^{\tilde f^{-1}(\nabla)}&
\tilde f^{-1}\calp_* \otimes_{\tilde f^{-1}\calr_*} \tilde
f^{-1}(\Omega^1_{\calr_*/Y})
&\rInto &
\tilde f^{-1}\calp_*\otimes_{\tilde
f^{-1}\calr_*}\Omega^1_{\calr_*'/Y'}
\end{diagram}
can be extended to a connection $\nabla'$ on $\calp_*':=
\tilde f^{-1}\calp_*\otimes_{\tilde f^{-1}\calr_*} \calr_*'
$. Hence under the natural map
$$
\Ext^k_{\calr_*}(\calp_*,
\calp_*\otimes \Omega^k_{\calr_*/Y})
\lto
\Ext^k_{\calr_*'}(\calp_*',
\calp_*'\otimes\Omega^k_{\calr_*'/Y'} )
$$
the Atiyah class $\At^k(\calp_*)$ maps onto
$\At^k(\calp_*')$. Since the module on the left
is isomorphic to
$\Ext^k_X(\calf,\calf\dotimes\Lambda^k\bbbl_{X/Y})$ and the
module on the right is isomorphic to
$\Ext^k_{X'}(Lf^*\calf, Lf^*\calf\dotimes
\Lambda^k\bbbl_{X'/Y'})$, the result follows.
\end{proof}

\begin{cor}\label{2.22}
Let
\begin{diagram}[h=7mm]
X' & \rTo^f & X\\
\dTo && \dTo\\
Y' & \rTo & Y
\end{diagram}
be a diagram of complex spaces and $\alpha:
Lf^*\calf\to\calf'$ a morphism of complexes of sheaves, where
$\calf\in D^-_{coh}(X)$ and
$\calf'\in D_{coh}^-(X')$. Then the diagram
\begin{diagram}[h=7mm]
Lf^*\calf & \rTo^{\alpha} & \calf'\\
\dTo^{Lf^*\At^k(\calf)} &&\dTo_{\At^k(\calf')}\\
Lf^*\calf\dotimes Lf^*(\Lambda^k\bbbl_{X/Y}) &
\rTo^{\alpha\otimes can} &
\calf'\dotimes
\Lambda^k\bbbl_{X'/Y'}
\end{diagram}
commutes.
\end{cor}

\begin{proof}
This is an immediate consequence of \ref{3.14} and the fact that
$Lf^*\At^k(\calf)=\At^k(Lf^*\calf)$.
\end{proof}

In the final result of this section we explain how to compute
the Atiyah class using the second fundamental form. Slightly
more generally it is convenient to show the following
result.

\begin{prop}\label{2.23}
Let $X$ be smooth over $Y$ and let
$$
0\to \calf'\stackrel{j}{\lto}
\calf\stackrel{p}{\lto} \calf''\to 0
$$
be an exact sequence of coherent sheaves
on $X$. Assume that there is a  map
$\nabla:\calf\to\calf''\otimes \Omega_{X/Y}$ satisfying the
product rule $\nabla(fa)=\nabla(f)a+p(f)\otimes da$
for local sections $f$ in $\calf$ and $a$ in $\calox$.
The linear map
$$
\sigma:=\nabla\circ j:\calf'\lto
\calf''\otimes \Omega^1_{X/Y}.
$$
is then $\calox$-linear, and if
$$
\begin{array}{c}
\delta': \Hom(\calf',\calf'' \otimes \Omega^1_{X/Y})\to
\Ext^1(\calf',  \calf' \otimes\Omega^1_{X/Y})\\
\delta'': \Hom(\calf', \calf'' \otimes\Omega^1_{X/Y})\to
\Ext^1(\calf'',  \calf''\otimes\Omega^1_{X/Y})
\end{array}
$$
denote the boundary operator in the respective long exact 
$Ext$-sequences then 
$$
\delta'(\sigma)=\At(\calf')\quad\mbox{and}\quad
\delta''(\sigma)=-\At(\calf'').
$$
\end{prop}

If $\calf$ itself admits a connection, say,
$\nabla_1:\calf\to\calf\otimes\Omega^1_{X/Y}$, then
this result applies to $\nabla:= p\otimes1\circ\nabla_1$, and
the map $\sigma$ becomes the usual second fundamental form. 
This shows thus in particular how to compute the Atiyah class 
using the second fundamental form.

\begin{proof}
Let $(X_*,W_*, \calr_*)$ be a resolvent of $X$
as in \ref{2.17}. There are projective
resolutions $\calp'_*$, $\calp_*\cong\calp_*'\oplus
\calp_*''$, $\calp''_*$ of $\calf'_*$, $\calf_*$, $\calf''_*$,
respectively that fit into a commutative diagram
\begin{diagram}[h=7mm]
0 & \rTo & \calp'_* &\rTo^{inj\quad} &
\calp_*\cong\calp_*'\oplus
\calp_*''  & \rTo^{\quad \tilde p} &  \calp''_* &\rTo & 0\\
&& \dTo>{\pi'} && \dTo>\pi && \dTo>{\pi''}\\
0 & \rTo & \calf'_* & \rTo^j & \calf_*
 & \rTo^p & \calf''_* & \rTo & 0
\end{diagram}
Consider connections $\nabla'$, $\nabla''$ on $\calp_*'$,
$\calp_*''$, respectively, and equip $\calp_*$ with the
connection $\nabla'\oplus\nabla''$. Using the isomorphisms
$$
\begin{array}{rcl}
\Hom(\calf', \calf'' \otimes\Omega^1_{X/Y})
&\cong &
\Hom(\calp_*', \calf_*'' \otimes\Omega^1_{\calr_*/Y})\\
\Ext^1(\calf'',  \calf''\otimes\Omega^1_{X/Y})
& \cong &
H^1(\Hom(\calp_*'',  \calf_*''\otimes\Omega^1_{\calr_*/Y})),
\end{array}
$$
$\delta''(\sigma)$
can be computed as follows. The linear map $\sigma\circ\pi'$
can be extended to a linear map $\tilde \sigma:
\calp_*\to\calf'' \otimes\Omega^1_{\calr_*/Y}$, and $\tilde
\sigma\circ\partial$ is zero on $\calp_*'$ and so defines a
map $\sigma'':\calp_*''\to\calf_*''
\otimes\Omega^1_{\calr_*/Y}$ that represents $-\delta''
(\sigma)$. Taking as extension the map $\tilde\sigma:= \nabla
\circ
\pi-(\pi''\otimes 1)\circ\nabla''\circ\tilde p$ we get
$\tilde\sigma\circ\partial=-(\pi''\otimes 1)\circ\nabla''
\partial\circ\tilde p $, and so $(\pi''\otimes 1)\circ\nabla''
\partial$ represents $\delta'' (\sigma)$. By construction it
also represents $-\At(\calf'')$, proving the second part of
the result. The equality  $\delta'(\sigma)=\At(\calf')$
follows with a similar argument and is left to the reader.
\end{proof}

\begin{rems}
1. If $X$ is smooth over $Y$ and if $\calf$ is a coherent
$\calox$--module then the Atiyah class of
$\At(\calf)$ is the cohomology class that is represented
by the extension
$$
0\to\calf \otimes_\calox\Omega^1_{X/Y}  \stackrel{j}{\lto}
\calp^{1}(\calf):=p_1^*(\calf)
\otimes_{\calo_{X\times_Y X}} \calo_{X\times_Y X}/\calj^2
\stackrel{p}{\lto} \calf\to 0\,,\leqno (*)
$$
where $\calj\subseteq\calo_{X\times_Y X}$ is the ideal of the
diagonal and where we consider $\calp^1(\calf)$ as an
$\calox$--module via the second projection $p_2$.
Using \ref{2.23} a simple proof can be given as follows. The map
$\nabla':p_1^*(\calf)\to
\calf\otimes_\calox \Omega^1_{X/Y}$
given by $\nabla(f\otimes a):= f\otimes da$ for local sections
$f$ of $\calf$, $a$ of $\calox$, is easily seen to factor
through a map $\nabla:\calp^{1}(\calf)\to
\calf\otimes_\calox \Omega^1_{X/Y}$ that satisfies the product
rule required in \ref{2.23}. The map $\sigma:=\nabla\circ
j$ is just the identity on $\calf\otimes_\calox
\Omega^1_X$ as $\Omega^1_{X/Y}$ is identified with
$\calj/\calj^2$ via $da\mapsto 1\otimes
a-a\otimes 1$. Hence \ref{2.23} shows that
$\At(\calf)$ is represented by $-\delta''(\sigma)$, where
$\delta'':\Hom(\calf \otimes\Omega^1_{X/Y}, \calf
\otimes\Omega^1_{X/Y})\to
\Ext^1(\calf,  \calf\otimes\Omega^1_{X/Y})$ is the boundary
map. In view of \cite[X.126, Cor.1(b)]{Alg} this proves
the remark.

2. In order to check that the sign of our Atiyah classes is correct,
consider the case that $X=\bbbp^n$ and $\calm=\calo_{\bbbp^n}(1)$. It
is well known that the extension class of the Euler sequence
$$\textstyle
0\lto \Omega^1_{\bbbp^n}(1)\stackrel{j}{\lto}
\calo_{\bbbp^n}^{n+1}=\bigoplus_{i=0}^n e_i \calo_{\bbbp^n}
\stackrel{p}{\lto} \calo_{\bbbp^n}(1) \lto 0
$$
in $\Ext^1_{\bbbp^n}(\calo_{\bbbp^n}(1),\Omega^1_{\bbbp^n}(1))\cong
H^1(\bbbp^n,\Omega^1_{\bbbp^n})\cong \bbbc$ represents the first Chern
class $c_1(\calo_{\bbbp^n}(1))$.  Explicitly, the maps $j,p$ are
given by $j(x_i dx_j-x_jdx_i)=e_jx_i-e_ix_j$ and $p(e_i)=x_i$ 
respectively, with $x_0,\ldots,x_n$ the homogeneous coordinates 
of $\bbbp^n$. The module
$\calo_{\bbbp^n}^{n+1}$ admits a unique connection with
$\nabla(e_i)=0$, and the corresponding
map $\sigma$ in \ref{2.23} is easily seen to be $-\id$ on
$\Omega^1_{\bbbp^n}(1)$.  Hence it follows from \ref{2.23} that
$c_1(\calo_{\bbbp^n}(1))=- \At(\calo_{\bbbp^n}(1))$.
\end{rems}

Summing up we have a naturally defined class
$$
\exp(-\At(\calf))=\sum_k (-1)^k\At^k(\calf)/k!\in \prod_k
\Ext^k_X(\calf, \calf\dotimes\Lambda^k\bbbl_{X/Y})\,,
$$
the {\em Atiyah-Chern character\/}, for every complex
$\calf$ in $D^-_{coh}(X)$.

\section{The semiregularity map}

\subsection*{A semiregularity map for modules}
Recall that a complex $\calf$ over a complex space $X$ is called {\em
perfect}\/ if it admits locally a quasiisomorphism to a bounded
complex of free coherent $\calox$--modules. For instance, every
$\calox$--module on a complex manifold when
considered as a complex concentrated in degree 0  is perfect.
If $\calg\in D(X)$ is a complex then for every perfect complex
$\calf$ there is a  natural trace map
$$
\Tr: \Ext_{X}^k(\calf, \calf\dotimes\calg)\lto
H^k(X,\calg),\quad k\ge 0,
$$
see \cite{Ill}.
These maps are compatible with taking cup products: if
$\xi\in\Ext_{X}^j(\calg,\calg')$, where $\calg'$ is another complex
in $D(X)$, then the diagram
\begin{diagram}[h=7mm]
\Ext_{X}^k(\calf,\calf\dotimes\calg)
&\rTo^\Tr  &
H^k(X,\calg)\\
\dTo>\xi && \dTo>\xi\\
\Ext_{X}^{k+j}(\calf,\calf\dotimes\calg')
&\rTo^\Tr  &
H^{k+j}(X,\calg')
\end{diagram}
commutes.
For instance, applying the trace map to the Atiyah classes we
obtain for every perfect complex $\calf$ well defined classes
$$
\ch_k( \calf):=\Tr((-1)^k\At^k(\calf))/k!\in
H^k(X,\Lambda^k\bbbl_{X/Y}),
$$
that are the components of the {\em Chern character}
$\ch(\calf)=\Tr \exp(-\At(\calf))$  of $\calf$. If
$X$ is a manifold and $\calf$ is a vector bundle then this gives
the usual Chern character of $X$, see \cite{At,OTT}, and in the 
general algebraic case it is Illusie's \cite{Ill} description.

\begin{defn}\label{3.1}
Let $X\to Y$ be a morphism of complex spaces and let
$\calf$ be a perfect complex of $\calox$--modules. The map
$$
\sigma:=\Tr( *\cdot \exp(-\At(\calf))):
\Ext_{X}^2(\calf,\calf)\lto
\prod_k H^{k+2}(X, \Lambda^k\bbbl_{X/Y})
$$
is called the {\em semiregularity map} for $\calf$.
\end{defn}

Slightly more generally, for every coherent
$\calox$--module $\caln$ and every $r\ge 0$ there are maps
$$
\sigma=\sigma_{\caln}:=\Tr( *\cdot \exp(-\At (\calf))):
\Ext_{X}^r(\calf,\calf\dotimes\caln)\lto
\prod_k H^{k+r}(X, \caln\dotimes\Lambda^k\bbbl_{X/Y}),
$$
to which we will also refer as the semiregularity map.

To formulate the next proposition, set $\bbbl:=\bbbl_{X/Y}$ 
and note that the group
$A:=\bigoplus_{i}A^{i}$ with
$$
A^{i}:= \bigoplus_{j}\Ext^{i+j}_X(\calf ,
\calf\dotimes\Lambda^j\bbbl)
\cong
\bigoplus_{j}\Ext^i_X(\calf ,
\calf\dotimes\bbbs^j(\bbbl[1]))
$$
carries a natural algebra structure that is associative but in
general not graded commutative. Moreover $M:=\bigoplus_i M^i$ with
$M^i:=\bigoplus_{j}\Ext^i_X(\calf , \calf\dotimes
\caln\dotimes\bbbs^j(\bbbl[1]))$ is a graded $A$-bimodule.
Every element $\xi\in T^{r-1}_{X/Y}(\caln)\cong
\Ext^{r}_{X}(\bbbl[1], \caln)$ defines a derivation
$\langle\xi,*\rangle:A\to M$ of degree $r$ which is
induced by the composition
\begin{diagram}
\Lambda^\sbullet\bbbl &\rTo^{\Delta} & \bbbl\dotimes
\Lambda^\sbullet\bbbl &\rTo^{\xi\otimes \id} & \caln\dotimes
\Lambda^\sbullet\bbbl\,,
\end{diagram}
where $\Delta$ is the indicated component of the 
comultiplication on $\Lambda^\sbullet\bbbl$.
Thus, for elements $\omega_{1}\in A^i$, $\omega_{2}\in A$
we have
$$
\langle\xi,\omega_{1}\omega_{2}\rangle=
\langle\xi,\omega_{1} \rangle\omega_{2} + (-1)^{ir}
\omega_{1}\langle\xi,\omega_{2}\rangle.
$$
In particular,
for $\calf=\calox$, this gives a derivation
$\langle\xi,*\rangle$ from the cohomology algebra
$\bigoplus H^{i+j}(X,\Lambda^{j}\bbbl)$ into $\bigoplus
H^{i+j}(X,\caln\dotimes\Lambda^{j}\bbbl)$. Note that using
resolvents one can check all this very explicitly.

\begin{prop}  
    \label{3.2}
Let $X\to Y$ be a morphism of complex spaces. If $\calf$ is a perfect
complex on $X$ and $\caln$ is an $\calox$--module then the diagram
\begin{diagram}[w=7mm,h=9mm]
T^{r-1}_{X/Y}(\caln)&& \rTo^{\langle *, -\At(\calf)\rangle} &&
\Ext_{X}^r(\calf, \calf\dotimes\caln)\\
& \rdTo(3,2)<{\langle *,\ch_{k+1}(\calf)\rangle }&& %
\ldTo(3,2)>{\sigma_{k}}\\ %
&&H^{k+r}(X,\caln\dotimes\Lambda^k\bbbl_{X/Y})
\end{diagram}
commutes, where $\sigma_{k}$ denotes the $k^{th}$ component of
$\sigma_{\caln}$.
\end{prop}

\begin{proof}
We need to show for $\xi \in T^{r-1}_{X/Y}(\caln)$ that
$$
\Tr\big(\langle\xi, \At(\calf)\rangle \cdot
\frac{\At^k(\calf)}{k!}\big)
= \big\langle {\xi}\,,\, \Tr\big( \frac{\At^{k+1}(\calf)}{
(k+1)!}\big)\big\rangle\ .
$$
As the trace map is compatible with taking cup
products, the diagram
\begin{diagram}[h=7mm]
\Ext_{X}^{k+1}(\calf,  \calf\dotimes
\Lambda^{k+1}\bbbl) &\rTo^{\Delta} &
\Ext_{X}^{k+1}(\calf, \calf\dotimes\bbbl\dotimes \Lambda^k\bbbl
)
&\rTo^{\xi} &
\Ext_{X}^{k+r}(\calf, \calf\dotimes\caln\dotimes\Lambda^k\bbbl
)\\
\dTo>\Tr && \dTo>\Tr  && \dTo>\Tr \\
H^{k+1}(X,\Lambda^{k+1}\bbbl)
&\rTo^{\Delta} &
H^{k+1}(X,\bbbl\dotimes \Lambda^k\bbbl)
&\rTo^\xi &
H^{k+r}(X, \caln\dotimes\Lambda^k\bbbl)
\end{diagram}
commutes, where as before $\bbbl=\bbbl_{X/Y}$.  Therefore
$$
\Tr \big\langle \xi,\frac{\At^{k+1}(\calf)}{(k+1)!} \big\rangle =
\big\langle \xi, \Tr
\big(\frac{\At^{k+1}(\calf)}{(k+1)!}\big)\big\rangle
\ .
\leqno (1)
$$
As $\langle\xi, *\rangle$ is a derivation of degree $r$ on
$\Ext^\sbullet_{X}(\calf,
\calf\dotimes \bbbs^\sbullet (\bbbl[1]))$ and 
$\At^k(\calf)\in A^0 \cong \Ext^0_X(\calf,
\calf\dotimes \bbbs^\sbullet (\bbbl[1]))$, 
we obtain that
$$
\langle \xi, \At^{k+1}(\calf) \rangle =\sum
\At^i(\calf)\langle \xi, \At(\calf) \rangle
\At^{k-i}(\calf).
$$
For homogeneous endomorphisms $f, g$ the trace satisfies
$\Tr(fg)=(-1)^{|f||g|}\Tr(gf)$, whence
taking traces yields
$$
\Tr\big \langle \xi, \frac{\At^{k+1}(\calf)}{(k+1)!}\big \rangle =
\Tr \left( \langle
\xi, \At(\calf)\rangle\cdot \frac{\At^k(\calf)}{k!} \right).
$$
Comparing with (1), the result follows.
\end{proof}

In case that $X$ is smooth and $Y$ is a reduced point the
theorem above specializes to the following corollary.

\begin{cor} \label{3.3}
For every complex manifold $X$ and every coherent
$\calox$--module $\calf$ there is a commutative diagram
\begin{diagram}[w=7mm,h=8mm,midshaft]
H^{r-1}(X,\Theta_X)
&& \rTo^{\langle *, -\At(\calf)\rangle} && \Ext_{X}^r(\calf,\calf)\\
&\rdTo<{\langle *, \ch_{k+1}(\calf)\rangle}&&
\ldTo>{\sigma_{k}}\\
&& H^{k+r}(X,\Omega^k_X)\,.
\end{diagram}
\end{cor}

In case of compact algebraic manifolds the map $\langle *,
\ch_{k+1}(\calf)\rangle$ on $H^{1}(X,\Theta_X)$ has the
following geometric interpretation; see \cite[4.2]{Blo},
or, for a more general statement, \ref{A.7}, \ref{A.8}.
Given an infinitesimal deformation of $X$ represented 
by a class $\xi\in H^1(X,\Theta_X)$, the unique 
horizontal lift of $\ch_{k+1}(\calf)$ relative to 
the Gau\ss-Manin connection stays
of Hodge type $(k+1, k+1)$ if and only if $\langle \xi ,
\ch_{k+1}(\calf)\rangle = 0$.  In the next result we will
show that $\langle \xi, -\At(\calf)\rangle$ gives the
obstruction for deforming $\calf$ in the
direction of $\xi$.  Thus the semiregularity map $\sigma_k$
relates the obstruction to deform $\calf$
along $\xi$ with the obstruction that its Chern class
$\ch_{k+1}(\calf)$ stays of pure Hodge type along $\xi$.

\begin{prop}\label{3.4}
Let $f:X\to Y$ be a morphism of complex spaces and let 
$X \subseteq X'$ be a $Y$--extension of $X$ by a  
coherent $\calox$--module $\caln$ so
that $X'$ defines a  class $[X']\in T^1_{X/Y}(\caln)$.  If
$\calf$ is a coherent $\calox$--module then under the composed map
\begin{diagram}
\ob: T^1_{X/Y} (\caln)
& \rTo^{\langle *,-\At(\calf)\rangle} &
\Ext_{X}^2(\calf, \calf\dotimes\caln ) & \rTo^{can} &
\Ext_{X}^2(\calf, \calf\otimes\caln )
\end{diagram}
one has $\ob ( [X'])=0$ if and only if there is an
$\calo_{X'}$--module $\calf'$ extending $\calf$ to $X'$, 
i.e.\ there is an exact sequence of $\calo_{X'}$--modules
$$
0\to\calf \otimes\caln\to\calf'\to \calf\to 0\,.
$$
\end{prop}

\begin{proof}
In the algebraic case this is shown in \cite[IV.3.1.8]{Ill}. In
the analytic case we can proceed as follows. Let $(X_*,W_*, \calr_*)$
be a resolvent of $X$ over $Y$ as in \ref{2.17} and let
$\gamma:\calp_*\to\calf_{*}$ be a projective resolution of $\calf_*$
as an $\calr_*$--module. A $Y$-extension $[X']$ of $X$ gives rise to
an extension $[X'_*]$ of $X_*$ by $\caln_*$. Since $W_{*}$ is smooth
over $Y$, the embedding $X_*\hto W_*$ can be lifted to a $Y$-map
$X'_*\hto W_*$, and the surjection of algebras $\pi:\calr_*\to
\calo_{X_*}$ to a map of $\calo_{W_*}$-algebras
$\pi':\calr_*\to\calo_{X'_*}$. With $\partial$ the differential on
$\calr_*$, the map $\xi:=-\pi'\partial:\calr_*\to \caln_*$ is a
$Y$-derivation of degree 1 that represents the class of $[X']$ in
$$
T^1_{X/Y}(\caln) \cong H^1(\Der_Y(\calr_*,\caln_*))\cong
\Ext^1_{\calr_*}(\Omega^1_{\calr_*/Y}, \caln_*).
$$
If one equips the trivial extension $\calr_*[\caln_*]$ 
with the differential $(r,n)\mapsto (\partial(r),\xi(r))$,
the map $\pi'+\id_{\caln_*}:\calr_*[\caln_*]\to
\calo_{X'_*}$ becomes a quasiisomorphism of DG algebras that
restricts to the identity on $\caln_*$.

Let now $\nabla:\calp_*\to\calp_*
\otimes\Omega^1_{\calr_*/Y}$ be a connection.  Contracting
with $\xi$ and projecting onto $\calf_{*}$ gives a map
$\nabla_\xi:\calp_*\to \calf_*\otimes\caln_*$ of degree 1
satisfying the product rule
$\nabla_\xi(pr)=\nabla_\xi(p)r+(-1)^{|p|}\gamma(p)\otimes \xi(r)$
for local sections $r$ in $\calr_*$ and $p$ in $\calp_*$ .  The
class $\ob(\calf)$ is  represented by the map
$$
-(\gamma\otimes \xi)\circ [\partial,\nabla]= [\partial,\nabla_\xi
]=\nabla_{\xi}\partial:\calp_*\lto \calf_*\otimes\caln_*
$$
of degree 2; note that by \ref{2.16}, \ref{2.11}(\ref{2.11.4}) and 
\ref{2.22.a}
$$
\Ext_{X}^2(\calf,\calf\otimes\caln)\cong
H^2(\Hom_{\calr_*}(\calp_*,\calf_*\otimes\caln_*))\,.\leqno (*)
$$
If the class $\ob(\calf)$ vanishes then
$[\partial,\nabla_\xi]=[\partial,h]$ for some $\calr_*$-linear
map $h:\calp_*\to\calf_*\otimes\caln_*$ of degree 1. The 
differential
$$
\left(p, f\otimes n\right)\mapsto
\left(\partial(p),(\nabla_\xi-h)(p) \right)
$$
defines then on $\calp_*\oplus \calf_*\otimes\caln_*$ the structure of
a DG module over $\calr_*[\caln_*]$ and we denote this DG module
$\calp_*[\calf_*\otimes\caln_*]$. Using the exact cohomology sequence
associated to the exact sequence of DG modules
$$
0\to \calf_*\otimes\caln_*\to \calp_*[\calf_*\otimes\caln_*]
\to \calp_*\to 0\,,
$$
it follows that the $\calo_{X'_*}$ module $\calh^0(
\calp_*[\calf_*\otimes\caln_*])$ is an extension of $\calf_*$ by
$\calf_*\otimes\caln_*$. Gluing yields an $\calo_{X'}$ modules
$\calf'$ which is an extension of $\calf$ by $\calf\otimes\caln$.

Conversely assume that there exists such an extension $\calf'$ of
$\calf$ by $\calf\otimes\caln$. As $\calp_*$ is projective the map
$\gamma:\calp_*\to \calf_*$ can be lifted to a map of
$\calr_*$--modules $\gamma':\calp_*\to \calf'_*$. A simple calculation
shows that $\delta:=-\gamma'\partial:\calp_*\to \calf_*\otimes\caln_*$
satisfies the product rule
$\delta(pr)=\delta(p)r+(-1)^{|p|}\gamma(p)\otimes \xi(r)$ for sections
$r$ in $\calr_*$ and $p$ in $\calp_*$. Hence $h=\nabla_\xi-\delta$ is
$\calr_*$-linear and satisfies $[\partial,h]=[\partial,\nabla_\xi]$,
whence the cohomology class of $[\partial,\nabla_\xi]$ in
$H^2(\Hom(\calp_*,\calf_*\otimes\caln_*))$ vanishes. As this class
represents $\ob(\calf)$ under the isomorphism $(*)$, the result
follows.
\end{proof}

Later on we will apply the semiregularity map to modules on the total
space of a deformation of a complex space. In order to verify that
such a module has locally finite projective dimension, the following
standard criterion is useful.

\begin{prop}
    \label{3.4.a}
Let $f:X\to S$ be a flat morphism of complex spaces and $\calf$ a
complex in $D_{coh}^-(X)$. If the restriction
$\calf\dotimes_{\calox}\calo_{X_{s}}$ to every fibre
$X_{s}:=f^{-1}(s)$ is a perfect complex on $X_{s}$, then $\calf$ is a
perfect complex on $X$.
\end{prop}

This is an immediately consequence of the following simple
lemma from commutative algebra.

\begin{lem}\label{7.21}
Let $A\to B$ be a flat morphism of local noetherian rings and let
$F^\sbullet$ be a complex of $B$--modules with finite cohomology that
is bounded above. If $F^\sbullet\dotimes_A A/\fm_A$ is a perfect
complex of $B/\fm_AB$--modules then $F^\sbullet$ is a perfect complex
of $B$--modules.
\end{lem}

\begin{proof}
For the convenience of the reader we include the simple argument. We
may assume that $F^\sbullet$ is a complex of finite free $B$--modules
with $F^i=0$ for $i\gg 0$. By assumption $F^\sbullet\otimes A/\fm_A$
is a perfect complex and so for $k\ll 0$ the complex
$$
F^\sbullet_{(k)}:\qquad\quad
\ldots\to  F^{k-1}
\stackrel{\partial}{\lto}
F^k \to F^k/\partial F^{k-1} \to 0
$$
has the property that $F^\sbullet_{(k)}\otimes A/\fm_A$ is exact with
$(F^k/\partial F^{k-1})\otimes A/\fm_A$ a free $B/\fm_AB$--module.
Using induction on $n$ and the long exact cohomology sequences 
associated
to the exact sequences of complexes
$$
0\lto F^\sbullet_{(k)}\otimes \fm_A^n/\fm_A^{n+1} \lto
F^\sbullet_{(k)}\otimes A/\fm_A^{n+1}\lto
F^\sbullet_{(k)}\otimes A/\fm_A^n \lto 0
$$
it follows that
$F^\sbullet_{(k)}\otimes A/\fm_A^n$ is exact and that
$(F^k/\partial F^{k-1})\otimes A/\fm_A^n$ is a free
$B/\fm_A^nB$--module. Hence $F^\sbullet_{(k)}$ is exact and
$F^k/\partial F^{k-1}$ is free as a $B$--module, proving the
lemma.
\end{proof}

\begin{rems}\label{3.5}
1. We note that the construction of the semiregularity map is
compatible with morphisms. More precisely, given a diagram of complex
spaces as in \ref{2.21} and a perfect complex $\calf$ on $X$, for any
coherent $\calox$--module $\caln$ the diagram
\begin{diagram}[h=7mm,midshaft]
\Ext^r_X(\calf,\calf\dotimes\caln)
& \rTo^\sigma &
\textstyle\prod_k H^{k+r}(X,
\caln\dotimes \Lambda^k\bbbl_{X/Y})\\
\dTo && \dTo \\
\Ext^r_{X'}(Lf^*\calf,Lf^*\calf\dotimes Lf^*\caln)
& \rTo^\sigma &
\textstyle\prod_k H^{k+r}(X',
Lf^*\caln\dotimes \Lambda^k\bbbl_{X'/Y'})
\end{diagram}
commutes. This follows from \ref{2.21} and the fact that the
trace map is compatible with taking inverse images.

2. If in \ref{3.1} the support of $\calf$ is contained
in a closed subset, say, $Z$ of $X$ then the semiregularity
map admits a factorization
$$
\Ext_{X}^r(\calf,\calf\dotimes\caln)\stackrel{\sigma_Z}{\lto}
\prod_k H_Z^{k+r}(X, \caln\dotimes\Lambda^k\bbbl_{X/Y})
\stackrel{can}{\lto}
\prod_k H^{k+r}(X, \caln\dotimes\Lambda^k\bbbl_{X/Y}).
$$
This follows as the trace map factors by \cite{Ill} through the local
cohomology.
\end{rems}
That the powers of the Atiyah classes are graded central elements by 
\ref{2.20} allows the following glimpse at the relevance of the 
semiregularity map for deformation problems. Let
$$
[\ ,\ ]:\Ext_{X}^i(\calf,\calf)\times
\Ext_{X}^j(\calf,\calf)\xto{\quad} 
\Ext_{X}^{i+j}(\calf,\calf)\quad,\quad
[\xi,\zeta] := \xi\zeta - (-1)^{ij}\zeta\xi\,,
$$
denote the graded Lie algebra structure underlying the Yoneda product 
on the $\Ext$--algebra. Centrality of the Atiyah classes together 
with the fact that the trace vanishes on commutators implies then the 
following result.

\begin{cor}
    \label{cor:4.8}
The family of semiregularity maps
$$
\sigma: \Ext_{X}^{r}(\calf,\calf)\lto \prod_{k}
H^{k+r}(X,\Lambda^k\bbbl_{X/Y})\quad,\quad r\ge 0\,,
$$
vanishes on 
$$
[\Ext_{X}^{\sbullet}(\calf,\calf),
\Ext_{X}^{\sbullet}(\calf,\calf)]\subseteq
\Ext_{X}^{\sbullet}(\calf,\calf)\,.
$$
\qed
\end{cor}
As is well known, and will be recalled in Section 6 below, the vector
space $\Ext_{X}^{1}(\calf,\calf)$ is the tangent space to a
semi\-universal deformation of $\calf$, if it is finite dimensional.
The obstructions to lift such tangent directions to second order lie
in $[\Ext_{X}^{1}(\calf,\calf), \Ext_{X}^{1}(\calf,\calf)]\subseteq
\Ext_{X}^{2}(\calf,\calf)$, and these obstructions are thus
annihilated by the semiregularity map.

\subsection*{A semiregularity map for subspaces}
Let $X$ be a complex space, $Z\subseteq X$ a closed complex
subspace and $\caln$ a coherent $\calox$--module.  In this
part we will show how to define a semiregularity map on
$T^2_{Z/X}(\caloz\dotimes_{\calox} \caln)$ provided that
$\caloz$ has locally  finite projective dimension as an
$\calox$--module.  In particular, this will give a
generalization of Bloch's semiregularity map to arbitrary
subspaces of manifolds.  The idea is to define first a map
from $T^2_{Z/X}(\caloz\dotimes \caln)$ into
$\Ext_{X}^2(\caloz,\caloz\dotimes \caln)$ and then to
compose this with the semiregularity map for $\calf=\caloz$
as defined in the previous part.  The key technical lemma is
as follows.

\begin{lem}\label{5.8}
Let  $Z\subseteq X$ be a closed embedding of complex spaces.

{\em 1.} For each complex of $\caloz$--modules $\calm$ there are 
natural maps 
$$
\epsilon^{(k)}:T^k_{Z/X}(\calm)\lto \Ext^k_X(\caloz,\calm)\,,\quad 
k\in\bbbz\,.
$$
In case $\calm=\caloz$, the map
$\epsilon^{(\sbullet)}:T^\sbullet_{Z/X}(\caloz)\to
\Ext^\sbullet_{X}(\caloz,\caloz)$ is a morphism of graded Lie
algebras.

{\em 2.}  Let $X\to Y$ be a morphism of complex spaces and let
$\caln$  be a coherent $\calox$--module. With
$T^{k-1}_{X/Y}(\caloz\dotimes \caln)\to T^k_{Z/X}(\caloz\dotimes
\caln)$ the boundary map in the long tangent cohomology
sequence associated to the triple $Z\to X\to Y$, the composition
$$
T^{k-1}_{X/Y}(\caln)\stackrel{can}{\lto}
T^{k-1}_{X/Y}(\caloz\dotimes \caln)\xto{}
    T^k_{Z/X}(\caloz\dotimes \caln) \xto{\epsilon^{(k)}}
\Ext^k_X(\caloz, \caloz\dotimes \caln)
$$
is given by $\xi\mapsto\langle\xi, -(-1)^{k}\At(\caloz)\rangle$.
\end{lem}

\begin{proof}
Let $(X_*, W_*, \calr_*)$ be a resolvent for $X$ over $Y$ as in
\ref{2.17} and choose a quasiisomorphism
$\calm_{*}\to\widetilde\calm_{*}$ into a $W_{*}$-acyclic complex of
$\calo_{Z_{*}}$--modules, where $\calm_{*}$ and $\calo_{Z_{*}}$ denote
the simplicial sheaves on $W_{*}$ associated to $\caloz$ and $\calm$,
respectively. We choose an algebra resolution 
$\calr_{*}\xto{i}\cals_*\to
\calo_{Z_*}$ of the composition $\calr_{*}\to\calo_{X_{*}}\to
\calo_{Z_*}$ so that $\cals_{*}$ is a graded free algebra over
$\calr_*$ and $\cals_*\to \calo_{Z_*}$ is a quasiisomorphism of DG
algebras. In particular, $\cals_*$ is a projective approximation of
$\calo_{Z_*}$ as an $\calr_*$--module, and
$\cals_{*}\to\cals_{*}\otimes_{\calr_{*}}\calo_{X_{*}}$ is a
quasiisomorphism by \ref{2.11}. It follows that
$\cals_{*}\otimes_{\calr_{*}}\calo_{X_{*}}$ provides a projective
resolution of $\calo_{Z_{*}}$ as $\calo_{X_{*}}$--module. Hence
\begin{align}
    \Ext^i_X(\caloz, \calm) &\cong
    H^i(\Hom_{\calr_*}(\cals_*,\widetilde\calm_*))
    \quad \quad\text{and}\tag{$*$}\\
    T^i_{Z/X}(\calm) &\cong 
    H^i (\Hom_{\calr_*}(\Omega^1_{\cals_*/\calr_{*}}, \widetilde 
\calm_*))
    \notag
\end{align}
for all $i$.
Composing the natural inclusions, see \ref{2.34},
$$
\Hom_{\calr_*}(\Omega^1_{\cals_*/\calr_{*}}, \widetilde\calm_*)
\hto \Der_{\calr_*}(\cals_*, \widetilde\calm_*)\hto
\Hom_{\calr_*}(\cals_*,\widetilde \calm_*)\leqno (**)
$$
gives the desired map in (1). If $\calm=\caloz$, then
$T^\sbullet_{Z/X}(\caloz)\cong H^\sbullet( \Der_{\calr_*}(\cals_*,
\cals_*))$ and $\Ext^\sbullet_X(\caloz, \calm)\cong
H^\sbullet(\End_{\calr_*}(\cals_*))$, and the inclusion of derivations
into endomorphisms is a morphism of DG Lie algebras that induces a
morphism of graded Lie algebras in cohomology.

To show (2), note first that $\calo_{Z_{*}}\dotimes\caln_{*}$ is
represented by $\cals_{*}\otimes_{\calr_{*}}\caln_{*}$. Consider a
derivation $\delta\in \Der_Y(\calr_*,\caln_*)$ of degree $k-1$ that
represents the cohomology class $\xi$ in $T^{k-1}_{X/Y}(\caln)$. Its
image in $T^{k-1}_{X/Y}(\caloz\dotimes \caln)$ is then represented by
$1\otimes \delta\in \Der_Y(\calr_*,\cals_{*} \otimes\caln_*)$. Under
the isomorphism in $(*)$, the image of the latter element in
$T^k_{Z/X}(\caloz\dotimes \caln)$ is represented by a
$\calr_*$-derivation $[\partial, \tilde\delta]$, where
$\tilde\delta:\cals_* \to \cals_{*}\otimes \caln_*$ is a derivation
restricting to $1\otimes\delta$ on $\calr_{*}$. Let now
$\nabla:\cals_*\to \cals_*\otimes\Omega^1_{\calr_*/Y}$ be a connection
on $\cals_*$, as exists by \ref{2.12}. If $h$ is one of the maps
$(1\otimes \delta)\circ \nabla$ or $\tilde \delta$ from $ \cals_*$ to
$\cals_*\otimes\caln_*$ then the product rule
$$
h(sr)=h(s)r+(1)^{|s|(k-1)} s\delta(r)
$$
is satisfied for local sections $s$ in $\cals_*$ and $r$ in $\calr_*$.
It follows that the difference $\tilde\delta- (1\otimes \delta)\circ
\nabla$ is $\calr_*$-linear and so $[\partial,\tilde\delta]$ and
$[\partial,(1\otimes \delta)\circ \nabla]$ represent the same
cohomology class in $H^k(\Hom_{\calr_*}(\cals_*,\cals_{*}\otimes
\caln_*))$. As 
$$\langle [\delta],\At(\caloz)\rangle = 
(1\otimes \delta)\circ[\partial,\nabla]
=(-1)^{k-1}[\partial,(1\otimes \delta)\circ \nabla]\,,
$$
by definition, and as $\delta$ is of degree $k-1$, we have
$$
\xi \mapsto [\partial,\tilde\delta] =
(-1)^{k-1}\langle [\delta],\At(\caloz)\rangle = \langle \xi, 
-(-1)^{k}\At(\caloz)\rangle\,,
$$
as required.
\end{proof}

\begin{defn}\label{5.9}
Let $X\to Y$ be a morphism of complex spaces and let
$Z\subseteq X$ be a closed complex subspace of $X$ such that
$\caloz$ has locally finite projective dimension over
$\calox$.  The composition of the canonical map
$T^2_{Z/X}(\caloz) \to \Ext^2_{X}(\caloz,\caloz)$ of
\ref{5.8} (1) with the semiregularity map $\sigma$ defined
in \ref{3.1} yields a map
$$
\tau:T^2_{Z/X}(\caloz)
\lto \prod_{k\ge 0}
H^{k+2}(X,\Lambda^k\bbbl_{X/Y}),
$$
which we call the {\em semiregularity map\/} for $Z$.
\end{defn}

Again we have such a semiregularity map more generally for
any coherent $\calox$--module $\caln$  and any $r\ge 0$,
$$
\tau_\caln:= \sigma_{\caln}\circ \epsilon^{r}:
T^r_{Z/X}(\caloz\dotimes\caln) \lto \prod_{k\ge 0} H^{k+r}(X,\caln
\dotimes_\calox\Lambda^k\bbbl_{X/Y})\,.
$$
Note that this gives in particular a semiregularity map for
every closed subspace $Z$ of a complex manifold $X$.
Combining \ref{3.2} and \ref{5.8} we obtain the following
result.

\begin{prop}\label{5.9.a}
With the notations and assumptions as in {\em \ref{5.9}\/}, the
diagram
$$
\begin{diagram}[w=7mm,h=8mm,midshaft]
T^{r-1}_{X/Y}(\caln)
&& \rTo^{can} &&  T^r_{Z/X}(\caloz\dotimes\caln)\\
&\rdTo<{\langle *, \ch_k(\caloz)\rangle}&&
\ldTo>{(-1)^{r}\tau_{k-1}}\\
&&   H^{k+r-1}(X,\caln\dotimes\Lambda^{k-1}\bbbl_{X/Y})\,,
\end{diagram}
%
$$
commutes, where $\tau_{k-1}$ is the
$(k-1)^{st}$ component of the  semiregularity map $\tau_\caln$.
\end{prop}

For instance, for $Y$ a reduced point, $X$ smooth, and $Z$ a closed
subscheme of codimension $k$ in $X$, we can rephrase \ref{5.9.a} as
follows.

\begin{cor} \label{5.10}
For a complex manifold $X$ and a
subspace $Z$ of codimension $k$ the diagram
$$
\begin{diagram}[w=7mm,h=8mm,midshaft]
H^{r-1}(X,\Theta_X)
&& \rTo^{can} && T^r_{Z/X}(Z)\\
&\rdTo<{\langle*,[Z]\rangle}&&
\ldTo>{(-1)^{r}\tau_{k-1}}\\
&&   H^{k+r-1}(X,\Omega^{k-1}_X)\,,
\end{diagram}
$$
commutes, where $[Z]$ denotes the fundamental class of $Z$ in
$H^k(X,\Omega^k_X)$.
\end{cor}

\begin{proof}
For a manifold, $T^{r-1}_{X/Y}(\calox)\cong H^{r-1}(X,\Theta_X)$, and
according to Grothendieck, \cite[4(16)]{Gro}, see also
\cite[2.18]{Mur}, the fundamental class of $Z$ is given by
$$
[Z]=\frac{(-1)^{k-1}}{(k-1)!} c_k(\caloz) =\ch_k(\caloz)\,,
$$
whence the result follows from \ref{5.9.a}.
\end{proof}

Let $X\to Y$ be a morphism of complex spaces and let
$Z\subseteq X$ be a closed complex subspace of $X$.
In the final result of this section we consider for a
coherent $\calox$--module $\caln$ the boundary map
$\delta$ in the long exact tangent cohomology sequence
$$
\cdots \lto T^1_{Z/Y}(\caloz\otimes \caln)\lto
T^1_{X/Y}(\caloz\otimes \caln)\stackrel{\delta}{\lto}
T^2_{Z/X}(\caloz\otimes \caln)\to\cdots
$$
for the triple $Z\to X\to Y$.  The exactness of this sequence
gives immediately the following interpretation of the
composed map
$$
\gamma: T^1_{X/Y}(\caln)\stackrel{can}{\lto}
T^1_{X/Y}(\caloz\otimes \caln)\stackrel{\delta}{\lto}
T^2_{Z/X}(\caloz\otimes \caln)
$$
in terms of extensions.

\begin{lem}\label{5.7.a}
Let $X'$ be an extension of $X$ by $\caln$. The image of the class
$[X']\in T^1_{X/Y}(\caln)$ under $\gamma$ vanishes, $\gamma([X'])=0$,
if and only if there is an extension $Z'$ of $Z$ by $\caloz\otimes
\caln$ that fits into a commutative diagram
\begin{diagram}[s=7mm]
0 & \rTo & \caln & \rTo & \calo_{X'} & \rTo& \calox &\rTo 
&0\hphantom{\,.}\\
&& \dTo^{can} && \dTo && \dTo\\
0 & \rTo & \caloz\otimes\caln & \rTo & \calo_{Z'}
& \rTo & \caloz &\rTo& 0\,.
\end{diagram}
\qed
\end{lem}

\begin{rems}\label{3.5a}
1. The construction of the canonical map $\epsilon^{(\sbullet)}$ in
\ref{5.8} (1) is compatible with morphisms of complex spaces. More
precisely, assume given a diagram of complex spaces as in \ref{2.21} 
and
a closed subspace $Z$ of $X$. With $Z':=f^{-1}(Z)\subseteq X'$, the 
diagram
\begin{diagram}[h=7mm,midshaft]
T^r_{Z/X}(\calm)
& \rTo^{\epsilon^{(r)}} &
\Ext_X^r(\caloz,\calm)\\
\dTo && \dTo \\
T^r_{Z'/X'}(\calo_{Z'}\dotimes_{\caloz}\calm)
& \rTo^{\epsilon^{(r)}} &
\Ext^r_{X'}(\calo_{Z'},\calo_{Z'}\dotimes_{\caloz}\calm)
\end{diagram}
commutes for every coherent $\caloz$--module $\calm$.

2. In analogy with \ref{3.5}, the construction of the semiregularity
map $\tau$ is also compatible with morphisms of complex spaces.
This follows from \ref{3.5}, using the preceding remark. We
leave the straightforward formulation and its proof to the reader.

3. It follows from \ref{3.5}(2) that the semiregularity
map $\tau$ factors through local cohomology,
$$
T^r_{Z/X}(\caloz\dotimes \caln)\xto{\tau_Z}
\prod_k H_Z^{k+r}(X, \caln\dotimes\Lambda^k\bbbl_{X/Y})
\xto{\ can\ }
\prod_k H^{k+r}(X, \caln\dotimes\Lambda^k\bbbl_{X/Y})\,.
$$
\end{rems}

Finally, we wish to point out that \ref{cor:4.8} carries over as well
to the family of semiregularity maps for subspaces, the comparison map
$\epsilon^{(\sbullet)}:T^\sbullet_{Z/X}(\caloz)\to
\Ext^\sbullet_{X}(\caloz,\caloz)$ being a morphism of graded Lie
algebras.

\section{Applications to the variational Hodge conjecture}

Let $X$ be a compact complex algebraic manifold so that its
cohomology admits a Hodge decomposition
$$
H^k(X,\bbbc)\cong \bigoplus_{p+q=k}H^q(X,\Omega_X^p)\,.
$$
Recall that a cohomology class
$$
\alpha \in H^{p,p}(X ,\bbbq):= H^p(X,\Omega_X^p)\cap H^p(X,\bbbq)
$$
is called {\em algebraic} if an appropriate multiple
$k\alpha , k \in \bbbn$, is represented by an algebraic
cycle $Z$ of codimension $p$ in $X$, in the sense that $k\alpha = 
[Z]$, the cohomology class of the cycle.
The famous and so far unsolved Hodge conjecture asks wether
every class in $H^{p,p}(X,\bbbq)$ is algebraic.

In \cite {Gro}, Grothendieck proposed the following weaker version
that is called the {\em variational Hodge conjecture}. Let $\pi:X\to
S$ be a deformation of a compact algebraic manifold $X_0=\pi^{-1} (0)$
over a smooth germ $(S,0)$. The local system $R^{2p} f_*(\bbbc)\otimes
\calos$ carries then the natural {\em Gau\ss-Manin connection\/}.
Assume that $\alpha$ is a horizontal section of $R^p
f_*(\Omega^p_{X/S})$ in the sense that $\alpha$ can be lifted locally
to a horizontal section in $R^{2p} f_*(\Omega^{\ge p}_{X/S})\subseteq
R^{2p} f_*(\bbbc)\otimes \calos$, see also
\ref{A.4}. 

The variational Hodge conjecture asks now: If the restriction of
$\alpha$ to the special fibre, $\alpha(0) \in H^p
(X_0,\Omega_{X_0}^p)$, is algebraic, is then $\alpha (s)\in
H^p(X_s,\Omega^p_{X_s})$ algebraic for all $s\in S$ near $0$, where
$X_s:=\pi^{-1}(s)$?

In this section we will give an affirmative answer to this
problem if the class $\alpha(0)$ is the $p^{th}$ component
$\ch_p(\cale_0)$ of the Chern character of some coherent
sheaf $\cale_0$ for which the $p^{th}$ component of the 
semiregularity 
map,
$$
\Ext^2_{X_{0}} (\cale_0, \cale_0)\longrightarrow
H^{p+1}(X_0,\Omega _{X_0}^{p-1})\,,
$$
is injective. We will then call $\cale_0$ in brief a {\em
$p$-semiregular sheaf}. Slightly more generally, with $n$ the
dimension of $X_0$ it is convenient to introduce for any subset $I
\subseteq \{{0,\ldots n}\}$ the following notion: $\cale_0$ is called
$I$-{\em semiregular} if the part of the semiregularity map
$$
\sigma_I: \Ext_{X_{0}}^2 (\cale_0, \cale_0)\longrightarrow
\prod_{p\in I} H^{p+1}(X_0,\Omega _{X_0}^{p-1})
$$
is injective. The main result of this section is the
following theorem.

\begin{theorem}\label{A.1}
Let $\pi:X\to S$ be a deformation of a compact complex
algebraic manifold $X_0$ over a smooth germ $S=(S,0)$ and
set $X_s:= \pi^{-1}(s)$ for $s\in S$.  Assume that
$(\alpha_p)_ {p \in I}$ is a horizontal section in
$\prod_{p\in I} R^p \pi_*(\Omega^p_{X/S})$.  If there is an
$I$-semiregular sheaf $\cale_0$ on $X_0$ with
$$
\alpha_p(0)=\ch_p(\cale_0)\in H^p(X_0, \Omega^p_{X_0})\,,\quad
p\in I,
$$
then $\alpha_p(s)\in
H^p(X_s,\Omega^p_{X_s})$ is algebraic for all
$s\in S$ near $0$ and each $p\in I$.
\end{theorem}

In analogy with the notion of a $I$-semiregular sheaf on
$X_0$, a complex subspace $Z_0\subseteq X_0$ will be called
$I$-{\em semiregular} if the part of the semiregularity map
$$
\tau_I:T^2_{Z_0/X_0}(\calo_{Z_0})\lto \prod_{p \in
I}H^{p+1}(X_0, \Omega^{p-1}_{X_0})
$$
is injective.  We will also derive the following variant of
\ref {A.1} that generalizes a result of S.~Bloch
\cite{Blo}.

\begin{theorem}\label{A.2}
Let $\pi : X \to S$ and $(\alpha _p)_{p \in I} $ be as in {\em 
\ref{A.1}\/}.
If there is an I-semiregular subspace $Z_0\subseteq X_0$ with
$\alpha_p(0)=\ch_p(\calo_{Z_0})$ for $p\in I$ then $\alpha_p(s)$
is  algebraic for all $s \in S$ near $0$ and each $p\in I$.
\end{theorem}

For the proof of these results we need a few preparations.
\begin{lem}\label{A.3}
Let $X \subseteq X'$ be an extension of a complex space $X$
by a coherent $\calo_X$--module $\calm$ and assume that the
morphism $\xi : \bbbl_{X/X'} \to\calm$ of degree $1$ in the
derived category $D(X)$ represents the class
$$
[X'] \in T^1_{X/X'}(\calm)\cong \Ext^1_X (\bbbl_{X/X'}, \calm)\,.
$$
If\/  $\bbbl_{X/X'}\stackrel{+1}{\lto}\bbbl_{X'} \dotimes\calo_X$
denotes the canonical map of degree 1, then the diagram
\begin{diagram}[s=7mm,midshaft]
\bbbl_{X/X'}  &\rTo^{+1}  & \bbbl_{X'} \dotimes
\calo_X
\\
 \dTo<{+1}>\xi &&  \dTo>{can} \\
\calm  &  \rTo^d   &\Omega^1_{X'} \otimes \calo_X
\end{diagram}
commutes in $D (X)$.
\end{lem}

\begin{proof}
Let $(X'_*, W'_*, \calr'_*)$ be a resolvent of $X'$ and choose a free
graded DG $\calr'_*$--algebra $\calr_*$ that provides a DG algebra
resolution $p:\calr_*\to \calo_{X_{*}}$ of the composition 
$\calr'_*\xto{p'}
\calo_{X'_{*}}\to \calo_{X_{*}}$. By \ref{2.11}, the induced map
$\calr_*\otimes_{\calr'_*}\calo_{X'_*}\to \calo_{X_*}$ is a
quasiisomorphism and so $(X_*, W_*:=X'_*,
\calr_*\otimes_{\calr'_*}\calo_{X'_*})$ constitutes a resolvent for
$X$ over $X'$. Thus
$$
T^1_{X/X'}(\calm)\cong
H^1(\Der_{\calo_{X'_*}}(\calr_*\otimes_{\calr'_*}\calo_{X'_*},
\calm_*)) \cong
H^1(\Der_{\calr'_*}(\calr_*,\calm_*))\,,
\leqno (*)
$$
and moreover
$$
\bbbl_{X/X'}\cong
C^\sbullet(\Omega^1_{\calr_*\otimes_{\calr'_*}\calo_{X'_*}
/\calo_{X'_*}} \otimes\calo_{X_*}) \cong
C^\sbullet(\Omega^1_{\calr_*/\calr'_*} \otimes\calo_{X_*})\,,
$$
where $C^\sbullet$ is the \v{C}ech functor as in \ref{Cech functor}.
By construction, $\calr_{*}$ is a free $\calr'_*$-algebra and so
admits an augmentation, say, $\beta:\calr_* \to \calr'_*$. Clearly,
the $\calr'$-linear map $\beta$ will not be a morphism of DG algebras
in general, but the commutative diagram
$$
\begin{diagram}[small]
\calr_*  &\rDotsto^\beta     & \calr'_*&\rTo^{p'} 
&\calo_{X'}\hphantom{\,.}\\
 &&\dTo&&\dTo\\
 &&\calr_*   &\rTo^p    &\calox\,.    
\end{diagram}
$$
shows that the $\calr'_*$-derivation $-p'\beta\partial : \calr_*\to
\calo_{X'}$ factors through a derivation $\xi:\calr_*\to\calm_*$ that
in turn represents $[X']\in T^1_{X/X'}(\calm)$ under the isomorphism
$(*)$, see \cite[(3.13)]{Fle1}. The map
$\bbbl_{X/X'}\xto{+1}\bbbl_{X'}\otimes\calo_X$ is induced
by the map of complexes
$$
[\partial, d \beta]:\Omega^1_{\calr_*/\calr'_*}\lto
\Omega^1_{\calr_*'}\otimes \calr_*\,,
$$
and composition with the canonical map $\Omega^1_{\calr'_*}
\otimes\calr_*\to \Omega^1_{X_*'} \otimes \calo_{X_*}$ yields the map
$-d p' \circ d \beta \circ \partial$. As this map coincides with the
$\calr'_{*}$-derivation $d\circ \xi: \calr_{*}\to \calm_{*}\xto{d}
\Omega^1_{X_*'} \otimes \calo_{X_*}$ under the identification
$$
\Hom_{\calr_*}
(\Omega^ 1_{\calr_*/\calr'_*}, \Omega^1_{X_*'} \otimes 
\calo_{X_*})\cong 
\Der_{\calr'_*}(\calr_*, \Omega^1_{X_*'} \otimes \calo_{X_*})\,,
$$
the result follows.
\end{proof}

\begin{sit}\label{A.4}
Let $\pi : X \to S$ be a deformation of a K\"ahler manifold $X_0$ over
an artinian germ $S=(S,0)$. The canonical exact sequence
$$
0\lto \calox\otimes_{\calo_S}\Omega^1_S\lto\Omega^1_X \lto
\Omega^1_{X/S}\lto 0
$$
induces a morphism $\Omega^1_{X/S}\xto{+1} \calox\otimes\Omega^1_S$ of
degree $1$ in the derived category $D(X)$ and as well, by taking 
exterior
powers, morphisms
$$
\nabla_{X/S}: \Omega^p_{X/S}\xto{+1}
\Omega^{p-1}_{X/S}\otimes_{\calo_S}\Omega^1_S\,,\quad p\ge 1,
$$
that we call the {\em Gau\ss-Manin connections\/}. These maps induce
the classical Gau\ss-Manin connections
$$
\nabla_{X/S}:H^q(X, \Omega^p_{X/S})\lto H
^{q+1}(X,\Omega^{p-1}_{X/S}) \otimes_{\calo_S}\Omega^1_S
.$$
Recall that the spectral sequence
$$
E^{pq}_1=H^q(X,\Omega^p_{X/S})\Rightarrow H^{p+q} (X,
\Omega^\sbullet_{X/S})\cong H^{p+q} (X_0, \bbbc)\otimes_\bbbc \calo_S
$$
degenerates and that, by Griffiths' transversality theorem, the
canonical connection $\nabla$ on
$H^{p+q} (X, \Omega^\bullet_{X/S})$ satisfies
$$
\nabla (H^q (X, \Omega^{\ge p}_{X/S}))\subseteq H^{q+1}(X,\Omega^{\ge
p-1}_{X/S})
\otimes_{\calo_S}\Omega^1_S
$$
and induces just the map $\nabla_{X/S}$ above (see
\cite{Gri2}).  We will call a class $\alpha$ in
$H^q(X,\Omega^p_{X/S})$ in brief {\em horizontal} if it can
be lifted to a horizonal section in
$$
H^{p+q}(X, \Omega^{\ge p}_{X/S})\subseteq
H^{p+q}(X_0,\bbbc)\otimes_\bbbc \calo_S.
$$
In particular, the following well known result holds.
\end{sit}

\begin{lem}\label{A.5}
A horizontal class $\alpha\in H^q(X, \Omega^p_{X/S})$ satisfies
$\nabla_{X/S}(\alpha) = 0$. \qed
\end{lem}

\begin{sit}\label{A.6}
Let $\pi:X\to S$ be as in \ref{A.4} and assume that $S\hto S'$ is an
extension of the artinian germ $S$ by a coherent
$\calo_S$--module $\caln$ such that
there is a smooth morphism $\pi' : X' \lto S'$ restricting to $\pi$
over $S$. {\em We will suppose henceforth that the map\/}
$$
d : \caln \lto \Omega^1_{S'} \otimes_{\calo_{S'}} \calo_S
\leqno (*)
$$
{\em is injective\/}. Let
$$
\nabla' : =1_{\calos}\otimes_{\calo_{S'}} \nabla_{X'/S' }:
\Omega^q_{X/S} \lto \Omega^{q-1}_{X/S}
\otimes_{\calo_{S'}}\Omega^1_{S'}
$$
be the map induced by the Gau\ss- Manin connection for $X'\to S'$
and denote by the same symbol the induced map in cohomology,
$$
\nabla': H^p (X, \Omega^q_{X/S})\lto H^{p+1}
(X,\Omega^{q-1}_{X/S}) \otimes_{\calo_{S'}}\Omega^1_{S'}\,.
$$
\end{sit}

\begin{lem} \label{A.7}
{\em (\cite {Blo})\/} With notation and assumptions as in {\em
\ref{A.6}\/}, a horizontal class $\alpha\in H^q(X, \Omega^p_{X/S})$ 
can be lifted to a horizontal section in $H^q
(X', \Omega^p_{X'/S'})$ if and only if $\nabla' (\alpha)= 0$ in
$H^{q+1}(X, \Omega^{p-1}_{X/S})\otimes_{\calo_{S'}}\Omega^1_{S'}$. 
\qed
\end{lem}

Now return to the notation and assumption in \ref {A.4}. The extension
$X'$ of $X$ by $\caln _X : =\calo_X \otimes_\calos\caln$ gives a class
$\xi \in T^1_X (\caln_X)$ or, equivalently, a morphism $\xi : \bbbl_X
\lto \caln_X$ of degree 1 in the derived category $D(X)$. Taking
exterior powers yields a map, denoted by the same symbol,
$$
\xi : \Lambda^p \bbbl_X \lto \Lambda^{p-1}\bbbl_X\dotimes\caln_X
\xto{can}\Omega^{p-1}_{X/S} \otimes \caln_X\,.
$$
The next result describes how this map relates to the map
$\nabla'$ introduced in \ref{A.6}.

\begin {lem}\label{A.8}
The diagram
\begin{diagram}[s=7mm,midshaft]
\Lambda^p\bbbl_X &\rTo^{can} &
\Omega^p_{X/S} \\
\dTo>{\xi}<{+1}&& \dTo>{\nabla'}<{+1}\\
\Omega^{p-1}_{X/S} \otimes\caln &\rTo^{1\otimes d} &
\Omega^{p-1}_{X/S} \otimes_{\calo_{S'}}\Omega^1_{S'}
\end{diagram}
in $D(X)$ is commutative.
\end{lem}

\begin{proof}
For suitable representative of the cotangent complexes involved there
is a commutative diagram of exact sequences of complexes of 
$\calox$--modules
\begin{diagram}[s=7mm,midshaft]
0  &\rTo & \bbbl_{X/X'}[-1] & \rTo &
\bbbl_{X'}\dotimes \calo_X & \rTo &
\bbbl_X &\rTo  & 0\hphantom{\,.}\\
&&\dTo^\xi && \dTo && \dTo\\
 0 &\rTo & \caln_X & \rTo^d
& \Omega^1_{X'}\otimes\calo_X & \rTo &\Omega^1_X &\rTo &0\,,
\end{diagram}
by \ref{A.3} and our assumption that $d$ is injective on $\caln$. In
$D(X)$ this gives rise to a commutative diagram
\begin{equation}
    \begin{diagram}[s=7mm] \bbbl_X &
    \rTo & \bbbl_{X/X'}\hphantom{\,.}\\
    \dTo && \dTo>\xi \\
    \Omega^1_X &\rTo  & \caln_X[1]\,.
    \end{diagram}
   \tag{$\dagger$}
\end{equation}
As well, the commutative diagram of exact sequences of
$\calox$--modules
\begin{diagram}[s=7mm,midshaft]
0&\rTo &\caln_X &\rTo^{1\otimes d}
&\Omega^1_{X'}\otimes\calo_X &\rTo &\Omega^1_X &
\rTo& 0\\
&&\dTo<{1\otimes d} && \dEqual &&\dTo\\
0 &\rTo & \calo_X \otimes_{\calo_{S'}} \Omega^1_{S'}&\rTo
&\Omega^1_{X'}\otimes
\calo_X  &\rTo & \Omega^1_{X/S} &\rTo & 0
\end{diagram}
yields in $D(X)$ a commutative diagram 
\begin{diagram}[s=7mm,midshaft]
\Omega^1_X & \rTo & \caln_X[1]\hphantom{\,.} \\
\dTo && \dTo>{1 \otimes d}\\
 \Omega^1_{X/S}\cong\Omega^1_{X'/S'}\otimes \calox
& \rTo & \calox\otimes_{\calo_{S'}}\Omega^1_{S'}[1]\,.
\end{diagram}
Combining this square with the one in ($\dagger$) gives the result 
for $p=1$.
The general case follows from this by taking appropriate exterior
powers.
\end{proof}

The key observation in proving \ref{A.1} is the following
proposition.  We keep the notation and assumptions as
introduced in \ref{A.4} and \ref{A.6}.

\begin{prop}\label{A.9}
Assume $\cale$ is a coherent $S$--flat sheaf on $X$ such
that $\cale_0 :=\cale|X_0$ is $I$-semiregular, where $X_0
= \pi^{-1}(0)$.  The following conditions are then equivalent.

{\em 1.}\ The sheaf $\cale$ can be extended to a
deformation $\cale'$ on $X'$ over $S'$.

{\em 2.}\ The partial Chern character $\ch_I
(\cale):=(\ch_{p}(\cale))_{p\in I}\in \prod_{p \in I} H^p (X,
\Omega^p_{X/S})$ can be lifted to a horizontal section in $\prod_{p
\in I} H^p(X', \Omega^p_{X'/S'})$.
\end{prop}

\begin{proof}
Let $\xi \in T^1_{X}(\caln_X)$ be the class corresponding to the
extension $X\hto X'$, given by a morphism $\xi :
\bbbl_{X/X'}\to\caln_X$ in the derived category. The result \ref {A.8}
induces a commutative diagram
\begin{diagram}[w=6mm,h=9mm,midshaft]
{\prod_{p \in I}} H^p(X,\Lambda^p \bbbl_{X})&\rTo
&{\prod\limits_{p \in I}} H^p(X, \Omega^p_{X/S})  \cong
\prod_{p \in I} H^p(X, \Omega^p_{X'/S'})
\otimes_{\calo_{S'}}\calo_S \\
\dTo_{\xi}&& \dTo_{\nabla'}\\
{\prod_{p \in I}}  H^{p+1}(X,
\Omega^{p-1}_{X/S})\otimes_{\calo_S} \caln&
\rInto^{1\otimes d}&
{\prod_{p \in I}} H^{p+1}(X, \Omega^{p-1}_{X/S})
\otimes_{\calo_{S'}}\Omega^1_{S'}\ ,
\end{diagram}
where $\nabla'$ is the Gau\ss-Manin connection.  

The relative partial Chern character
$\ch_I (\cale) \in\prod_{p\in I} H^p(X,\Omega^p_{X/S})$
is the image of the absolute partial Chern character $\ch'_I(\cale):=
(\Tr((-\At(\cale))^p/p!))_{p\in I}$ in $\prod_{p\in I}
H^p(X,\Lambda^p \bbbl_{X})$.  Using \ref {A.7} and the
injectivity of $1 \otimes d$ in the diagram above, it follows
that (2) is equivalent to

3.\ Contracting against $\xi$ sends the absolute partial Chern 
character to 
zero, 
$$
\langle\xi ,\ch'_I (\cale)\rangle  = 0\,.
$$
Now, \ref{3.2} yields a commutative diagram
\begin{diagram}[w=6mm,h=10mm]
T^1_X(\caln_X)&& \rTo^{\ob:=\langle *,-\At(\caln_X)\rangle}
&& \Ext^2_X(\cale, \cale
\otimes_{\calo_S}\caln)\\
&\rdTo(3,2)_{\langle * , \ch'_I(\cale)\rangle}  &&
\ldTo(3,2)_{\sigma_{I,\caln}}\\ 
&&{\prod\limits_{p\in I}}
H^{p+1}(X, \Omega^{p-1} _{X/S}) \otimes_{\calo_S}\caln
\end{diagram}
and $\sigma_{I,\caln}$ is injective by \ref{A.10} below. It follows 
that
$\ob(\xi) = 0$ if and only if (3) holds. On the other hand, by
\ref{3.4}, the condition $\ob(\xi) = 0$ is equivalent to (1), so the
result follows.
\end{proof}

\begin{lem}\label{A.10}
With notation as above, if $\cale_0$ is I-semiregular
then the map $\sigma_{I,\caln}$ is injective.
\end{lem}
 
The {\em proof }follows by a simple induction on the length
of $\caln$. The initial step for $\caln\cong \bbbc$ is equivalent to 
the 
assumption as $\sigma_{I,\bbbc}$ is the corresponding partial 
semiregularity map for $\cale_0$ in view of the isomorphisms
\begin{align*}
    \Ext^2_X(\cale, \cale\otimes_{\calo_S}\bbbc)&\cong 
    \Ext^2_{X_{0}}(\cale_{0}, \cale_{0})\quad,\quad\text{as $\cale$ 
is 
    $S$--flat,}\\
    H^{p+1}(X, \Omega^{p-1} _{X/S}) \otimes_{\calo_S}\bbbc &\cong
    H^{p+1}(X_{0}, \Omega^{p-1} _{X_{0}}) \quad,\quad\text{by base 
change.}
\end{align*}
The induction step is left to the reader.
\qed

\begin{proof}[Proof of Theorem {\em \ref{A.1}\/}]
Let $S_{n}$ be the $n^{th}$ infinitesimal neighbourhood  of  $0$ in
$S$ so that $S_n\hto S_{n+1}$ is an extension of $S_n$
by $\caln_n = \fm^{n+1}/\fm^{n+2}$, where $\fm\subseteq\calo_{S,
0}$ is the maximal ideal. The map
$$
d: \caln_n\lto \Omega^1_{S_{n+1}} \otimes \calo_{S_n}
$$
is injective and applying \ref {A.9} repeatedly we see that $\cale_0$
can be lifted to a deformation $\cale_n$ on $X_n$ for all $n$. Let
$\cale$ be a versal deformation of $\cale_0$ which is a coherent
module on $X\times_ST, $ where $(T,0)$ is a complex space germ over
$(S,0)$. Using versality there are $(S,0)$--morphisms
$\varphi_n:(S_n,0)\to(T,0)$ with $(1\times \varphi_n)^*(\cale)\cong
\cale_n$ and $\varphi_{n+1}|S_n=\varphi_n$. Hence $(T,0)\to (S,0)$
admits a formal section, namely
$\bar\varphi:=\lim\limits_{\longrightarrow}\varphi_n: (\hat S,0)\to
(T,0)$. By Artin's approximation theorem, we can find a convergent
section $\varphi:(S,0) \to (T,0)$. Now $\calf:= (1\times
\varphi)^*(\cale)$ is a coherent S-flat module on $X$ that induces
$\cale _0$ on the special fibre. The uniqueness of the horizontal
lifting gives that $\alpha_p= \ch_p(\calf)$ as sections in $R^p
f_*(\Omega^p_{X/S})$ for each $p\in I$. Hence $\alpha_p(s)=
\ch_p(\calf|X_{s})$ is algebraic for all $s \in S$ near $0$ and each
$p\in I$.
\end{proof}

The {\em proof} of \ref{A.2} is similar.  The sole
difference is that in order to derive the analogue of
\ref{A.9}, one has to use \ref{5.9.a} and \ref{5.7.a} instead of
\ref{3.2} and \ref{3.4}.  We leave the technical details to the
reader.
\qed

\section{Deformation theory}

Before formulating the main results we review some basic notation and
facts about deformation theories. In contrast to \cite{Sch} we will
not use the language of deformation functors but instead employ
deformation groupoids as in \cite{Rim, Fle2, BFl}. Most deformations
will take place over $\Ans$, $\Arts$, or $\Formals$, the categories of
germs of complex spaces, Artinian complex spaces, or formal complex
spaces respectively, over a fixed germ $\Sigma=(\Sigma,0)$.

\begin{sit}
Let $p\colon \bE\to \bB$ be a functor between categories. We will
denote objects of $\bB$ by capital letters whereas the objects of
$\bE$ will be written in lower case. To indicate that $a$ is an object
of $\bE$ over $S\in \bB$ we write simply $a\mapsto S$, although this
is {\em not\/} a morphism! 
    
Recall that a morphism $a'\to a$ over $f:S'\to S$ is {\em cartesian}
if every morphism $b\to a$ over $f$ factors uniquely into
$b\xrightarrow{\tilde{g}} a'\to a$ with $p(\tilde{g})=\id_{S'}$. If
$a'\to a$ is a cartesian morphism over $f\colon S'\to S$, one sets,
slightly abusively, $a\times_S S':= a'$.
    
Following \cite{Rim}, a {\em fibration in categories} is a functor
$p\colon \bE\to \bB$ with the following properties:
\begin{itemize}
\item[]{\bf (FC1)}\quad For every morphism $f\colon S'\to S$ in $\bB$
and every object $a$ over $S$ there is a morphism $ a'\to a$ over $f$
that is cartesian.
    
\item[]{\bf (FC2)}\quad Compositions of cartesian morphisms are
cartesian.
\end{itemize}
    
The category $\bB$ is often called the {\em basis\/} of the fibration.
For $S\in \bB$ the {\em fibre} $\bE(S)$ is the subcategory of $\bE$
whose objects are those $a\in\bE$ with $p(a)=S$, and whose morphisms
$\varphi$ are the ones over $\id_S$, that is $p(\varphi)=\id_S$.
    
Recall that a groupoid is a category in which all morphisms are
isomorphisms. A fibration in categories $p\colon \bE\to \bB$ is called
a {\em fibration in groupoids\/} if each fiber $\bE(S)$ is a groupoid.
\end{sit}

These notions provide a natural framework for deformation theory.
As an example let us consider deformations of complex spaces.

\begin{exam}
\label{4.1}
Let $\bE$ be the category whose objects are the germs of flat
holomorphic maps $f\colon (X,X_0)\to(S,0)$, where $X_0=f^{-1}(0)$. The
morphisms into a second object given by $g\colon (Y,Y_0)\to(T,0)$ are
all cartesian squares
$$
\begin{diagram}[h=7mm]
(X,X_0)&\rTo^{f} &(S,0)\\
\dTo&&\dTo\\
(Y,Y_0)&\rTo^{g} &(T,0)\,.
\end{diagram}%
$$
The functor $p\colon \bE\to\Ans$ assigns to $f$ its basis $(S,0)$.

In a similar way one treats (flat) deformations of singularities,
coherent sheaves, or embedded deformations.
\end{exam}

\begin{sit}\label{4.2}
Let $\bB$ be one of the categories $\Ans$, $\Arts$ or $\Formals$.
Abusing notation again, we write $a\hto a'$ when the underlying
morphism in $\bB$ is an embedding. A fibration in groupoids $p\colon
\bE\to\bB$ is a {\em deformation theory\/} if the following
homogeneity property is satisfied.

\begin{itemize}
    \item[] {\bf (H)} \quad For every diagram in \bE,
\end{itemize}
\begin{diagram}[h=7mm]
a     &\rInto &a'   &      &&S&\rInto^{i}  &S' \\
\dTo  &       &     & \qquad\text{over}&&\dTo^\alpha \\
b     &       &     &      && T&&\\
\end{diagram}
with $i:S\hto S' $ an extension by a coherent $\calo_{S,0}$--module
$\calm$ and $\alpha:S\to T$ a finite map of germs, the fibred sum $b'
= a'\coprod_a b$ exists in $\bE$.

We remark that $b'$ lies necessarily over $S'\coprod_S T$, which in
turn exists as an analytic germ by \cite{Schu}.

The condition of homogeneity can be weakened to so called {\em
semihomogeneity\/}, see \cite{Rim}. We remark that the main
applications of this section remain true under this weaker condition
in view of the results of \cite{Fle2}. Note, however, that in all
reasonable geometric situations condition (H) above is satisfied.

If $a_0\in \bE(0)$ is a specific object over the reduced point, then a
{\em deformation\/} of $a_0$ over a germ $S=(S,0)$ is an object $a\in
\bE(S)$ together with a morphism $a_0\to a$ that lies necessarily over
$0\hookrightarrow (S,0)$.

Let $\bB$ again be one of the categories $\Arts$, $\Ans$ or $\Formals$
and let $p\colon \bE\to\bB$ be a fibred groupoid, $S=(S,0)$ a germ of
a complex space and $a\in \bE(S)$ an object over $S$. For a coherent
$\calo_S$--module $\calm$ an {\em extension of $a$ by $\calm$} is a
morphism $a\hto b$ such that the underlying morphism $S\hookrightarrow
T:= p(b)$ is an extension of $S$ by $\calm$. Two extensions
$a\hookrightarrow b$ and $a\hookrightarrow b'$ are said to be
isomorphic if there is a morphism $b\to b'$ that is compatible with
$a\hookrightarrow b$, $a\hookrightarrow b'$ and induces the identity
on $\calm$. We denote by $\Ex(a,\calm)$ the set of such isomorphism
classes.

In contrast, consider extensions $a\hookrightarrow b$ with
$p(b)=S[\calm]$, the trivial extension. Two such extensions
$\beta=(a\hookrightarrow b)$ and $\beta'= (a\hookrightarrow b')$ will
be called $S[\calm]$-isomorphic if there is a morphism of extensions
$b\to b'$ over $\id_{S[\calm]}$. The corresponding set of isomorphism
classes will be denoted $\Ex(a/S,\calm)$. This vector space is
commonly called the space of {\em infinitesimal deformations\/} (of 
first order).
\end{sit}

We recall the following facts, see e.g.\ \cite{Fle1,Fle2}.

\begin{prop}\label{4.3}
{\em 1.} The vector spaces $\Ex(a,\calm)$, $\Ex(a/S,\calm)$ define
(covariant) functors with respect to $\calm\in\Coh(S)$. They are
compatible with direct products: for coherent $\calo_S$--modules
$\calm,\caln$ and for $F$ one of these functors, one has naturally
$F(\calm\times \caln)\stackrel{\sim}{\to} F(\calm)\times F(\caln).$

{\em 2.} The functors in {\em (1)} are $\calo_S$-linear, i.e.\ they
carry natural $\calo_S$--module structures. Moreover they are
half-exact.

{\em 3.} There is a functorial exact sequence
$$
T^0_{S/\Sigma}(\calm)
\stackrel{\delta_{KS}}{\longrightarrow}
\Ex(a/S,\calm) \to \Ex(a,\calm)\to T^1_{S/\Sigma}(\calm)\,.
$$
\qed 
\end{prop}

Note that a functor $G\colon \Coh(S)\to\Sets$ compatible with finite
direct products and satisfying $G(0)\neq \emptyset$ always carries a
natural $\calo_S$--module structure, whence the first part of (2) is a
consequence of (1). Moreover, if $G'\colon \Coh(S)\to\Sets$ is a
second such functor, and if $G\to G'$ is a morphism of functors,
then the maps $G(\calm)\to G'(\calm)$ are necessarily $\calo_S 
$-linear.

We will refer to the sequence in (3) as the Kodaira-Spencer sequence.
Moreover, $\delta_{KS}$ is the so called {\em Kodaira-Spencer map.}

\begin{sit}
Let $p\colon \bE\to\Arts$ be a deformation theory and consider its
completion $\hat p\colon \hat\bE\to\Formals$, as described (dually) in
\cite{Rim}. Recall that a deformation $a\in\hat\bE(S)$ is called {\em
formally versal} if it satisfies the following lifting property: for
every morphism $b\hookrightarrow b'$ in $\hat\cate$ lying over a
closed embedding $T\hto T'$, and for every map $f\colon b\to a$, there
is a morphism $f'\colon b'\to a$ lifting $f$. Moreover, $a$ is said to
be {\em formally semiuniversal} if the induced map of tangent spaces
$T_{T',0}\to T_{S,0}$ is independent of the lifting.
    
By the theorem of Schlessinger (see \cite{Rim}), if $\Ex(a_0, \bbbc)$
is a vector space of finite dimension then a formally versal
deformation of $a_0$ exists. Moreover, we have the following criterion
for formal versality.
\end{sit}

\begin{prop}\label{4.4}{\em (\cite{Fle2})}
Let $a$ be a deformation of $a_0$ over the base $S\in\Formals$.
The following statements are equivalent.
\begin{enumerate}
\item The deformation $a$ is formally versal.
\item $\Ex(a, \calm)=0$ for every finite $\calo_S$--module $\calm$.
\item $\Ex(a, \calo_S/\fm_S)=0$.
\end{enumerate}
\end{prop}

In case $S/\Sigma$ is smooth, $T^1_{S/\Sigma}(\calm)$ vanishes for
every $\calm$ and so, in view of the Kodaira-Spencer sequence, $a$ is
formally versal if and only if the Kodaira-Spencer map
$\delta_{KS}\colon T^0_{S/\Sigma}(\calo_S/\fm_S)\to
\Ex(a/S,\calo_S/\fm_S)$ is surjective.

Next, we give a simple proof of a result by Z.~Ran \cite{Ran1}. We
state more generally a relative version over an arbitrary base
$\Sigma$. To formulate the last part of it, recall that an 
artinian germ $T\in\Ans$ is {\em curvilinear\/} if $\calot\cong
\bbbc[\![X]\!]/(X^n)$ as local $\bbbc$--algebras.

\begin{theorem}\label{4.5}
If $a_0 \in \cate(0)$ admits a formally semiuniversal
deformation $a\in \hat \cate(S)$ over some formal germ $S$,
then the following conditions are equivalent.

\begin{enumerate}
\item The germ $(S,0)$ is smooth over a closed subspace of the 
completion
$(\hat\Sigma, 0)$.
\item The functor $\calm \longmapsto \Ex (a/S, \calm)$
is right exact on $\Coh(S)$.
\item For every $b \in \cate(T)$ over an artinian germ $T\in\Ans$, the
map of infinitesimal deformations
$$
\Ex (b/T, \calot)\lto \Ex (b/T, \calot/\fm_T)
$$
is surjective.
\end{enumerate}
Moreover, if $\Sigma$ is a reduced point, then these are equivalent
to the following condition.
\begin{enumerate}
\item[4.] The map in {\em (3)} is surjective for every $b \in
\cate(T)$ over an artinian curvilinear germ $T\in\Ans$, .
\end{enumerate}
\end{theorem}

\begin{proof}
For (1)$\Rightarrow$(2) observe that the map $T^0_{S/\Sigma}(\calm)\to
\Ex(a/S,\calm)$ is surjective due to the versality of $a$. In turn, by
\ref{4.6} below, the functor $\calm\mapsto T^0_{S/\Sigma}(\calm)\cong
\Der_\Sigma(\calos, \calm)$ is right exact as $S$ is smooth over a
subspace of $\hat\Sigma$. Thus, $\Ex(a/S,\calm)$ is right exact in
$\calm$ as well. To show (2)$\Rightarrow$(3), let $b$, $T$ be as in
(3). By versality, there is a morphism $b\to a$ that lies over some
$\Sigma$--morphism $T\to S$. If now $b\hto b'$ is an extension of $b$
over $T\hto T[\calm]$, then the homogeneity condition yields an
extension $a\hto a':= a \coprod_b b'$ over $S\hto S[\calm]$ that
satisfies $b' = a'\times_{S}T$. Thus, there is a natural isomorphism
$$
\Ex (a/ S, \calm) \cong \Ex (b/T, \calm)
$$
for every artinian $\calo_{T}$--module $\calm$, whence (2) implies
(3). In order to show (3)$\Rightarrow$(1) consider the $n^{\rm th}$
infinitesimal neighbourhood $S_n$ of 0 in $S$ and set
$a_n:=a\times_SS_n$. Repeating the argument just given, we have
$\Ex(a/ S, \calm) \cong \Ex (a_n/S_n, \calm)$ for
any $\calm\in\Coh(S_n)$. Hence there is a commutative diagram
\begin{diagram}[h=7mm]
T^0_{S/\Sigma}(\calo_{S_n})&
\rTo^{(\delta_{KS})_n} & \Ex (a/S,\calo_{S_n})
&\cong &\Ex ( a_n/S_n, \calo_{S_n})\\
\dTo<{\alpha_n}  &&\dTo<{\beta_n}
&&\dTo<{\beta_n}\\
T^0_{S/\Sigma}(\calo_{S_0})&
\rTo^{(\delta_{KS})_0} & \Ex ( a/ S, \calo_{S_0})
&\cong& \Ex ( a_n/S_n, \calo_{S_0})\,.
\end{diagram}
Since $ a$ is formally semiuniversal, the map $(\delta_{KS})_0$ is
bijective and $(\delta_{KS})_n$ is surjective.  By
assumption $\beta_n$ is surjective and so $\alpha_n$ is surjective 
too.
Now the result follows from the Jacobian criterion \ref{4.6}
below.

It is obvious that (3) implies (4). Finally,
(4)$\Rightarrow$(1) follows with the same reasoning as above from the
smoothness criterion given in \ref{4.6} (5).
\end{proof}

In the proof above we have referred to the following smoothness
criteria that are essentially reformulations of the Jacobian 
criterion.

\begin{lem}\label{4.6}
For a morphism $A \to B$ of complete analytic
$\bbbc$-algebras the following conditions are equivalent.
\begin{enumerate}
\item There is an ideal $\fa \subseteq A$ and an isomorphism of
$A$--algebras
$$
B \cong (A/\fa)[\![T_1, \dots T_k]\!]\quad\text{for some } k
\ge 0\,;
$$
\item The $B$--module $\Omega^1_{B/A}$ is free; 

\item The functor $M\mapsto\Der_A(B, M)$ is right exact on finite
$B$--modules; 

\item With $B_n := B/\fm_B^{n+1}$, the natural map $\Der_A(B,B_n) \to
\Der_A(B, \bbbc)$ is surjective for all $n$.
\end{enumerate}
If $A=\bbbc$, then these conditions are also equivalent to:
\begin{enumerate}
\item[5.] The natural map $\Der_\bbbc(B,C)\to \Der_\bbbc(B, \bbbc)$ is
surjective for any artinian curvilinear $B$-algebra
$C\cong\bbbc[\![X]\!]/(X^n)$.
\end{enumerate}
\end{lem}

\begin{proof}
The implications (1)$\Rightarrow$(2)$ \Rightarrow$(3)$
\Rightarrow$(4)$ \Rightarrow$(5) are obvious. To show
(4)$\Rightarrow$(1), note first that
$$
\Der_A (B, B)\cong \lim_\leftarrow \Der_A(B,B_n)
$$
as $B$ is complete. Hence (4) implies that $\Der_A (B, B)\to
\Der_A(B,\bbbc)$ is surjective. This means that there are derivations
$\delta_1,\ldots, \delta_r:B\to B$ and elements
$x_1,\ldots,x_r\in\fm_B$ such that $\det(\delta_i(x_j))\not\in \fm_B$,
where $r:=\dim_\bbbc\Der_A(B,\bbbc)$. By the criterion of Lipman and
Zariski, see \cite[30.1]{Mat}, there is then an isomorphism $B\cong
C[\![X_1,\ldots,X_r]\!]$, where $C$ is an $A$--subalgebra of $B$. By
construction, $\Der_A(C,\bbbc)=0$, and so $C$ must be a quotient of 
$A$,
which gives (1).

Finally, assume that $A=\bbbc$ and that (5) is satisfied. Writing
$B=R/I$ with $R:=\bbbc[\![X_1,\ldots,X_r]\!]$ and $I\subseteq
\fm_R^2$, we need to show that $I=0$. If not, choose a power series
$f\ne 0$ of minimal order in $I$. After a linear change of coordinates
we may assume that $f=X_1^n+g$ with $g\in (X_2,\ldots,X_r)
\fm_R^{n-1}$. Now consider the curvilinear $B$-algebra
$C=B/(X_2,\ldots,X_r)\cong \bbbc[\![X_1]\!]/(X_1^n )$. The derivation
$\partial/\partial X_1$ on $R$ induces a derivation $B\to \bbbc$
that by assumption can be lifted to a derivation $B\to C$. Composing
with $R\to B$ yields a derivation, say, $\delta:R\to C$ with
$\delta(I)=0$. Using the product rule and the fact that
$\delta(X_i)\equiv 0 \bmod \fm_C$ for $i\ge 2$, it follows that 
$\delta
(g)=0$. On the other hand, $\delta(X_1)\equiv 1 \bmod \fm_C$, whence
$$
\delta(f) = nX_1^{n-1}\delta(X_1)=n X_1^{n-1} \cdot 1_C\ne 
0\quad\text{in $C$}\,,
$$
and this is a contradiction.
\end{proof}

\begin{rem}\label{4.7}
In case $A\ne \bbbc$, the condition (5) above is no longer equivalent
to the other ones as the following example shows. Consider a
$\bbbc$-algebra $A\cong R/I$, with $R\cong
\bbbc[\![X_1,\ldots,X_k]\!]$, and assume that there is an element
$a\in R$ that is integral over $I$ but not in $I$. If $\bar a$
denotes the residue class of $a$ in $A$, then the reader may verify
that the $A$-algebra $B:=A[\![T]\!]/(\bar aT)$ satisfies (5) although
it is not smooth over a quotient of $A$.
\end{rem}

\begin{defn}\label{4.8}
An {\em obstruction theory for} $a\in \hat\cate(S)$ over a formal germ
$S\in \Formals$ consists in a functor
$$
\Ob (a, -): \Coh_{art}(S)\lto \Coh_{art}(S),
$$
satisfying the following condition:
\begin{description}
\item[Ob1] There is a morphism of functors
$$
\ob : T^1_{S/\Sigma}(\calm) \lto \Ob(a,\calm)\quad,\quad \calm\in 
\Coh_{art}(S)\,,
$$
so that the sequence
$$
\Ex(a,\calm) \lto T^1_{S/\Sigma}(\calm)\xto{\ \ob\ }
\Ob (a,\calm)
$$
is exact for each $\calm$.
\end{description}
\end{defn}

In other words, if a $\Sigma$-extension $S \hto S'$ of $S$ by
$\calm\in \Coh_{art}(S)$ is given, then $\ob([S'])=0$ in
$\Ob(a,\calm)$ if and only if we can find an extension $a\hto a'$ of
$a$ by $\calm$ over $S\hto S'$. As an immediate consequence of the
Kodaira-Spencer sequence we obtain the following standard estimate.

\begin{cor}\label{4.9}
Let $a\in\hat\cate(S)$ be a formally semiuniversal deformation of
$a_0$ and assume that there is an obstruction theory for $a$. With
$$
k:=\dim_\bbbc\Ex(a_0,\bbbc)\quad\mbox{and}\quad
t:=\dim_\bbbc\Ob(a,\bbbc)\,,
$$
the basis $S$ of $a$ can be realized as a subspace of the formal germ
$(\Sigma\times\bbbc^k,0)\hat{\,}$ cut out by at most $t$ equations. In
particular, 
$$
\dim S\ge \dim\Sigma+k-t\,.
$$
\end{cor}

\begin{proof}
As was already observed in the proof of \ref{4.5}, we have
$\Ex(a_0,\bbbc)\cong\Ex(a/S,\bbbc)$, and, as $a$ is formally
semiuniversal, the map $T^0_{S/\Sigma}(\bbbc)\to \Ex(a/S,\bbbc)$ is
bijective. Hence we can write $\calo_{S,0}\cong
\Lambda[\![X_1,\ldots,X_k]\!]/I$ with $\Lambda:=\calo_{\hat\Sigma,0}$.
By definition, $T^1_{S/\Sigma}(\bbbc)\cong \Hom(I/\fm I,\bbbc)$, and
so the dimension of $T^1_{S/\Sigma}(\bbbc)$ is just the minimal number
of generators of $I$. In the extended Kodaira-Spencer sequence
$$
\ldots\lto \Ex(a,\bbbc)\lto T^1_{S/\Sigma}(\bbbc) \lto \Ob(a,\bbbc)
$$
the module $\Ex(a,\bbbc)$ vanishes by the versality of $a$. Hence
$T^1_{S/\Sigma}(\bbbc)$ injects into  $\Ob(a,\bbbc)$ and the
minimal number of generators of $I$ is bounded by $t$.
\end{proof}

Let $S$ be a formal germ over $(\Sigma, 0)$ and let
$\Lambda:=\calo_{\Sigma, 0}$, $A:=\calo_{S,0}$ denote the associated
local rings. It is well known that for a coherent $\calos$--module
$\calm$ the group $T^1_{S/\Sigma}(\calm)$ is canonically isomorphic to
$T^1_{A/\Lambda}(M)$, the group of $\Lambda$-algebra extensions of $A$
by $M:=\calm_{0}$. This group always contains
$\Ext^1_{A}(\Omega^1_{A/\Lambda},M)$ that can be identified as the
subgroup of those extensions whose associated Jacobi map is injective,
see \ref{1.1}. In case $M=\bbbc$, it contains in turn a further
distinguished subspace, $\Ex^c_{A/\Lambda}(\bbbc)$, the space of
curvilinear extensions, see \ref{9.1} and \ref{1.1}. With these
notations, the following improvement of \ref{4.5}, due to Kawamata,
can also easily be deduced.

\begin{cor}\label{4.10}
Let $a\in\hat\cate(S)$ be a formally semiuniversal deformation
of $a_0$ that admits an obstruction theory $\Ob(a,-)$.  If $V:=
\ob(\Ex^c_{A/\Lambda}(\bbbc)) \subseteq \Ob(a,\bbbc)$ is
the subspace of curvilinear obstructions, then
$$
\dim S\ge \dim\Sigma+\dim_\bbbc \Ex(a_0,\bbbc)-\dim_\bbbc V\,.
$$
\end{cor}

\begin{proof}
Note that $V\cong\Ex^c_{A/\Lambda}(\bbbc)$ as $T^1_{S/\Sigma}(\bbbc)$
injects into $\Ob(a,\bbbc)$. By \ref{1.1}, the vector space
$\Ex^c_{A/\Lambda}(\bbbc)$ is dual to $I/(J+\fm I)$, where $I$ is as
in the proof of \ref{4.9} and $J$ is the integral closure of
$\fm I$ in $I$. Now the result follows from \ref{1.2}.
\end{proof}

In the proofs of the preceding two corollaries we only used the module
$\Ob(a,\bbbc)$ so that it would have been sufficient to have just this
module of obstructions at our disposal. However, the next result
requires the full strength of the notion of obstruction theory.

\begin{prop} 
    \label{4.11}
Let $p: \cate \to\An_{\Sigma}$ be a deformation theory and $a
\in\hat\cate(S)$ an object over a germ $S \in \Formals$ that admits
an obstruction theory $\calm \mapsto \Ob(a,\calm)$. If there exists a
transformation $\psi : \Ob(a, \calm) \to G(\calm)$ into a left exact
functor $G$ on $\Coh_{art}(S)$ then the following hold.

{\em 1.} The restriction of $\psi\circ\ob$ to $\Ext_{S}^1(
\Omega^1_{S/\Sigma}, \calm) \hto T^1_{S/\Sigma}(\calm)$ is the zero
map. In other words, if $S\hookrightarrow S'$ is an extension of $S$
by $\calm\in \Coh_{art}(S)$ such that the Jacobi map $\calm\lto
\Omega^1_{S'/\Sigma} \otimes \calo_S$ is injective, then $\psi
\left(\ob[S']\right) = 0.$

{\em 2.}  If $a$ is formally versal then $\dim S\ge \dim_\bbbc\Ex
(a_0,\bbbc)-\dim_\bbbc K$ with $K:=\ker (\Ob(a,\bbbc)\to
G(\bbbc))$.
\end{prop}

\begin{proof}
The injective hull $\calj$ of $\calm$ can be obtained as a limit
$\lim\limits_{\to}\calj_\alpha$ with $\calm\subseteq
\calj_\alpha\subseteq\calj$ and each $\calj_\alpha$ finite artinian.
As
$$
\Ext_{S}^1(\Omega^1_{S/\Sigma},\calm)\to
\lim_{\to}\Ext_{S}^1(\Omega^1_{S/\Sigma},\calj_\alpha)\cong
\Ext_{S}^1(\Omega^1_{S/\Sigma},\calj)=0
$$
is the zero map, there is an index $\alpha$ so that
$\Ext_{S}^1(\Omega^1_{S/\Sigma},\calm)\to
\Ext_{S}^1(\Omega^1_{S/\Sigma}, \calj_\alpha)$ is already zero.
Restricting the map $\ob$ to $\Ext_{S}^1(\Omega^1_{S/\Sigma},-)$ 
yields a
commutative diagram
\begin{diagram}[h=7mm]
\Ext_{S}^1(\Omega^1_{S/\Sigma}, \calm)
&\rTo^\ob & \Ob(a,\calm)
&\rTo^\psi& G(\calm)\\
\dTo>0 &&\dTo && \dInto \\
\Ext_{S}^1(\Omega^1_{S/\Sigma}, \calj_\alpha)
&\rTo^\ob & \Ob(a,\calj_\alpha)
&\rTo^\psi & G(\calj_\alpha)\,,
\end{diagram}
and (1) follows.

In order to derive (2), note that the map $\ob$ embeds
$T^1_{S/\Sigma}(\bbbc)$ into $\Ob(a, \bbbc)$, as $a$ is formally
versal, and under this map $\Ext_{S}^1(\Omega^1_{S/\Sigma}, \bbbc)$
becomes thus isomorphic to a subspace of $K$ by (1). As
$\Hom_{S}(\Omega^1_{S/\Sigma}, \bbbc)$ is isomorphic to
$\Ex(a_{0},\bbbc)$, the space of infinitesimal deformations, the claim
follows from \ref{1.3}.
\end{proof}

In practise, it is cumbersome to construct obstruction theories for
objects over an arbitrary formal germ $S$. The following
considerations show that it is essentially sufficient to check the
existence of obstruction theories over an artinian base.

\begin{defn}
    \label{4.12}
Let $p:\cate\to \Ans$ be a deformation theory. An {\em obstruction
theory for $\cate$} consists in a collection of obstruction theories
$\Ob(a,-)$ for every $a\in \cate(S)$ over an artinian germ $S\in\Arts$
satisfying the following condition.
\begin{description}
\item[Ob2] For every inclusion $T\hto S$ in
$\Arts$ and every $a\in\cate(S)$ there are functorial
isomorphisms
$$
\Ob(a\times_ST,\calm)\xto{\ \cong\ } 
\Ob(a,\calm),\quad \calm\in\Coh(T)\,,
$$
such that the diagram
\begin{diagram}[h=7mm]
T^1_{T/\Sigma}(\calm) &\rTo^\ob & \Ob(a\times_ST,\calm)\\
\dTo && \dTo>\cong\\
T^1_{S/\Sigma}(\calm) &\rTo^\ob & \Ob(a,\calm)
\end{diagram}
commutes.
\end{description}
\end{defn}

If $\cate$ admits an obstruction theory and if $a\in \cate(S)$ over $S
\in \Formals$ is a formally versal deformation, then for every module
$\calm\in\Coh_{art}(S)$ the sequence of groups $\Ob(a_n,\calm)_{n\gg 
0}$,  where $a_n:=a|S_n$ is the
restriction to the $n^{th}$ infinitesimal neighbourhood $S_n$ of $0\in
S$, is essentially constant. Accordingly,
$$
\Ob(a,\calm):=\Ob(a_n,\calm),\quad n\gg 0,
$$
constitutes an obstruction theory for $a$ as the sequences
$$
\Ex(a_n,\calm) \lto T^1_{S_n/\Sigma}(\calm)\lto
\Ob (a_n,\calm)
$$
are exact for $n\gg 0$, and taking the direct limit results in the 
exact sequence
$$
\Ex(a,\calm)\cong\lim_\rightarrow \Ex(a_n,\calm)\lto
T^1_{S/\Sigma}(\calm)\cong \lim_\rightarrow
T^1_{S_n/\Sigma}(\calm)\lto
\Ob (a,\calm)\,.
$$

Assume now that $\cate$ admits an obstruction theory and that for
every $a\in\cate(S)$ over an artinian germ $S$ there is a
transformation $\Ob(a,-)\to G(a,-)$ into a left exact functor
$G(a,-)$ on  $\Coh_{art}(S)$. Furthermore, suppose that for every
inclusion $T\hto S$ in $\Arts$ there is an isomorphism
$G(a\times_ST,\calm)\xto{\ \cong\ }G(a,\calm)$,
$\calm\in\Coh T$ such that the diagram
\begin{diagram}[h=7mm]
\Ob(a\times_ST,\calm) &\rTo & G(a\times_ST,\calm)\\
\dTo_\cong && \dTo_\cong\\
\Ob(a,\calm) &\rTo & G(a,\calm)
\end{diagram}
commutes. Under these assumptions we get the following
corollary.

\begin{cor}\label{4.13}
If $a\in\hat\cate$ is a formally versal deformation of $a_0$, then
$$
\dim S\ge \dim_\bbbc\Ex (a_0,\bbbc)-\dim_\bbbc\ker (\Ob(a_0,\bbbc)\to
G(a_0,\bbbc))\,.
$$
\end{cor}

\begin{proof}
Indeed, the functor $G(a,\calm):=G(a_n,\calm), n\gg 0,$ is left exact
on $\Coh_{art}(S)$ and satisfies the assumptions of \ref{4.11}.
\end{proof}

\section{Applications}

\subsection*{Deformations of modules}
Let $f:X\to\Sigma$ be a flat holomorphic map and $0\in \Sigma$ a fixed
point. In the following we will consider deformations of coherent
modules on $X$. These deformations constitute a deformation theory
$p:\cate\to\Ans$, where the objects over a germ $S=(S,0)\in\Ans$ are
coherent $S$--flat modules $\calf$ on $X_S:=X\times_\Sigma S$. A
morphism into another module $\cale\in \cate(T)$ is given by a
morphism of $\Sigma$-germs $g:S\to T$ together with an isomorphism
$(\id_X\times g )^*(\calf) \stackrel{\sim}{\to} \cale$.

For a germ $S\in\Ans$ let $f_S:X_S\to S$ denote the projection. There
is a canonical map $T^1_{S/\Sigma}(\caln)\to T^1_{X_S}(f_S^*\caln)$
that assigns to an extension $S'$ of $S$ by $\caln\in\Coh S$ the
extension $X_{S'}$ of $X_S$ by $f_S^*\caln$. If $\calf\in\cate(S)$ is
a deformation over $S$ then by \ref{3.4} $\langle
[X_{S'}],\At(\calf)\rangle =0$ if and only if there is a module
$\calf'$ on $X_{S'}$ that forms an extension of $\calf$ by
$\calf\otimes_\calos\caln:=\calf\otimes_{f_S^{-1}\calos}
f_S^{-1}\caln$; observe that $\calf\dotimes_{\calo_{X_S}}f_{S}^*\caln
\cong \calf\otimes_{\calo_{X_S}}f_{S}^*\caln$ as $\calf$ is flat over
$S$. By a standard result in commutative algebra, see \cite[7.7]{Mat},
the $S'$--module $\calf'$ is then automatically flat. Hence we obtain
the following result.

\begin{lem}\label{5.1}
The composed maps
\begin{diagram}
\ob: T^1_{S/\Sigma}(\caln)
& \rTo^{can} &
T^1_{X_S}(f_S^*\caln) &
\rTo^{\langle *,-\At(\calf)\rangle}&
\Ext^2_{X_S}(\calf, \calf\otimes_\calos \caln)
\end{diagram}
for $\caln\in \Coh S$ define an obstruction theory for the deformation
theory of coherent modules on $X$ in the sense of {\em \ref{4.8},
\ref{4.12}\/}.
\end{lem}

It is well known that the space of infinitesimal deformations
of $\calf$ over $S[\caln]$ is just $\Ext^1_{X_S}(\calf,
\calf\otimes_\calos \caln)$.  Applying the results of the
previous sections we obtain the following theorem.

\begin{theorem}\label{5.2}
Assume that $X$ is proper and smooth over $\Sigma$ and that
$X_0:=f^{-1}(0)$ is bimeromorphically equivalent to a K\"ahler
manifold. Let $\calf_0$ be a coherent module on $X_0$ and let
$S=(S,0)\in\Ans$ be the basis of a semiuniversal deformation of
$\calf_0$. If
$$
\sigma :\Ext^2_{X_0} (\calf_0, \calf_0) 
\xto{\quad} 
\prod_{n \ge 0} H^{n+2}(X_0, \Omega^n_{X_0})
$$
denotes the semiregularity map as in {\em \ref{3.1}}, then the
following hold.

{\em 1.}  $\dim S \ge \dim_\bbbc \Ext^1_{X_0}(\calf_0, \calf_0) -
\dim_\bbbc \ker \sigma$.

{\em 2.}  If $\sigma$ is injective then $S$ is smooth at $0$ over a
closed subspace of $\hat\Sigma$.
\end{theorem}

\begin{proof}
Let $\calf$ be a deformation of $\calf_0$ over an artinian germ
$T\in\Arts$, whence $\calf$ is an $\calo_{X_T}$--module that is flat
over $T$ and restricts to $\calf_0$ on $X_0$. By \ref{3.1}, for every
coherent module $\caln$ on $T$ there is a semiregularity map
$$
\sigma_\caln:
\Ext^2_{X_T} (\calf, \calf \otimes_\calot \caln)
\lto \prod_{n \ge 0}
H^{n+2}(X_T,\caln\otimes_\calot\Omega^n_{X_T/T}),
$$
and by \ref{3.5}(1), this map is compatible with base change.
According to \cite{Del1}, the functor
$$
\caln\longmapsto H^p(X_T,
\caln\otimes_\calot\Omega^q_{X_T/T}),\quad  \caln\in\Coh T\,,
$$
is exact; note that using the results of \cite{Fuj}, Deligne's 
original
result extends to the case of compact manifolds that are
bimeromorphically equivalent to K\"ahler manifolds with the same 
proofs
as in (loc.cit.). Applying \ref{4.11}, claim (1) follows.

Finally assume that $\sigma=\sigma_\bbbc$ is injective.
Using induction on the length of $\caln$ it follows easily
that $\sigma_\caln$ is injective for all $\caln\in\Coh T$.
In particular, the functor $\caln\mapsto \Ext^2_{X_T} (\calf,
\calf \otimes_\calot \caln)$ is exact on the left and
therefore $\caln\mapsto \Ext^1_{X_T} (\calf, \calf
\otimes_\calot \caln)$ is exact on the right. Now
\ref{4.5} shows that $S$ is smooth over a closed subspace of 
$\hat\Sigma$,
as claimed.
\end{proof}

In case of deformations of modules on a fixed complex space
$X=X_0$, that is, when $\Sigma$ is a reduced point, the
result above holds without assuming that $X$ is smooth.

\begin{theorem}\label{5.3}
Let $\calf_0$ be a coherent module on $X$ with a finite
dimensional space of infinitesimal deformations
$\Ext^1_X(\calf_0, \calf_0)$, and let $S\in\Formals$ be the
basis of a formally semiuniversal deformation of $\calf_0$.
If $\calf_{0}$ has locally finite projective
dimension as an $\calox$--module then the semiregularity map
$$
\sigma :
\Ext^2_X (\calf_0, \calf_0) \lto
\prod_{n \ge 0} H^{n+2}(X, \Lambda^n\bbbl_X)
$$
is defined and
$$
\dim S \ge \dim_\bbbc \Ext^1_X(\calf_0, \calf_0) -
\dim_\bbbc \ker\sigma\,.
$$
In particular, if $\sigma$ is injective then $S$ is smooth.
\end{theorem}

\begin{proof}
Let $\calf$ be a deformation of $\calf_0$ over an artinian germ
$T\in\Ans$, whence $\calf$ is an $\calo_{X\times T}$--module that is
flat over $T$ and induces $\calf_0$ on $X$. The functor
$$
\caln\longmapsto H^p(X\times T,\caln \otimes_\calot
\Lambda^q \bbbl_{X\times T/T})\cong
\caln \otimes_\bbbc H^p(X, \Lambda^q\bbbl_X), \quad
\caln\in\Coh T,
$$
is exact.  As the semiregularity map is compatible with base
change $T\to S$, see \ref{3.5}, the result follows as before
from \ref{4.11}.
\end{proof}

As a special case this contains the result of Artamkin-Mukai,
\cite{Art,Muk}, that the injectivity of the trace map $\Ext^2_X
(\calf_0, \calf_0) \lto H^2(X, \calox)$ implies smoothness of the
basis of the semiuniversal deformation of $\calf_0$.

The proof of the following variant is similar and left to the
reader.

\begin{prop} \label{5.4}
Let $X$ be a complex space embedded into a
complex manifold $M$.  Let $\calf_0$ be a coherent
$\calox$--module with $\dim_\bbbc \Ext^1_X(\calf_0,
\calf_0)<\infty$. The dimension of the basis $(S,0)$ of a formally
semiuniversal deformation of $\calf_0$ satisfies then
$$
\dim S \ge \dim_\bbbc \Ext^1_X(\calf_0, \calf_0) - \ker 
\dim_\bbbc\sigma'\,,
$$
where
$$
\sigma':\Ext^2_X (\calf_0, \calf_0)
\xto{\ can\ }
\Ext^2_M(\calf_0, \calf_0)
\xto{\ \sigma\ }
\prod_{n\ge 0} H^{n+2}(M, \Omega^n_M).
$$
is the composition of the semiregularity map for $\calf_0$ as coherent
$\calo_{M}$--module with the canonical map between $\Ext$--functors.
If $\sigma'$ is injective then $S$ is smooth. \qed
\end{prop}

\begin{rems} \label{5.5}
1. If in \ref{5.3} the module $\calf$ is supported on a closed
subspace $Z$ of $X$ then the map $\sigma$ factors through a map
$$
\sigma_Z :
\Ext^2_X (\calf_0, \calf_0) \lto
\prod_{n \ge 0} H^{n+2}_Z(X, \Lambda^n\bbbl_X)\,,
$$
see \ref{3.5} (2). It is clear from the proof that the conclusion of
\ref{5.3} also holds with $\sigma_Z$ instead of $\sigma$.

2. Ideally, the map $\tau$ in \ref{5.4} should factor through a map
$$
\Ext^2_X(\calf_0, \calf_0) \lto\prod_{n\ge 0}  IH^{n+2,n},
$$
where $IH$ denotes intersection cohomology.
\end{rems}

\subsection*{The Hilbert scheme}
If $f:X\to \Sigma$ is a holomorphic map, let $H_{X/\Sigma}$ be the
relative {\em Douady space\/} of $X$ that represents the Hilbert
functor $\Hilb:\Ans\to\Sets$, where $\Hilb(S)$ is the set of all
closed subspaces of $X_S:=X\times_\Sigma S$ that are proper and flat
over $S$. In this section we will give several smoothness criteria for
the Douady space that generalize results of \cite{Blo,Kaw1,Ran2}.

For this, it is convenient to consider the deformation theory
associated to the Hilbert functor. More generally, we will study the
deformation theory $p:\cate\to \Ans$, where an object of $\cate$ over
a germ $S\in\Ans$ is a subspace $Z\subseteq X_S$ that is flat over
$S$. Note that $Z$ is not required to be proper over $S$. A morphism
of $Z$ into another object, say, $Z'\subseteq X_{S'}$ consists in a
morphism $g:S'\to S$ such that $(\id_X\times g)^{-1}(Z)=Z'$. It is
well known and easy to see that this constitutes a deformation theory
as explained in Section 6. Let $\Ex(Z/S,\caln)$, $\caln\in \Coh(S)$,
be the space of infinitesimal deformations of $Z\subseteq X_S$. The
following lemma is well known.

\begin{lem}\label{5.6}
If $Z\subseteq X_S$ is an $S$--flat subspace with ideal sheaf
$\calj\subseteq\calo_{X_S}$, then there is an isomorphism
$\Ex(Z/S,\caln)\cong
\Hom_{X_S}(\calj,\caloz\otimes_\calos
\caln)$, where $\caloz\otimes_\calos\caln:=
\caloz\otimes_{f_S^{-1}\calos}f_S^{-1}\caln$.\qed
\end{lem}

Note that $\Hom_{X_S}(\calj,\caloz\otimes_\calos
\caln)$ is just $T^1_{X/Z}(\caloz\otimes_\calos
\caln)$. In case that $Z$ is locally a complete
intersection in $X$ this is just the space of section of the
normal bundle $\caln_{Z/X}$ of $Z$ in $X$.

In order to describe an obstruction theory for $Z$ we will assume for
simplicity that $f$ is flat; see remark \ref{5.14} (3) for the general
case. In the flat case, there is for any coherent $\calos$--module
$\caln$ a canonical map $T^1_{S/\Sigma} (\caln)
\stackrel{\alpha}{\lrarrow} T^1_{X_S/\Sigma}
(\calo_{X_S}\otimes_\calos\caln)$ that assigns to an extension $[S']$
of $S$ by $\caln$ the extension $[X_{S'}]$ of $X$ by $\caln_{X_{S}}:=
\calo_{X_S} \otimes_\calos \caln$. Composing this with the map
$\gamma$ considered in \ref{5.7.a}, we get a map
$$
\ob:T^1_{S/\Sigma} (\caln)\lrarrow
T^2_{Z/X_S}(\caloz\otimes_\calos\caln)\,.
$$
Applying \ref{5.7.a}, gives the following lemma.

\begin{lem}\label{5.7}
If $f:X\to \Sigma$ is flat, then the  map
$\ob$ constitutes an obstruction theory for $Z$.\qed
\end{lem}

Let $f:X\to \Sigma$ be a flat map, $(S,0)$ a germ over
$\Sigma=(\Sigma,0)$, and let $Z\subseteq X_{S}$ be an $S$--flat
subspace with special fibre $Z_{0}$ over $0\in S$. If $\calo_{Z_{0}}$
has locally finite projective dimension as a module on
$X_{0}:=f^{-1}(0)$ then by \ref{3.4.a} the sheaf $\caloz$ has finite
projective dimension over $X$. Hence applying \ref{5.9} yields a
semiregularity map
$$
\tau_{\caln} : T^2_{Z/X_S} (\caln_{Z}) \lto \prod_{n \ge
0} H^{n+2} (X_S, \caln\dotimes_{\calos}
\Lambda^n\bbbl_{X_{S}/S}), \quad\caln\in \Coh S\,.
$$
In case of smooth maps $f$ we get the following result.

\begin{theorem}\label{5.11}
Assume that $X$ is proper and smooth over $\Sigma$ and that
$X_0:=f^{-1}(0)$ is bimeromorphically equivalent to a K\"ahler
manifold. Let $Z_0\subseteq X_0$ be a closed subspace. If
$$
\tau :
T^2_{Z_0/X_0} (\calo_{Z_0}) \lto \prod_{n \ge 0} H^{n+2}
(X_0, \Omega^n_{X_0})\ .
$$
denotes the semiregularity map for $Z_0\subseteq X_0$, then
the following hold.

{\em 1.} If $\tau$ is injective then the Douady space $H_{X/\Sigma}$ 
is
smooth over a closed subspace of $\Sigma$ in a neighbourhood of its
point $[Z_0]$.

{\em 2.} The dimension of $H_{X/\Sigma}$ at $[Z_{0}]$ satisfies
$$
\dim_{[Z_0]} H_{X/\Sigma} \ge 
\dim_\bbbc T^1_{Z_0/X_0}(\calo_{Z_0}) -
\dim_\bbbc \ker \tau\,.
$$
\end{theorem}

The {\em proof\/} follows again from \ref{4.11} along the same line
of arguments as in \ref{5.2}.\qed
\medskip

We note that in general the injectivity of $\tau$ does not imply that
$H_{X/\Sigma}$ is smooth over $\Sigma$. For instance, let $X\to
\Sigma$ be a smooth family of surfaces over a germ $(\Sigma,0)$ and
consider a curve $C\subseteq X_0$. The semiregularity map
$H^1(C,\caln_{C/X_0})\to H^2(X_0,\calo_{X_0})$ is dual to the
restriction map $H^0(X_0,\omega_{X_0})\to H^0(C,
\omega_{X_0}\otimes_{\calo_{X_0}}\calo_C)$ as we will show in Section
8. If $X\to \Sigma$ is a versal family of K3-surfaces and $C$ is a
connected reduced curve on $X_0$, then this restriction map is
bijective. Hence $H_{X/\Sigma}$ is smooth over a closed subspace of
$\Sigma$ at $[C]$. However, this subspace cannot be all of $\Sigma$
since there are no curves on the general K3-surface.

In the absolute case, when $\Sigma$ is a reduced point, the
following stronger result holds.

\begin{theorem}\label{5.12}
Let $Z$ be a closed subspace of a complex space
$X$ with a finite dimensional space of infinitesimal embedded
deformations $T^1_{Z/X}(\caloz)$ and let
$S\in\Formals$ be the basis of a formally
semiuniversal deformation of $Z$. Assume that $\caloz$ has
locally finite projective dimension as $\calox$--module.
With
$$
\tau :
T^2_{Z/X} (\caloz) \lto
\prod_{n \ge 0} H^{n+2}(X, \Lambda^n\bbbl_X)
$$
the semiregularity map, the following hold.

{\em 1.} The dimension of $S$ satisfies
$$\dim S \ge \dim_\bbbc T^1_{X/Z}(\caloz) -
\dim_\bbbc \ker\tau\,.
$$
In particular, if $\tau$ is injective then $S$ is smooth.

{\em 2.} If $Z$ is compact then $\dim_{[Z]} H_X \ge \dim_\bbbc
T^1_{X/Z}(\caloz) - \dim_\bbbc \ker\tau\ .$ In particular, if
$\tau$ is injective then $H_X$ is smooth at $[Z]$.\qed
\end{theorem}

In the smooth case this specializes further to the following
corollary.

\begin{cor}\label{5.13}
Let $Z\subseteq X$ be a compact subspace of a complex manifold $X$.
With
$$
\tau :
T^2_{Z/X} (\caloz) \lto \prod_{n \ge 0} H^{n+2}
(X, \Omega^n_{X})
$$
the semiregularity map, we have $\dim_{[Z]} H_X \ge \dim_\bbbc
T^1_{Z/X}(\caloz) - \dim_\bbbc \ker \tau$. In particular, if $\tau$ is
injective then $H_X$ is smooth at $[Z]$.\qed
\end{cor}

As well, the result \ref{5.4} can be formulated in the case of
deformations of subspaces. We leave the straightforward
formulation and its proof to the reader.

\begin{rems}\label{5.14}
1. Note that for a locally complete intersection $Z\subseteq X$ with
normal bundle $\caln_{Z/X}$ there is a canonical isomorphism
$T^k_{Z/X}(\calm)\cong H^{k-1}(Z,\caln_{Z/X}\otimes \calm)$ for every
$\caloz$--module $\calm$ and $k\ge 0$. Hence in this case the
statements \ref{5.11}--\ref{5.13} above hold with $T^2_{Z/X}(\caloz)$
replaced by $H^1(Z,\caln_{Z/X})$.

2. In analogy with \ref{5.5}(1), the map $\tau$ in \ref{5.13} factors
through a map $\tau_{Z}:T^2_{Z/X} (\caloz) \to \prod_{n \ge 0}
H^{n+2}_{Z} (X, \Omega^n_{X})$, see \ref{3.5a} (3). It is clear from
the proof that the conclusion of \ref{5.13} remains true for
$\tau_{Z}$ instead of $\tau$. A similar remark applies to \ref{5.12}.

3. We note that there is also an obstruction theory for embedded
deformations if $f:X\to \Sigma$ is not flat. To show this, let $S$ be
a space over $\Sigma$ and $Z\subseteq X_S$ an $S$--flat subspace.
Consider $X_S$ and $Z$ as subspaces of $X\times S$ via the diagonal
embedding. For a coherent $\calos$--module $\caln$ there are natural
maps
$$
T^1_{S/\Sigma}(\caln)
\stackrel{\alpha}{\lrarrow}
T^1_{X\times S/\Sigma\times \Sigma}(\calo_{X\times
S}\otimes_\calos \caln)
\stackrel{\beta}{\lrarrow}
T^1_{X\times S/\Sigma\times
\Sigma}(\caloz\otimes_\calos\caln)\,,
$$
In terms of extensions, if $S'$ is an extension of $S$ by $\caln$ then
$\alpha([S'])=[X\times S']$ and $\beta([X\times S'])$ is the induced
extension of $X\times S$ by $\caloz\otimes_\calos\caln$. The embedded
deformation $Z\subseteq X_S$ can be extended to an embedded
deformation $Z\subseteq X_{S'}$ if and only if $\alpha([X\times S'])$
is in the image of the natural map
$$
\gamma: T^1_{Z/\Sigma}(\caloz\otimes_\calos\caln) \to
T^1_{X\times S/\Sigma\times
\Sigma}(\caloz\otimes_\calos\caln)\,.
$$
This map embeds into a long exact cohomology sequence: namely, if
$\calc^\sbullet$ denotes the mapping cone of $\bbbl_{X\times
S/\Sigma\times \Sigma}\dotimes_{\calo_{X\times S}} \caloz \to
\bbbl_{Z/\Sigma}$, then the cokernel of $\gamma$ embeds into
$\Ext^2_Z(\calc^\sbullet, \caln\otimes_\calos\caloz)$. Thus the
composition
$$
T^1_{S/\Sigma}(\caln)\stackrel{\beta\circ\alpha}{\llrarrow}
 T^1_{X\times S/\Sigma\times
\Sigma}(\caloz\otimes_\calos\caln)
\to
\Ext^2_Z(\calc^\sbullet , \caloz\otimes_\calos\caln)
$$
gives an obstruction theory.

Note that in case of a flat map $f:X\to \Sigma$, the complexes
$$
\bbbl_{X_S/S}\dotimes_{\calo_{X_S}}\caloz
\quad\text{and}\quad
\bbbl_{X\times S/\Sigma\times
\Sigma}\dotimes_{\calo_{X\times S}} \caloz
$$
are quasiisomorphic. Hence $\calc^\sbullet$ becomes quasiisomorphic to
the mapping cone of $\bbbl_{X_S/S}\dotimes_{\calo_{X_S}}\caloz \to
\bbbl_{Z/S}$ and so is quasiisomorphic to $\bbbl_{Z/X_S}$. This is
the same obstruction theory we described before.

In the general case, the reader may easily verify that there is
a natural map $\bbbl_{Z/X_S}\to\calc^\sbullet$ that induces
a map
$$
\Ext^2_Z(\calc^\sbullet, \caloz\otimes_\calos \caln)
\lto
T^2_{Z/X_S}(\caloz\otimes_\calos \caln)\,.
$$
Taking the composition with the map $\tau$ as defined in
\ref{5.9}, we arrive at a semiregularity map also in this case.
\end{rems}

\subsection*{The $\Quot$--functor}

Let $f:X\to \Sigma$ be a holomorphic map and $\cale$ a coherent sheaf
on $X$. Generalizing partly the results of the previous part we will
study the $\Quot$-space $Q_{\cale/\Sigma}$ that was constructed as a
complex space by Douady in the absolute case and by Pourcin in the
relative case, see \cite{Dou,Pou}. We remind the reader that it
represents the functor $\Quot_{\cale/\Sigma}:\Ans\to\Sets$, where
$\Quot_{\cale/\Sigma}(S)$ is the set of all quotients of
$\cale_S:=p_S^*(\cale)$ that are proper and flat over $S$; here
$p_S:X_S:=X\times_{\Sigma}S \to X$ denotes the projection.

Again it is convenient to consider the associated
deformation theory, say, $p:\cate\to \Ans$.  An object of
$\cate$ over a germ $S\in\Ans$ is a quotient $\calq$ of
$\cale_S$ that is flat over $S$.  A morphism of $\calq$ into
another object, say, $\calq'$ defined over the germ $S'$ is
given by a morphism $g:S'\to S$ such that $(\id_X\times
g)^*(\calq)=\calq'$ as quotients of $\cale_{S'}$.  It is
well known and easy to see that this constitutes a
deformation theory as in Section 3.  The space of
infinitesimal deformations $\Ex(\calq/S,\caln)$, $\caln\in
\Coh(S)$, is described in the following well known lemma.

\begin{lem}\label{5.15}
If $\cale_S\to\calq$ is an $S$--flat quotient with kernel
$\calf:=\ker(\cale_S\to\calq)$ then there is a natural
isomorphism  $\Ex(\calq/S,\caln)\cong
\Hom_{X_S}(\calf,\calq\otimes_\calos
\caln)$.\qed
\end{lem}

In case $\cale$ is flat over $\Sigma$, there is furthermore a well
known obstruction theory for $\calq$. Since there seems to be no
explicit reference for this in the relative case, we describe in
brief the construction. First note the following simple lemma whose
proof is left to the reader.

\begin{lem}\label{5.15.a}
Let $X$ be a complex space and $X\subseteq X'$ an extension of $X$
by a coherent $\calox$--module $\cali$. For coherent $\calox$--modules
$\calg$, $\calh$ consider the map
$$
\mu:\Ext^1_{X'}(\calg, \calh)\to
\Hom_{\calox}(\calg\otimes\cali,\calh)
$$
that assigns to an $X'$-extension $[\calg']\in \Ext^1_{X'}(\calg,
\calh)$ the homomorphism $\mu(\calg')$ obtained from the
multiplication map $\calg'\otimes \cali\to \calg'$ through the natural
factorization
\begin{diagram}[s=5mm]
\calg'\otimes \cali &\rTo^{proj} & \calg\otimes\cali
&\rTo^{\mu(\calg')} & \calh &\rInto &\calg'\,.
\end{diagram}
The map $\mu$ is functorial in $\calg$ and $\calh$ and thus is
$\Gamma(X',\calo_{X'})$-linear. Moreover, $\ker \mu$ is
canonically isomorphic to $\Ext^1_{X}(\calg, \calh)$.
\qed
\end{lem}

We note that $\mu$ can also be described as the
boundary map in the spectral sequence
$E^{pq}_{2}=\Ext^p_{X} (\cTor_{q}^{\calo_{X'}}
(\calg,\calox),\calh) \Rightarrow \Ext^{p+q}_{X'}
(\calg,\calh)$.  However, the more explicit description
given above is better suited for our needs.

Let us return to the situation as described before
\ref{5.15.a}, and consider the composition of the canonical
maps
$$
T^1_{S/\Sigma}(\caln)
\stackrel{\alpha}{\lrarrow}
\Ext^1_{X_{S'}}(\cale_{S},\cale_{S}\otimes_\calos \caln)
\stackrel{\beta}{\lrarrow}
\Ext^1_{X_{S'}}(\cale_{S}, \calq\otimes_\calos\caln),
$$
where $\alpha$ maps an extension $[S']$ of $S$ by $\caln$ onto the
class of $\cale_{S'}$. There is furthermore an exact $Ext$--sequence
$$
\Ext^1_{X_{S'}}(\calq,\calq\otimes_\calos \caln)
\stackrel{\gamma}{\to}
\Ext^1_{X_{S'}}(\cale_{S}, \calq\otimes_\calos \caln)
\stackrel{\delta}{\to}
\Ext^1_{X_{S'}}(\calf, \calq\otimes_\calos\caln)\,.
$$
For $\calg=\cale_{S}, \calh= \cale_{S}\otimes_\calos\caln$, the map
$\mu$ described in \ref{5.15.a} associates to $\cale_{S'}$ just the 
identity on
$\cale_{S}\otimes_\calos\caln$. Using again \ref{5.15.a}, the diagram
\begin{diagram}[s=7mm]
\Ext^1_{X_{S'}}(\cale_{S},\cale_{S}\otimes_\calos\caln) &
\rTo^{\delta\circ\beta} &
\Ext^1_{X_{S'}}(\calf,\calq\otimes_\calos \caln)\\
\dTo^{\mu} && \dTo^\mu\\
\Hom_{X_{S}}(\cale_{S}\otimes_\calos\caln 
,\cale_{S}\otimes_\calos\caln)
& \rTo^{\ can\ \ } &
\Hom_{X_{S}}(\calf\otimes_\calos\caln ,\calq\otimes_\calos \caln)
\end{diagram}
commutes, whence the extension $\delta\circ\beta([\cale_{S'}])$ maps
to $0$ under $\mu$ and so can be identified with an element of
$\Ext^1_{X_{S}}(\calf ,\calq\otimes_\calos \caln)$. In other words,
$\delta\circ\beta\circ\alpha$ factors through a map
$$
\ob: T^1_{S/\Sigma}(\caln)\lto
\Ext^1_{X_S}(\calf,\calq\otimes_\calos\caln)\,.
$$
If $\ob([S'])$ vanishes then $\beta([\cale_{S'}])$ is in the
image of $\gamma$, and the corresponding extension $\calq'$ gives
a lifting of $\calq$ to $S'$. This establishes the
following result.

\begin{prop}
    \label{5.16}
If $\cale$ is flat over $\Sigma$ then the map $\ob$ just defined 
provides an obstruction theory for $\cale$.\qed
\end{prop}

It is now immediate how to define a semiregularity map on
$\Ext^1_{X_S}(\calf,\calq\otimes_\calos \caln)$.

\begin{defn}\label{5.17}
Let $X\to\Sigma$ and $\cale$ be as above and let $\calq$ be an
$S$--flat quotient of $\cale_S$ over some germ $S=(S,0)$ over
$\Sigma$. Assume that $\calq$ has locally finite projective dimension
on $X$. The composition of the boundary map
$\Ext^1_{X_S}(\calf,\calq\otimes_\calos \caln) \to
\Ext^2_{X_S}(\calq,\calq\otimes_\calos\caln)$ with the semiregularity
map $\sigma$ defined in \ref{3.1} yields a map
$$
\tau:\Ext^1_{X_S}(\calf,\calq\otimes_\calos\caln)
\lto \prod_{n\ge 0}
H^{n+2}(X_S,\caln\otimes_\calos\Lambda^n\bbbl_{X_S/S})
$$
that we call the {\em semiregularity map for $\calq$}.
\end{defn}

Now it is possible to deduce results analogous to
\ref{5.11}--\ref{5.13} and \ref{5.4}. As a sample we restrict to the
following application; the proof is similar to that of \ref{5.11} and
left to the reader.

\begin{theorem}\label{5.18}
Assume that $X$ is proper and smooth over $\Sigma$ and that
$X_0:=f^{-1}(0)$ is bimeromorphically equivalent to a K\"ahler
manifold. Let $\cale$ be a $\Sigma$--flat coherent sheaf on $X$ and
$\cale_0\to\calq_0$ be a quotient of $\cale_0:=\cale|X_0$. With
$\calf_0$ the kernel of $\cale_0\to\calq_0$ and
$$
\tau :
\Ext^1_{X_0}(\calf_0,\calq_0) \lto
\prod_{n \ge 0} H^{n+2} (X_0, \Omega^n_{X_0})
$$
the semiregularity map for $\calq_0$, the following hold.

{\em 1.} If $\tau$ is injective then the Quot-space $Q_{\cale/\Sigma}$
is smooth over a closed subspace of $\Sigma$ in a neighbourhood of
$[\calq_0]$.

{\em 2.} The dimension of the Quot-space $Q_{\cale/\Sigma}$ at the
point $[\calq_0]$ satisfies
$$
\dim_{[\calq_0]} Q_{\cale/\Sigma} \ge \dim_\bbbc
\Hom_{X_0}(\calf_0,\calq_0) - \dim_\bbbc \ker \tau\,.
$$
\qed
\end{theorem}

Considering a subspace of $X$ as a quotient of $\calox$, the Douady
space becomes a special case of the Quot-space. However, note that
\ref{5.18} does not imply the corresponding result for the Douady
space. The reason is that the obstruction theory for the Quot-functor
does not specialize to the obstruction theory for the Hilb-functor.
For instance, if $Z\subseteq X$ is a complete intersection of
dimension 0 then $T^2_{Z/X}(\caloz\otimes_\calos \caln)$ vanishes
while $\Ext^1_X(\calj, \caloz)$ is isomorphic to the space of sections
of the second exterior power of the normal bundle of $Z$ in $X$. In
general, $T^2$ provides a much smaller obstruction theory than the
$\Ext^1$--functor. Hence it is worthwhile to treat these cases
separately.

\begin{rem}\label{5.19}
In analogy with \ref{5.14}(3), we note that $\cale$ need not be flat
over $\Sigma$ in order to obtain an obstruction theory for the
$\Quot$--functor. Indeed, with $S$, $\caln$, $\calq$ as in \ref{5.15}
and $S'$ an extension of $S$ by $\caln$, consider the composition
$$
T^1_{S/\Sigma}(\caln) \stackrel{\alpha}{\lrarrow}
\Ext^1_{X\times S'}(\pi_{S}^*(\cale),
\pi_{S}^*(\cale)\otimes_\calos \caln)
\stackrel{\beta}{\lrarrow} \Ext^1_{X\times
S'}(\pi_{S}^*(\cale), \calq\otimes_\calos\caln)\,,
$$
where $\pi_S:X\times S\to S$ denotes the projection.
Here $\alpha$ maps an extension $[S']$ of $S$ by $\caln$ onto
the class of $\pi_{S'}^*(\cale)$. Now
$$
\Ext^1_{X\times S'}(\pi_S^*(\cale),\calq\otimes_\calos\caln)
\cong
\Ext^1_{X_{S'}}(\cale\dotimes_{\calo_\Sigma}\calos
,\calq\otimes_\calos\caln)\,,
$$
with $\cale\dotimes_{\calo_\Sigma}\calos:=
\pi_S^*(\cale)\dotimes_{\calo_{X\times S'}}\calo_{X_{S'}}$. For
$\calq$ to be extendable to $S'$, the element $\beta\circ\alpha([S'])$
must be in the image of the map $\gamma$ in the following exact
sequence:
$$
\Ext^1_{X_{S'}}(\calq, \calq\otimes_\calos\caln)
\xto{\gamma}
\Ext^1_{X_{S'}}(\cale\dotimes_{\calo_\Sigma}\calos,
\calq\otimes_\calos\caln)
\xto{\delta}
\Ext^1_{X_{S'}}(\calf^\sbullet,\calq\otimes_\calos\caln)\,,
$$
where $\calf^\sbullet$ is the mapping cone of
$\cale\dotimes_{\calo_\Sigma}\calos[-1]\to\calq[-1]$. As before, one
can verify that $\delta\circ\beta\circ\alpha$ factors through a map
$$
\ob:T^1_{S/\Sigma}(\caln)
\lto
\Ext^1_{X_S}(\calf^\sbullet, \calq\otimes_\calos\caln)\,,
$$
and this map constitutes an obstruction theory. The reader may
easily check that this gives rise to a semiregularity map on
$\Ext^1_{X_S}(\calf^\sbullet ,\calq\otimes_\calos\caln)$.
\end{rem}

\subsection*{Deformations of mappings} In this part we
will generalize the semiregularity map for embedded deformations, see
\ref{5.9}, to deformations of holomorphic maps $X_{0}\to Y_{0}$. For
the special case that $X_{0}$ is a stable curve over an algebraic
manifold $Y_{0}$, such a semiregularity map was independently
constructed by K.\ Behrend and B.\ Fantechi in order to define refined
Gromov-Witten invariants in certain situations. We consider the
following setup.

\begin{sit}\label{7.16}
Let $\pi:Y\to \Sigma=(\Sigma,0)$ be a fixed germ of a flat holomorphic
mapping, set $Y_0:=\pi^{-1}(0)$ and let $f_0:X_0\to Y_0$ be a morphism
of complex spaces with $X_0$ compact. By a deformation of
$X_0/Y_{0}$ over a germ $(S,0)\in\Ans$ we mean a commutative diagram
$$
\begin{diagram}[s=7mm,midshaft]
X&& \rTo^f&& Y_S:= Y\times_\Sigma S\\
& \rdTo(2.3,2)_q && \ldTo(2.3,2)_{\pi_S}\\
&& S
\end{diagram}
\leqno(1)
$$
such that $q$ is flat and proper and $f$ induces $f_0$ on the special
fibre $X_0=q^{-1}(0)$. This deformation is, abusively, denoted by
$X/Y_{S}$. Such deformations form in a natural way a deformation
theory over $\Ans$.
\end{sit}

It is well known, and follows easily from the existence of the
(relative) Douady space, that there are always convergent versal
deformations of $X_0/Y_{0}$, see \cite{Fle2,BKo}.

The infinitesimal deformations of $X_0/Y_{0}$ can be
described as follows.  Let $X/Y_{S}$ be a deformation of
$X_0/Y_{0}$ over $S$ as in (1) and let $\caln$ be a coherent
$\calos$--module.  An $S':=S[\caln]$-extension of $f:X\to Y_S$ by
$\caln$ consists in a deformation $X'/Y_{S'}$ of $X_0/Y_{0}$ over
$S'$ as in the diagram
$$
\begin{diagram}[s=7mm]
X' && \rTo^{f'} && Y_{S'}\\
& \rdTo && \ldTo_{\pi_{S'}}\\
&& S'=S[\caln]
\end{diagram}
\leqno(2)
$$
that induces on $S$ the given deformation $X/Y_{S}$ as in (1). Let
$\Ex(X/Y_{S},\caln)$ denote the group of these extensions. An
$S[\caln]$-extension $X'/Y_{S[\caln]}$ of $X/Y_{S}$ corresponds to an
extension of $X$ by $\caln_X:=\caln\otimes_{\calos}\calox$ that is a
space over $Y_S$. Hence we obtain the following lemma.

\begin{lem}\label{7.17}
There is a canonical isomorphism $
\Ex(X/Y_{S},\caln)\cong T^1_{X/Y_S}(\caln_X)$.\qed
\end{lem}

Next we will describe obstructions for extending deformations. Let
$X/Y_{S}$ be as in (1) and let $S\hto S'$ be an extension of $S$ by
$\caln$. The extension $Y_{S'}$ of $Y_S$ by
$\caln_{Y_{S}}:=\caln\otimes_\calos \calo_{Y_S}$ gives an element
$[Y_{S'}]$ of $T^1_{Y_{S}}(\caln_{Y_{S}})$. The composition of the two
canonical maps
$$
\gamma:T^1_{Y_S}(\caln_{Y_{S}})
{\lto}
\Ext^1_X(Lf^*(\bbbl_{Y_S/S}), \caln_X)
{\lto}
T^2_{X/Y_S}(\caln_X)
$$
admits the following interpretation in terms of extensions.

\begin{lem}\label{7.18}
The class of the extension $Y_{S'}$ is mapped to $0$ under $\gamma$ if
and only if there is an extension $X'$ of $X$ by $\caln_{X}$ together
with a map $f':X'\to Y_{S'}$ that induces $f$ on $Y$ and the identity
on $\caln_{X}$. In particular, the composite map
$$
\ob: T^1_{S/\Sigma}(\caln)\stackrel{can}{\lto}
T^1_{Y_S}(\caln_{Y_{S}})
\stackrel{\gamma}{\lto} T^2_{X/Y_S}(\caln_X)
$$
provides an obstruction theory in the sense of {\em \ref{4.8}\/}.
\end{lem}

\begin{proof}
To prove this statement, we use the existence of tangent functors
$T^i_{f}(-)$ for holomorphic mappings as constructed in
\cite{Fle1}, see also \cite{Ill}. These functors fit into an exact 
sequence
$$
\cdots\lto T^1_{f}(\caln_{Y_{S}})\stackrel{\beta}{\lto}
T^1_{Y_{S}}(\caln_{Y_{S}})\stackrel{\gamma}{\lto}
T^2_{X/Y_{S}}(\caln_{X})\lto \cdots ,\leqno (*)
$$
see \cite[3.4]{Fle1}. The group $T^1_{f}(\caln_{Y_{S}})$ is
canonically isomorphic to the set of all isomorphism classes of
extensions of $f$ by $\caln_{Y_{S}}$. Such extension is a holomorphic
map $f':X'\to Y'$, where $X'$, $Y'$ are extensions of $X$, $Y_{S}$ by
$\caln_{X}$, $\caln_{Y_{S}}$, respectively, with $f'$ inducing the
map $f$ on $X$ and the identity on $\caln_{X}$; see \cite[3.16]{Fle1}
for details. Moreover, the map $\beta$ in $(*)$ assigns to $[f']$ the
extension $[Y']$. In view of the exactness of the sequence $(*)$ this
proves the lemma.
\end{proof}

In a next step, we generalize \ref{5.8}(1) to arbitrary mappings.

\begin{prop}\label{7.19}
For every morphism of complex spaces $f:X\to Y$ and every complex of
$\calox$--module $\calm$ bounded below there are canonical maps
$$
T^k_{X/Y}(\calm)\lto
\Ext^k_{Y}(Rf_*(\calox),Rf_*(\calm)), \quad k\in \bbbz\,.
$$
In case $\calm=\calox$, this specializes to a map of graded Lie
algebras $T^\sbullet_{X/Y}(\calox)\to
\Ext^\sbullet_{Y}(Rf_*(\calox),Rf_*(\calox))$.
\end{prop}

\begin{proof}
Let $(X_{*},W_{*}, \calr_{*})$ be a resolvent for $X$ over $Y$ and let
$\calm_{*}\to\tilde\calm_{*}$ be a quasiisomorphism into a
$W_{*}$-acyclic complex $\tilde\calm_{*}$ of $\calo_{X_{*}}$--modules
as in \ref{2.6}. We need to construct a map
$$
T^k_{X/Y}(\calm)\cong
H^k(\Hom_{\calr_{*}}(\Omega^1_{\calr_{*}/Y},\tilde
\calm_{*}))\lto
\Ext^k_{Y}(Rf_*(\calox),Rf_*(\calm))\,.
$$
Let $f^{-1}(\caloy)$ be the topological preimage of the sheaf $\caloy$
and let $f_{*}^{-1}(\caloy)$ be the associated simplicial sheaf of
rings on $X_{*}$. As the topological restriction $\calr_{*}|X_{*}$ is
a sheaf of abelian groups on $X_{*}$, we can form its associated
\v{C}ech complex and so we can consider the composed map
$$
\begin{array}{l}
\Hom_{\calr_{*}}(\Omega^1_{\calr_{*}/Y},\tilde
\calm_{*}) \hto
\Der_{Y}(\calr_{*},\tilde\calm_{*})\\
\qquad \hto
\Hom_{f_{*}^{-1}(\caloy)}(\calr_{*},\tilde\calm_{*}) \lto
\Hom_{f^{-1}(\caloy)}(C^\sbullet(\calr_{*}|X_{*})
,C^\sbullet(\tilde\calm_{*})),
\end{array}\leqno (*)
$$
where the first two maps are the natural inclusions and the last one
is given by the \v{C}ech functor. There is always a natural morphism
$\Hom(-,-)\to\RHom(-,-)$, thus taking cohomology we obtain a natural
map
$$
T^k_{X/Y}(\calm)\lto
\Ext^k_{f^{-1}(\caloy)}(C^\sbullet(\calr_{*}|X_{*})
,C^\sbullet(\tilde\calm_{*}))\cong
\Ext^k_{f^{-1}(\caloy)}(\calox,\calm)\,,\leqno (**)
$$
where the final isomorphism results from the fact that the complexes
$C^\sbullet(\calr_{*}|X_{*})$ and $C^\sbullet(\tilde\calm_{*})$ are
quasiisomorphic to $\calox$, resp.\ $\calm$, see \ref{2.16}. Composing
$(**)$ with
$$
Rf_{*}:\Ext^k_{f^{-1}(\caloy)}(\calox,\calm)\lto
\Ext^k_{\caloy}(Rf_{*}(\calox),Rf_{*}(\calm))
$$
gives the desired map. If $\calm=\calox$, then replacing $\tilde\calm$
by $\calr_{*}$ and $C^\sbullet(\tilde\calm_{*})$ by
$C^\sbullet(\calr_{*}|X_{*})$, the first inclusion in $(*)$ becomes
bijective, and so the argument shows the compatibility with the Lie
algebra structure.
\end{proof}

In the next lemma we provide a criterion as to when $Rf_*(\calox)$ is
a perfect complex on $Y_S$. We keep the notation introduced in
\ref{7.16}.

\begin{lem}\label{7.20}
If $Y$ is smooth over $\Sigma$, then $Rf_*(\calox)$ is a
perfect complex on $Y_S$.
\end{lem}

\begin{proof}
For every coherent
$\calos$--module $\caln$ the natural map
$$
Rf_*(\calox)\dotimes\caln \to Rf_*(\caln_X)
$$
is a quasiisomorphism and so the claim follows from \ref{3.4.a}.
\end{proof}

We are now able to apply the constructions of the previous section to
define a semiregularity map for deformations of mappings. Let
$Y\to\Sigma$ and $f:X\to Y_S$ be as in \ref{7.16} and assume moreover
that $Y\to\Sigma$ is smooth. We define a semiregularity map
$$
\tau_\caln: T^2_{X/Y_S}(\caln_X)\to
\prod_{p\ge 0} H^{p+2}(Y_S,\caln\otimes \Omega^p_{Y_S/S})
$$
as the composition of the two maps
$$
T^2_{X/Y_S}(\caln_X)
\stackrel{can}{\lto}
\Ext^2_{Y_S}(Rf_*(\calox),Rf_*(\calox)\dotimes\caln)
\stackrel{\sigma}{\lto}
\prod_{p\ge 0} H^{p+2}(Y_S,\caln\otimes\Omega^p_{Y_S/S})\,,
$$
where the first map is as in \ref{7.19} and $\sigma$ is the
semiregularity map defined in \ref{3.1}.

In analogy with \ref{5.12} we are now able to deduce the
following result.

\begin{theorem}\label{7.22}
Assume that $\pi:Y\to \Sigma$ is proper and smooth and that
$Y_0:=f^{-1}(0)$ is bimeromorphically equivalent to a K\"ahler
manifold. Let $f_0:X_0\to Y_0:=\pi^{-1}(0)$ be a proper holomorphic
map and denote $S=(S,0)$ the basis of the semi\-universal deformation
of $X_0/Y$. If
$$
\tau_0: T^2_{X_0/Y_0}(\calo_{X_0})\xto{\quad}
\prod_{p\ge 0} H^{p+2}(Y_0,\Omega^p_{Y_0})
$$
is the semiregularity map as above, then the following hold.

{\em 1.} $\dim S \ge \dim_\bbbc T^1_{X_0/Y_0}(\calo_{X_0}) -
\dim_\bbbc \ker \tau_0$.

{\em 2.} If $\tau_0$ is injective, then $S$ is smooth at $0$ over a
closed subspace of $\Sigma$.
\end{theorem}

\begin{rems}
1. In the special case that $\Sigma$ is a reduced point, $Y=Y_0$ is a
compact complex manifold of dimension $d$, and $f:X\to Y$ is a map
from a rational curve $X$ into $Y$, the group $T^2(X/Y,\calox)$ is
isomorphic to the group $H^1(X, f^*(\Theta_Y))$, where $\Theta_Y$ is 
the
tangent bundle of $Y$. Thus the top component of the semiregularity
map provides a map
$$
H^1(X, f^*(\Theta_Y))\lto H^d(Y,\Omega^{d-2}_Y)\,,
$$
and we recover thus the map constructed by Behrend and Fantechi.

2. The results \ref{5.3} and \ref{5.4} also admit generalizations to
the case of deformations of mappings. We leave the straightforward
formulation and  proof to the reader.
\end{rems}

\section{Comparison with Bloch's semiregularity map}

\begin{sit}\label{6.1}
Assume that $X$ is a compact complex manifold and $Z\subseteq X$ is a
locally complete intersection of (constant) codimension $q$ with ideal
sheaf $\calj\subseteq\calox$. In this section we will compare our
semiregularity map from \ref{5.9} with the semiregularity map defined
in \cite{Blo}. Observe that for a locally complete intersection
$T^k_{Z/X}(\caloz) \cong H^{k-1}(Z,\caln_{Z/X})$ for all $k$, where
$\caln_{Z/X}\cong (\calj/\calj^2)^\vee :=
\calhom_{\caloz}(\calj/\calj^2,\caloz)$ is the normal bundle of $Z$ in
$X$. Bloch's semiregularity map is constructed as follows. With
$m:=\dim Z$, there is a natural pairing
$$
\Omega_X^{m+1}\times \Lambda^{q-1}\caln_{Z/X}^\vee
\stackrel{1\times \Lambda^{q-1}\bar d}{\lllrarrow}
\Omega_X^{m+1}\times \Omega_X^{q-1}\otimes \caloz
\stackrel{\wedge}{\lto}
\omega_X\otimes\caloz,
$$
where $\bar d:\caln_{Z/X}^\vee\cong \calj/\calj^2\to 
\Omega_X^1\otimes \caloz$
is the map induced by the differential $d:\calj\to
\Omega_X^1$. Equivalently, this amounts to a
map
$$
\Omega_X^{m+1}\lto
\Lambda^{q-1}\caln_{Z/X}\otimes \omega_X
\cong \caln_{Z/X}^\vee  \otimes\omega_Z\,,\leqno (1)
$$
where we have used the adjunction formula $\omega_Z\cong
\det\caln_{Z/X}\otimes\omega_X$ and the isomorphism
$\Lambda^{q-1}\caln_{Z/X}\cong \caln_{Z/X}^\vee
\otimes\det\caln_{Z/X}$.  Dualizing the induced map in
cohomology $ H^{m-k}(X, \Omega_X^{m+1})\to H^{m-k}(Z,
\caln_{Z/X}^\vee\otimes \omega_Z)$ gives a map
$$
\tau_{B}:H^k(Z,\caln_{Z/X})\lto H^{q+k}
(X,\Omega^{q-1}_X), \leqno (2)
$$
and Bloch's semiregularity map is just this map for $k=1$.
We will compare it with the component
$$
\tau:
H^k(Z,\caln_{Z/X})\cong T^{k+1}_{Z/X}(\caloz) \lto H^{q+k}
(X,\Omega^{q-1}_X)
$$
of our semiregularity map defined in \ref{5.9}.
\end{sit}

\begin{prop}\label{6.2}
The maps $\tau_{B}$ and $\tau$ coincide.
\end{prop}

For the proof, we need a more explicit description of the
semiregularity map. First observe that
$$
H^p\left(X, \calh^q_Z(\Omega^{q-1}_X)\right)
\xto{\ \cong\ } H^{q+p}_Z(X,\Omega^{q-1}_X)\ ,\quad p\ge 0\,,
$$
as $\calh^i_Z(\Omega_X^{q-1}) = 0$ for $i \ne q$; for the algebraic
case, see \cite{Blo}, whereas in the analytic case this follows from
\cite{Scheja}. To proceed further we describe the local cohomology
sheaves $\calh^{q}_Z(\Omega_X^{q-1})$ in terms of a Cousin-type
complex.

\begin{sit}\label{6.3}
Let $\cale$ be a vector bundle on $X$ and let $Z$ be as above. Assume
$U\subseteq X$ is an open subset such that the ideal
$\calj\subseteq\calo_U$ of $Z\cap U\subseteq U$ is generated by
sections $f_1, \dots, f_q \in \Gamma (U, \calo_U)$. For an index
$\alpha=(\alpha_1,\ldots,\alpha_p)$ with $1\le\alpha_1<\ldots<
\alpha_p\le q$ set $|\alpha|=p$ and $U_\alpha:=\{x\in U|f_\alpha(x)\ne
0\}$, where $f_\alpha=\prod_{i=1}^p f_{\alpha_i}$. Consider basis
elements $\delta f_\alpha:=\delta f_{\alpha_1}\wedge\dotsb\wedge
\delta f_{\alpha_p}$ and the Cousin complex
$$
\calc_{Z}^\sbullet(\cale|U):\qquad
0\to \cale|U
=:(\cale|U_\emptyset)\delta f_\emptyset\to\prod_{|\alpha|=1}
\cale_\alpha
\delta f_\alpha\to\cdots\to
\prod_{|\alpha|=q} \cale_\alpha \delta f_\alpha \to 0\,,
$$
where $\cale_\alpha:=j_{\alpha*}(\cale|U_\alpha)$ with $j_\alpha$ the
inclusion $U_\alpha\hto U$. The differential on
$\calc_{Z}^\sbullet(\cale|U)$ is given by $\partial(\delta f_\alpha):=
-\sum_{i=1}^{q} \delta f_i\wedge\delta f_\alpha$. The minus sign
occurs in order to have an exact sequence of complexes
$$
0\to \calc^\sbullet(\{U_i\cap U\},\cale)[-1]\to
\calc_{Z}^\sbullet(\cale|U)\to \cale|U\to 0\,,
$$
where $\calc^\sbullet:= \calc^\sbullet(\{U_i\cap U\},\cale)$ is the
sheafified \v{C}ech complex; observe that the differential on
$\calc^\sbullet[-1]$ is the {\em negative} of the differential on
$\calc^\sbullet$! It is well known that $\calh_Z^p(\cale|U)\cong
\calh^p(\calc_{Z}^\sbullet(\cale|U))$. Note that
$\calc_{Z}^\sbullet(\calo_U)$ carries a graded algebra structure via
$\delta f_\alpha \delta f_\beta:=\delta f_\alpha\wedge \delta f_\beta$
and that $\calc_{Z}^\sbullet(\cale|U)$ is a graded module over
$\calc_{Z}^\sbullet(\calo_U)$.
\end{sit}

In \cite[p.61]{Blo}, Bloch defines a natural map
$$
\mu: \caln_{Z/X} \lto \calh^q_Z(X,\Omega^{q-1}_X)
$$
that is locally given as follows.  Multiplying the element
$$
\omega:= \left[\frac{\delta f_1}{f_1} \land \cdots \land
\frac{\delta f_q}{f_q}\right] \in\Gamma\left(U,
\calh^q_Z (\calox)\right)
$$
with $df_1\wedge\dotsb\wedge df_q$ gives a form $\omega \otimes df_1
\land \cdots \land df_q$ in $H^0\left(U,
\calh^q_Z(X,\Omega^q_X)\right)$ that is independent of the choice of
the equations $f_1,\dots,f_q$. The map $\mu$ is then given by
contraction against $\omega \otimes df_1 \land \cdots \land df_q$,
which yields for a linear map $\varphi: \calj \lto \caloz$ the
explicit expression
$$
\mu(\varphi) =  \sum (-1)^{i-1}  \omega\cdot\varphi
(f_i)\otimes df_1 \land \cdots\land
\widehat{df_i} \land \cdots \land df_q\,.\leqno (3)
$$
Now $\tau_{B}$ in \ref{6.1} (2) is the composition
$$
H^k(Z, \caln_{Z/X}) \stackrel{\mu}{\lto} H^k\left(X,
\calh^q_Z(X,\Omega^{q-1}_X)\right) \cong H^{q+k}_Z (X,
\Omega^{q-1}_X) \stackrel{can}{\lto} H^{q+k}(X,
\Omega^{q-1}_X)\,,
$$
see (loc.cit.) for a proof. To compare this map with our
semiregularity map, observe first that $\tau$ also admits a
factorization
$$
H^k(Z,\caln_{Z/X})\cong T^{k+1}_{Z/X} (\caloz)  \lto 
H^{q+k}_Z(X,\Omega^{q-1}_X)
\stackrel{can}{\lto} H^{q+k}(X,\Omega^{q-1}_X)\ ,
$$
see \ref{3.5}. As taking traces is compatible with localization, to
deduce \ref{6.2} it is sufficient to show the following lemma.

\begin{lem} \label{6.4}
The diagram of $\caloz$-modules
$$
\begin{diagram}[midshaft]
\cExt^1_X(\caloz, \caloz) \cong \caln_{Z/X}&& 
\rTo^{\quad*\cdot (-\At(\caloz))^{q-1}/(q-1)!\quad}&&
\cExt^q_X(\caloz, \caloz \otimes
\Omega^{q-1}_X)\\
&\rdTo_{\mu}&&\ldTo_{\Tr}\\
&&\calh^q_Z(X,\Omega^{q-1}_X)
\end{diagram}
%
$$
commutes.
\end{lem}

\begin{proof}
This is a local calculation, so we may suppose $X$ Stein and
$Z$ defined by equations $f_1, \dots, f_q$ so that the Koszul
complex $K_{\sbullet} (\underline f, \calox)$ is an 
$\calox$-resolution of the
sheaf $\caloz$. More explicitly, set
$$
K^{-p} := K_p (\underline f, \calox) =
\bigoplus_{i_1 <\cdots < i_p} \calox \gamma
f_{i_1}\land \cdots \land \gamma f_{i_p}\ ,
$$
with the Koszul differential given by $\partial(\gamma f_j) = f_j$.
Note that $K^{\sbullet}$ is the free graded algebra over $\calox$ with
generators $\gamma f_i\in K^{-1}$, $1\le i\le q$, in (cohomological)
degree $-1$ and that $\partial$ is a derivation of degree 1. In
particular,
$$
\Ext^k_X(\caloz,\caloz\dotimes_\calox\calm )\cong
H^k(\Hom_X(K^\sbullet, K^\sbullet\otimes_\calox\calm)),\quad
k\ge 0\,,
$$
for every coherent $\calox$-module $\calm$.
Consider on $K^{\sbullet}$ the connection
$$
\nabla:K^\sbullet \lto
K^\sbullet\otimes \Omega^1_X\quad\mbox{with}\quad
\nabla\left(\gamma f_{i_1} \land \cdots \land
\gamma f_{i_p}\right) =0\ .
$$
The Atiyah class of $\caloz$ is now the element of
$\Ext_X^1(\caloz,\caloz \otimes\Omega^1_X)$
represented by 
$$
[\partial, \nabla] : K^\sbullet \lto
K^\sbullet\otimes \Omega^1_X\,,
$$
an $\calox$-linear map of degree 1. Note that $\nabla$ as well as
$[\partial,\nabla]$ are derivations on the ring $K^\sbullet$ and that
$[\partial,\nabla](\gamma f_i)=-1\otimes df_i$. Thus, the map
$$
[\partial, \nabla]^{q-1}/(q-1)! : K^\sbullet \lto
K^\sbullet\otimes \Omega^{q-1}_X\,,
$$
which is only nonzero on $K^{-q+1}$ and $K^{-q}$, is given there by
\begin{align*}
&\gamma f_1 \land \cdots \land\widehat{\gamma f_i} \land \cdots\land
\gamma f_q \ \longmapsto\ (-1)^{\binom{q}{2}} 1\otimes (df_1 \land
\cdots \land \widehat{df_i} \land \cdots \land df_q)\\
&\gamma f_1 \land \cdots \land \gamma f_q\  \longmapsto\  \sum_{i}
(-1)^{\binom{q-1}{2}+i}\gamma f_i\otimes df_1 \land \cdots \land
\widehat{df_i} \land \cdots \land df_q\,.
\end{align*}

Consider now $\bar\varphi\in
\Hom_{X}(\calj,\caloz)\cong\Ext^1_X(\caloz,\caloz)$ and write
$\bar\varphi(f_i)=\bar\varphi_i$ with sections $\varphi_i$ of
$\calox$. Under the isomorphism $\Ext^1_X(\caloz,\caloz) \cong
H^1(\Hom_{X}(K^\sbullet,K^\sbullet))$, the element $\bar\varphi$
corresponds to the derivation $\varphi:K^\sbullet\to K^\sbullet$ with
$\varphi(\gamma f_i)=\varphi_i$. The composition $\varphi\circ 
[\partial,
\nabla]^{q-1}/(q-1)!$ is a map of degree $q$ and is therefore 
determined by
the component
$$
\begin{array}{rcl}
\varphi \circ [\partial, \nabla]^{q-1}/(q-1)! :  K^{-q} &
\longrightarrow &K^0 \otimes \Omega^{q-1}_X\\[3pt]
\gamma f_1 \land \cdots \land \gamma f_q & \longmapsto &
\sum(-1)^{\binom{q-1}{2}+i} \varphi_i\otimes
df_1 \land \cdots \land \widehat{df_i} \land \cdots
\land df_q\ .
\end{array}
$$
Note that $\varphi\circ [\partial, \nabla]^{q-1}$ represents
$\langle [\bar\varphi], \At^{q-1}(\caloz)\rangle$.

Now we wish to take traces to produce an element in
$\calh^q_Z(X,\Omega^{q-1}_X)$. Let
$\calc_{Z}^\sbullet=\calc_{Z}^\sbullet(\calox)$ be the Cousin complex,
see \ref{6.3}. Let $\hat\gamma f_i$ be the basis of $\Hom_{X}(K^{-1},
\calox)$ dual to $\gamma f_i$ and set
$$
\hat\gamma f_\alpha:= \hat \gamma
f_{\alpha_1}\wedge\dotsb\wedge
\hat\gamma f_{\alpha_p}= (-1)^{\binom{p}{2}}\hat \gamma
f_{\alpha_p}\wedge\dotsb\wedge
\hat\gamma f_{\alpha_1}.
$$
With this convention, $(-1)^{\binom{p}{2}}\hat\gamma
f_\alpha$ is the basis element dual to $\gamma f_\alpha:=\gamma
f_{\alpha_1}\wedge \dotsb\wedge\gamma f_{\alpha_p}$ and the
differential on $\cHom_{X}(K^\sbullet,\calox)$ is given by multiplying
from the left with $-\sum f_i\cdot \hat\gamma f_i$. The computation
above gives that
$$
\varphi\circ[\partial, \nabla]^{q-1}/(q-1)!=
\hat\gamma f_1 \land \cdots \land \hat\gamma f_q
\sum(-1)^{q+i} \varphi_i\otimes  df_1 \land \cdots \land
\widehat{df_i} \land \cdots \land df_q\,.
$$
Using \ref{6.5} below, this is mapped under the trace map to
$$
(-1)^{q-1}\frac{\delta f_1 \wedge \cdots \wedge\delta f_q}{f_1
\cdots f_q}
\cdot\sum (-1)^{i-1} \varphi_i\otimes df_1 \land
\cdots \land \widehat{df_i} \land \cdots \land
df_q\,.
$$
Comparing with $(3)$ above, the result follows.
\end{proof}

It remains to verify the following lemma.

\begin{lem}\label{6.5}
The trace map $\cExt_X^q(\caloz, \caloz) \lto \calh^q_Z (\calox)$
maps the class of $\hat\gamma f_1 \land \cdots \land \hat\gamma
f_q\in\Hom_{X}(K^{-q},\calox)\subseteq\Hom_X(K^\sbullet,K^\sbullet)$ 
onto
the class of
$$
\frac{\delta f_1\land \dotsb\land\delta f_q}{f_1\cdots f_q}\,.
$$
\end{lem}

\begin{proof}
As $\calc_{Z}^\sbullet$ is a sheaf of flat $\calox$--modules, the
complex $K^\sbullet \otimes \calc_{Z}^\sbullet$ is quasiisomor\-phic
to $\caloz\otimes\calc_{Z}^\sbullet\cong \caloz$. Therefore, the
canonical projection $\calc_{Z}^\sbullet\to\calc_{Z}^0=\calox$ induces
a quasiisomorphism
$$
K^{\sbullet} \otimes \calc_{Z}^\sbullet  \lto  K^{\sbullet}\,.
$$
In a first step, we construct a section of this projection,
$$
\begin{diagram}[small]
    K^{\sbullet}&&
    \rTo^\psi&&
     K^{\sbullet} \otimes \calc_{Z}^\sbullet\\
    &\rdTo_\id
    && 
    \ldTo_{proj}\\
    && K^{\sbullet}
\end{diagram}
$$
as follows. The map $\hat\gamma f_\alpha \mapsto \delta f_\alpha
/f_\alpha$, where $\hat\gamma f_\alpha = \hat\gamma f_{\alpha_1} \land
\cdots\land \hat\gamma f_{\alpha_p}$ is as above, realizes
$\cHom_{X}(K^\sbullet,\calox)$ as a subcomplex of
$\calc_{Z}^\sbullet$. We define $\psi$ to be the composition
$$
K^\sbullet \lto \cHom_X(K^\sbullet,K^\sbullet)\cong
K^\sbullet\otimes \cHom_X(K^\sbullet,\calox)
\hto K^\sbullet\otimes\calc_{Z}^\sbullet\,,
$$
where the first map is given by $k\mapsto k\cdot\id$. As the identity
corresponds to the element $\sum(-1)^{\binom{|\alpha|}{2}}\gamma
f_\alpha \otimes \hat\gamma f_\alpha$ in $K^\sbullet\otimes
\cHom_X(K^\sbullet,\calox)$, the map $\psi$ is given explicitly by
$$
k\longmapsto k\cdot \sum_{\alpha} (-1)^{\binom{|\alpha|}{2}}\gamma
f_{\alpha} \otimes \delta f_{\alpha}/f_\alpha\,.
$$
Now we can define the local trace map
$\cHom_{X}(K^\sbullet,K^\sbullet)\to \calc_{Z}^\sbullet$ as the
composition of
$$
\cHom_{X}(K^\sbullet,K^\sbullet)
\cong
\cHom_{X}(K^\sbullet,\calox)\otimes K^\sbullet
\stackrel{1\otimes \psi}{\lrarrow}
\cHom_{X}(K^\sbullet,\calox)\otimes K^\sbullet
\otimes\calc_{Z}^\sbullet
\stackrel{\Tr\otimes 1}{\lrarrow}
\calc_{Z}^\sbullet.
$$
The image of $\hat\gamma f_1\land\dotsb\land \hat\gamma f_q $ under
these maps is just given by
\begin{eqnarray*}
\hat\gamma f_1\land\dotsb\land \hat\gamma f_q
&\stackrel{1\otimes\psi}{\longmapsto}&
\hat\gamma f_1\land\dotsb\land \hat\gamma f_q\otimes
\sum_\alpha(-1)^{\binom{|\alpha|}{2}}
\gamma f_\alpha\otimes \delta f_\alpha /f_\alpha\\
& \stackrel{\Tr\otimes 1}{\longmapsto} &
\frac{\delta f_1\land\dotsb\land \delta f_q}{f_1\cdots
f_q}\,,
\end{eqnarray*}
as desired.
\end{proof}

\begin{sit}
Another application of the above construction concerns the
infinitesimal Abel-Jacobi map. Let $X$ be an $n$--dimensional compact
algebraic manifold and $Z\subseteq X$ a closed submanifold of
codimension $q$. The infinitesimal Abel-Jacobi map is the differential
at $[Z]$ of the Abel-Jacobi map
$$
\Hilb_X \lto  J^q(X)
$$
into the intermediate Jacobian $J^q(X)$; see, for example,
\cite{Gre,Gri2}, or \cite{Voi} in the analytic case. This differential
can be considered as a map
$$
\beta: H^0 (Z,\caln_{Z/X}) \lto  H^q(X, \Omega^{q-1}_X)
$$
that has the same homological description in terms of Serre duality as
Bloch's semiregularity map and is just the map $\tau_{B}$ in \ref{6.1}
(2) for $k=0$, see (loc.cit.).
\end{sit}

Applying \ref{6.2}, we obtain the following description of the
infinitesimal Abel-Jacobi map.

\begin{prop} \label{6.6}
The infinitesimal Abel-Jacobi map fits into the commutative diagram
\begin{diagram}[h=7mm,midshaft]
H^0(Z,\caln_{Z/X})
&\rTo^{\cong}
& \Ext^1_X (\caloz, \caloz)\\
\dTo<\beta && \dTo>{\langle *, (-\At(\caloz))^{q-1}/(q-1)!
\rangle}\\
H^q(X, \Omega^{q-1}_X) & \lTo^\Tr & \Ext^q_{X} (\caloz,
\caloz \otimes \Omega^{q-1}_X).
\end{diagram}
\qed
\end{prop}

This statement should generalize. There should be an Abel-Jacobi map
for deformations of arbitrary coherent sheaves on a compact algebraic
manifold such that its differential is essentially given by
multiplication with powers of the Atiyah class as above. More
precisely, let us pose the following problem.

\begin{prob}\label{6.7}
For a coherent sheaf $\calf_0$ on a compact algebraic manifold $X$ the
Chern character $\ch_k(\calf_0)$ is a well defined class in the Chow
group $CH^k(X)_\bbbq$. Assume that the sheaf $\calf$ on $X\times S$ is
a semiuniversal deformation of $\calf_0$ over a germ $S=(S,0)$. If
$\calf_s$ denotes the restriction of $\calf$ to the fibre $X\cong
X\times\{s\}$, then the following should hold:
\begin{enumerate}
\item The map $s \longmapsto \ch_k(\calf_s)$ provides a family of
$k$-dimensional cycles on $X$. 

\item Integrating over a (topological) $(2k+1)$-chain in $X$ with
boundary $\ch_k(\calf_s) -\ch_k(\calf_0)$ gives a holomorphic map
$S\stackrel{\varphi_k}{\lto} J^k(X)$. 

\item The derivative of $\varphi_k$ is given by
\begin{diagram}
T_0 S &\rTo^{\cong\quad} & \Ext^1_X(\calf_0, \calf_0) &
\rTo^{\Tr\langle *, (- \At(\calf_0))^k/k!\rangle}
 & H^{k+1}(X,\Omega^k_X)\,.
\end{diagram}
\end{enumerate}
\end{prob}

\section{Appendix: Infinitesimal deformations and integral
dependence}

We recall first the definition of integral dependence, see
\cite{ZSa1}. Let $R$ be a ring and $I \subseteq R$ an ideal. An
element $x \in R$ is {\em integral\/} over $I$ if there is an equation
$$
x^n + a_{1} x^{n-1} + \dotsb + a_n = 0
$$
with $a_{\nu}\in I^\nu$.

For instance, every element of $I$ is integral over $I$. The set
${\overline I}\subseteq R$ of all elements from $R$ that are integral
over $I$ is an ideal, the {\it integral closure\/} of $I$ in $R$.

We remind the reader of the following criterion for integral
dependence that we formulate for our purposes as follows. Let $k$ be
an algebraically closed field and let $A=(A,\fm)$ be a local
noetherian complete $k$-algebra with residue field $k$. If $I\subseteq
\fm$ is an ideal, then $f$ is in the integral closure of $I$ if and
only if for every ``arc'' $\alpha: A\to k[\![T]\!]$ the element
$\alpha(f)$ is contained in $\alpha(I) k[\![T]\!]$.

\begin{sit}\label{9.1}
For the main result of this section we consider the following setup.
Let $\Lambda\to A$ be a morphism of complete noetherian $k$-algebras
with residue field $k$ and assume that $A\cong R/I$ with
$R:=\Lambda[\![X_1,\ldots,X_s]\!]$ and $I\subseteq \fm_\Lambda
R+(X_1,\ldots,X_s)^2$. It is well known that for every finite
$A$--module $M$ there is a canonical isomorphism
$$
T^1_{A/\Lambda}(M)\cong \coker\left(
\Hom_{A}(\Omega^1_{A/\Lambda},M)\lto \Hom_{A}(I,M)\right),\leqno (*)
$$
where $T^1_{A/\Lambda}(M)$ is the first tangent cohomology. The
elements of $T^1_{A/\Lambda}(M)$ correspond to isomorphism classes of
algebra extensions $[A']$ of $A$ by $M$. On the level of such algebra
extensions the isomorphism above is given as follows. If $[A']\in
T^1_{A/\Lambda}$ is an algebra extension of $A$ by $M$, let $p':R\to
A'$ be a morphism of $\Lambda$--algebras lifting the given map $p:R\to
A$. Restricting $p'$ to $I$ gives a map $\varphi_{A'}:= p'|I:I\to M$,
and the correspondence $(*)$ assigns to $[A']$ the residue class of
this homomorphism. As $I\subseteq \fm_\Lambda R+(X_1,\ldots,X_s)^2$, 
we
get in particular that
$$
T^1_{A/\Lambda}(k)\cong  \Hom_{A}(I,k)\,.\leqno (**)
$$
There is always a canonical inclusion of
$\Ext^1_A(\Omega^1_{A/\Lambda}, M)$ into $T^1_{A/\Lambda}(M)$. In case
$M=k$, another important subspace of $T^1_{A/\Lambda}(k)$ is
$\Ex^c_{A/\Lambda}(k)$, the space of {\em curvi\-linear extensions\/}.
This is by definition the subspace generated by all curvilinear
extensions $[A']$, which are those extensions that fit into a
commutative diagram of $\Lambda$--algebras
$$
\label{diagram}
\begin{diagram}[h=7mm]
0&\rTo &  k &\rTo & A' & \rTo & A&\rTo &0\hphantom{\,.}\\
&& \dTo<\cong &&\dTo &&\dTo \\
0&\rTo &  k&\rTo &k[\![t]\!]/(t^{n+1}) & \rTo
&k[\![t]\!]/(t^n) &\rTo &0\,.
\end{diagram}\leqno
(\mathbf D)
$$
We will give the following characterizations of these
subspaces.
\end{sit}

\begin{theorem}
    \label{1.1}
{\em 1.}  If $k$ is algebraically closed then under the
isomorphism $(**)$ the subspace $\Ex^c_{A/\Lambda}(k)$ of
$T^1_{A/\Lambda}(k)$ corresponds to the subspace
$\Hom_{A}(I/J,k)$ of $\Hom_{A}(I,k)$, where $J$ is the integral
closure of $\fm I$ in $I$.

{\em 2.} The elements of $\Ext^1_A(\Omega^1_{A/\Lambda}, M)\hto
T^1_{A/\Lambda}(M)$ correspond to those extensions $[A']$ of $[A]$ by
$M$ for which the associated Jacobi map
$j:M\to\Omega^1_{A'/\Lambda}\otimes_{A'}A$ is injective. 

{\em 3.} If {\em Char} $k=0$ and $\Lambda=M=k$, then there are
inclusions
$$
\Ex^c_{A/k}( k)\subseteq
\Ext^1_A(\Omega^1_{A/k}, k)\subseteq T^1_{A/k}(k)\,.
$$
\end{theorem}

\begin{proof}
For (1), let $A'$ be a curvilinear extension of $A$ by $k$ and let
$p':R\to A'$ be a morphism of $\Lambda$--algebras lifting the given
map $p:R\to A$, so that $\varphi_{A'}=p'|I:I\to k$ corresponds under
$(**)$ to the extension $[A']$. By the valuative criterion of integral
dependence mentioned above, $J$ is in the kernel of $p'$, whence
$\varphi_{A'}\in \Hom_{A} (I/J,k)$. Thus
$\Ex^c_{A/\Lambda}(k)\subseteq \Hom_{A}(I/J,k)$. To show equality,
assume that $\alpha$ is a $k$-linear form on $\Hom_{A}(I,k)$ that
vanishes on $\Ex^c_{A/\Lambda}(k)$. We need to show that $\alpha$
vanishes on $\Hom_{A}(I/J,k)$. Such a linear form can be written as
$\Hom_{A}(I,k)\ni f\mapsto \alpha(f)= f(x)$ for some $x\in I$. By
assumption, for every curvilinear extension $A'$ of $A$ by $k$, the
element $\alpha(\varphi_{A'})=\varphi_{A'}(x)$ vanishes. Applying the
valuative criterion of integral dependence it follows that $x\in J$
and so $\alpha$ vanishes on $\Hom_{A}(I/J,k)$, as desired.

In order to deduce (2) note that the map
$\Ext^1_A(\Omega^1_{A/\Lambda}, M)\to T^1_{A/\Lambda}(M)$ assigns to
an extension $0\to M\to E \stackrel{q}{\to}\Omega^1_{A/\Lambda}\to 0$
the algebra extension $[A']$ of $A$ by $M$ that is the quotient of the
trivial extension $A[E]$ by the ideal $\ker ((d, -q):A[E]\to
\Omega^1_{A/\Lambda})$, where $d$ is the differential. The reader may
easily verify that then $E\cong \Omega^1_{A'/\Lambda}\otimes_{A'}A$
and that the Jacobi map $j$ becomes the inclusion of $M$ into $E$,
whence $j$ is injective. Conversely, if for an extension $[A']$ the
map $j$ is injective, then it is easily seen that the construction
just described recovers $[A']$ from the extension
$E:=\Omega^1_{A'/\Lambda}\otimes_{A'}A$ of $\Omega^1_{A/\Lambda}$ by
$M$.

Finally for (3), if $[A']$ is a curvilinear extension as in the
diagram $(\mathbf D)$ in \ref{diagram}, then the composed
map $k\stackrel{j}{\to} \Omega^1_{A'/\Lambda}\otimes_{A'}A
\to k[\![t]\!]/(t^n)\cdot dt$ is the map $1\mapsto d(t^n)$ and
so is injective.  Hence $j$ is also injective, proving the
inclusion $\Ex^c_{A/k}( k)\subseteq \Ext^1_A(\Omega^1_{A/k},
k)$.
\end{proof}

The following result shows how to bound the dimension of $A$ in terms 
of its curvilinear extensions, as
$$
\Ex^c_{A/\Lambda}(k) \cong \Hom_{A}(I/J,k) \cong \Hom_{k}(I/(J+\fm 
I),k)
$$
by the preceding result. Kawamata \cite{Kaw2} attributes the
corresponding geometric argument to Mori.

\begin{prop}\label{1.2}
Let $A=R/I$ be a quotient of a local ring
$(R,\fm,k)$ modulo an ideal $I\subseteq \fm$. If $J\subseteq I$
is integral over $\fm I$, then
$$
\dim A\ge \dim R - \dim_k(I/(J+\fm I))\,.
$$
\end{prop}

\begin{proof}
Replacing $J$ by $J  + \fm I$ we may assume that $J
\supseteq \fm I$.  Choose elements $x_1, \dots, x_s \in I$ that
form a basis of the $k$-vector space $I/J$ and consider the
natural ring  homomorphism
$$
k[X_1, \dots, X_s]\xrightarrow{\quad}
\bigoplus\limits^\infty_{\nu =0}
\left(I^\nu/\fm I^\nu\right)T^\nu = R[IT]/\fm R[IT]
$$
given by $X_ i \mapsto \bar x_iT \in (I/\fm I)T$. In a first step we
prove that this map is finite. In fact, the elements $\bar f T\,, f
\in J,$ generate the ring $R[IT]/\fm R[IT]$ as an algebra over $k[X_1,
\dots, X_s]$, and if $f^n + a_{1} f^{n-1} + \cdots + a_n = 0$ is an
equation of integral dependence for such an $f\in J$ over $\fm I$, the
coefficients satisfy $a_\nu \in (\fm I)^{\nu}$, whence $({\bar f} T)^n
= 0$ and finiteness follows. This implies
$$
\dim R[IT]/\fm R[IT] \le s\,.\leqno (1)
$$
As $R[IT]/\fm R[IT]$ appears as the special fibre of
$$ R/I\xrightarrow{\quad} G_I(R) := \bigoplus\limits^\infty_{\nu
=0} I^\nu  /I^{\nu +1}\,,
$$
we obtain
$$
\dim R[IT]/\fm R[IT] \ge \dim G_I(R) - \dim R/I = \dim R-\dim
R/I\,,
$$
see \cite[Thms.15.1, 15.7]{Mat}. Together with (1) the result
follows.
\end{proof}

For the next result we have to assume that the ground field $k$
has characteristic zero.

\begin{prop}\label{1.3}
If $\Lambda\to A$ is a morphism of local noetherian complete
$k$--algebras with residue fields $k$, then
$$
\dim A\ge \dim_k \Hom_A(\Omega^1_{A/\Lambda}, k) - \dim_k
\Ext^1_A(\Omega^1_{A/\Lambda}, k)\,.
$$
\end{prop}

\begin{proof}
In the absolute case, where $\Lambda =k$, this is a result due to
Scheja and Storch, see \cite[3.5]{SSt}. Alternatively, it follows from
the chain of inequalities
$$
\begin{array}{rcl}
\dim A &\ge& \dim R-\dim_k (I/J+\fm I)=\dim_k
\Hom_A(\Omega^1_{A/k}, k)-\Ex_{A/k}^c(k)\\
&\ge&
\dim_k\Hom_A(\Omega^1_{A/k}, k)-\Ext^1_A(\Omega^1_{A/k}, k),
\end{array}
$$
where we have applied \ref{1.2} and \ref{1.1}.

To deduce the general case, set $\bar A:=A/\fm_\Lambda A$. The
spectral sequence
$$
E^{pq}_2=\Ext^p_{\bar A}(\Tor_q^A(\Omega^1_{A/\Lambda}, \bar A),
k)\Rightarrow \Ext^{p+q}_A(\Omega^1_{A/\Lambda},k)
$$
yields
$$
\Hom_A(\Omega^1_{A/\Lambda},k)\cong\Hom_{\bar A}(\Omega^1_{\bar
A/k},k)\quad\text{and}\quad
\Ext^1_A(\Omega^1_{A/\Lambda},k)\supseteq\Ext^1_{\bar
A}(\Omega^1_{\bar A/k},k)\,.
$$
Hence the result follows from the chain of inequalities
$$
\begin{array}{rcl}
\dim A\ge \dim\bar A &\ge &\dim_k\Hom_{\bar A}(\Omega^1_{\bar
A/k},k)-\dim_k \Ext^1_{\bar A}(\Omega^1_{\bar A/k},k)\\
&\ge &\dim_k\Hom_A(\Omega^1_{A/\Lambda},k)-
\dim_k\Ext^1_A(\Omega^1_{A/\Lambda},k).
\end{array}
$$
\end{proof}


\end{document}